# Brownian excursion area, Wright's constants in graph enumeration, and other Brownian areas *

**Svante Janson**

*Department of Mathematics, Uppsala University*
*PO Box 480*
*SE-751 06 Uppsala, Sweden*
*e-mail:* svante.janson@math.uu.se
*url:* http://www.math.uu.se/∼svante/

**Abstract:** This survey is a collection of various results and formulas by different authors on the areas (integrals) of five related processes, viz. Brownian motion, bridge, excursion, meander and double meander; for the Brownian motion and bridge, which take both positive and negative values, we consider both the integral of the absolute value and the integral of the positive (or negative) part. This gives us seven related positive random variables, for which we study, in particular, formulas for moments and Laplace transforms; we also give (in many cases) series representations and asymptotics for density functions and distribution functions. We further study Wright's constants arising in the asymptotic enumeration of connected graphs; these are known to be closely connected to the moments of the Brownian excursion area.

The main purpose is to compare the results for these seven Brownian areas by stating the results in parallel forms; thus emphasizing both the similarities and the differences. A recurring theme is the Airy function which appears in slightly different ways in formulas for all seven random variables. We further want to give explicit relations between the many different similar notations and definitions that have been used by various authors. There are also some new results, mainly to fill in gaps left in the literature. Some short proofs are given, but most proofs are omitted and the reader is instead referred to the original sources.

## Contents

*This is an original survey paper

## 1. Introduction

This survey started as an attempt to better understand two different results, viz. the asymptotic enumeration of connected graphs with a given number of independent cycles by Wright [58; 59] and the formulas for moments of the Brownian excursion area by Louchard [34], Takács [44] and others, and the surprising connection between these two seemingly unrelated results found (and explained) much later by Spencer [43], who showed that the same sequence of constants appears in both results (see (36) below). The literature may seem confusing (to me at least), however, because different authors, and sometimes the same author in different papers, have used not only different notations for the same constants but also several different related constants. I therefore set out to collect the various definitions and notations, and to list the relations between them explicitly. I was further inspired by the paper by Flajolet and Louchard [17], studying further related properties of the distribution.

After doing this, I realized that many of the results for Brownian excursion area have close parallels for the integrals of (the absolute value of) other related processes, viz. Brownian bridge, Brownian motion, Brownian meander and Brownian double meander, and also for the positive parts of Brownian motion and Brownian bridge (see Sections 20–25). In particular, there are similar recursion formulas for moments and similar expressions for double Laplace transforms involving the Airy function, but the details differ between the seven different Brownian areas. I therefore decided to collect various results for these seven different areas of Brownian processes together so that both the similarities and the differences would become clear. One important source of results and inspiration for this is the series of papers by Takács [44; 45; 46; 47; 48; 49; 50]; another, treating different cases in a unified way, is Perman and Wellner [38]. In the course of doing this, I have also added some minor new results motivated by parallel results. For example, much of Section 29 and the related Appendix B, which both generalize results by Flajolet and Louchard [17], seem to be new, as do many of the explicit asymptotic results as $x \to 0$.

To faciliate comparisons and further studies in the literature, I have often listed the notations used by other authors to provide a kind of dictionary, and I have tried to give as many explicit relations as possible, at the cost of sometimes perhaps being boring; I hope that the reader that gets bored by these details can skip them. I have found draft versions of this survey very useful myself, and I hope that others will find it useful too.

There are two main methods to prove the results on moments for the various processes given here. The first studies Brownian motion using some form of the Feynman–Kac formula, see for example Kac [29; 30; 31], Louchard [34], Takács [49] and Perman and Wellner [38]; this leads to various expressions involving the Airy function. See also Section 27 for a probabilistic explanation of the two slightly different forms of formulas that arise for the different cases. A related method is the path integral method used in mathematical physics and applied to Brownian areas by Majumdar and Comtet [35]. The second main method is combinatorial and studies the corresponding discrete structures, combined with singularity analysis for the resulting generating functions; see for example Takács [44; 45; 46; 47; 48; 50]; and Nguyên Thê [36].

We have omitted most of the proofs, but we include or sketch some. In particular, we give for completeness in Appendix C proofs by the Feynman–Kac formula of the basic analytic identities (a double Laplace transform that is expressed using the Airy function) for all types of areas of Brownian processes studied here.

There are, of course, many related results not covered here. For example, we will not mention integrals of the supremum process or of local time [26], or processes conditioned on their local time as in [8]. See also the impressive collection of formulas for Brownian functionals by Borodin and Salminen [7].

The Brownian excursion area is introduced in Section 2 and some equivalent descriptions are given in Section 3. The results on graph enumeration are given in Section 4, the connection with the Brownian excursion area is given in Section 5, and various aspects of these results and other related results on the Brownian excursion area are discussed in Sections 6–18. Integrals of higher powers of a Brownian excursion are discussed briefly in Section 19. The other types of Brownian processes are studied in Sections 20–25; these sections are intended to be parallel to each other as much as possible. Furthermore, in Sections 24 and 25, results are given both for the positive parts of a Brownian bridge and motion, respectively, and for the joint distribution of the positive and negative parts. Sections 26 and 27 give more information on the relation between the seven different Brownian areas (the reader may prefer to read these sections before reading about the individual variables), and Section 28 compares them numerically. Finally, Section 29 discusses negative moments of the Brownian areas.

The appendices contain some general results used in the main part of the text.

The first few moments of the seven different Brownian areas studied here are given explicitly in 11 tables throughout the paper together with the first elements in various related sequences of constants. Most of these values, and some other expressions in this paper, have been calculated using `Maple`.

We let $B(t)$ be a standard Brownian motion on the $[0,\infty)$ with $B(0) = 0$ (usually we consider only the interval $[0,1]$). We remind the reader that a Brownian bridge $B_{\mathrm{br}}(t)$ can be defined as $B(t)$ on the interval $[0,1]$ conditioned on $B(1)=0$; a (normalized) Brownian meander $B_{\mathrm{me}}(t)$ as $B(t)$ on the interval $[0,1]$ conditioned on $B(t)\geq 0$, $t\in[0,1]$, and a (normalized) Brownian excursion $B_{\mathrm{ex}}(t)$ as $B(t)$ on the interval $[0,1]$ conditioned on both $B(1)=0$ and $B(t)\geq 0$, $t\in[0,1]$. These definitions have to be interpreted with some care, since we condition on events of probability 0, but it is well-known that they can be rigorously defined as suitable limits; see for example Durrett, Iglehart and Miller [14]. There are also many other constructions of and relations between these processes. For example, see Revuz and Yor [39, Chapter XII], if $T>0$ is any fixed time (for example $T=1$), and $g_T := \max\{t<T : B(t)=0\}$ and $d_T := \min\{t>T : B(t)=0\}$ are the zeros of $B$ nearest before and after $T$, then conditioned on $g_T$ and $d_T$, the restrictions of $B$ to the intervals $[0,g_T]$, $[g_T,T]$ and $[g_T,d_T]$ are, respectively, a Brownian bridge, a Brownian meander and a Brownian excursion on these intervals; the normalized processes on $[0,1]$ studied in this paper can be obtained by standard Brownian rescaling as $g_T^{-1/2}B(tg_T)$, $(T-g_T)^{-1/2}B(g_T+t(T-g_T))$ and $(d_T-g_T)^{-1/2}B(g_T+t(d_T-g_T))$. (See also Section 27, where we consider Brownian motion, bridges, meanders and excursions on other intervals than $[0,1]$, and obtain important result by this result and other rescalings.) Another well-known construction of the Brownian bridge is $B_{\mathrm{br}}(t)=B(t)-tB(1)$, $t\in[0,1]$. Further, the (normalized) Brownian excusion is a Bessel(3) bridge, i.e., the absolute value of a three-dimensional Brownian bridge, see Revuz and Yor [39].

While the Brownian motion, bridge, meander and excursion have been extensively studied, the final process considered here, viz. the Brownian double meander, has not been studied much; Majumdar and Comtet [35] is the main exception. In fact, the name has been invented for this paper since we have not seen any given to it previously. We define the Brownian double meander by $B_{\mathrm{dm}}(t) := B(t) - \min_{0\le u\le 1} B(u)$; this is a non-negative continuous stochastic process on $[0,1]$ that a.s. is 0 at a unique point $\tau \in [0,1]$, and it can be regarded as two Brownian meanders on the intervals $[0,\tau]$ and $[\tau,1]$ joined back to back (with the first one reversed), see Section 23.

## 2. Brownian excursion area and its moments

Let $B_{\mathrm{ex}}$ denote a (normalized) Brownian excursion and

$$\mathcal{B}_{\mathrm{ex}} := \int_0^1 B_{\mathrm{ex}}(t)\,\mathrm{d}t, \tag{1}$$

the Brownian excursion area. Two variants of this are

$$\mathcal{A} := 2^{3/2}\mathcal{B}_{\mathrm{ex}} \tag{2}$$

studied by Flajolet, Poblete and Viola [18] and Flajolet and Louchard [17] (their definition is actually by the moments in (8) below), and

$$\xi := 2\mathcal{B}_{\mathrm{ex}} \tag{3}$$

used in [21; 23; 22]. (Louchard [34] uses $\xi$ for our $\mathcal{B}_{\mathrm{ex}}$; Takács [44; 45; 46; 47; 48] uses $\omega^+$; Perman and Wellner [38] use $A_{\mathrm{excur}}$; Spencer [43] uses $L$. Flajolet and Louchard [17] use simply $\mathcal{B}$ for $\mathcal{B}_{\mathrm{ex}}$, but since we in later sections will consider areas of other related processes too, we prefer to use $\mathcal{B}_{\mathrm{ex}}$ throughout the paper for consistency.) Flajolet and Louchard [17] call the distribution of $\mathcal{A}$ the *Airy distribution*, because of the relations with the Airy function given later. (But note that the other Brownian areas in this paper also have similar connections to the Airy function.)

This survey centres upon formulas for the moments of $\mathcal{B}_{\mathrm{ex}}$ and various analogues of them. The first formula for the moments of $\mathcal{B}_{\mathrm{ex}}$ was given by Louchard [34] (using formulas in [33]), who showed (using $\beta_n$ for $\mathbb{E}\,\mathcal{B}_{\mathrm{ex}}^n$)

$$\mathbb{E}\,\mathcal{B}_{\mathrm{ex}}^k = (36\sqrt{2})^{-k}\frac{2\sqrt{\pi}}{\Gamma((3k-1)/2)}\gamma_k, \qquad k\ge 0, \tag{4}$$

where $\gamma_k$ satisfies

$$\gamma_r = \frac{12r}{6r-1}\,\frac{\Gamma(3r+1/2)}{\Gamma(r+1/2)} - \sum_{j=1}^{r-1}\binom{r}{j}\frac{\Gamma(3j+1/2)}{\Gamma(j+1/2)}\gamma_{r-j}, \qquad r\ge 1. \tag{5}$$

Takács [44; 45; 46; 47; 48] gives the formula (using $M_k$ for $\mathbb{E}\,\mathcal{B}_{\mathrm{ex}}^k$)

$$\mathbb{E}\,\mathcal{B}_{\mathrm{ex}}^k = \frac{4\sqrt{\pi}\,2^{-k/2}k!}{\Gamma((3k-1)/2)}K_k, \qquad k\ge 0, \tag{6}$$

where $K_0 = -1/2$ and

$$K_k = \frac{3k-4}{4} K_{k-1} + \sum_{j=1}^{k-1} K_j K_{k-j}, \qquad k \geq 1. \tag{7}$$

(See also Nguyên Thê [36], where $\mathbb{E}\,\mathcal{B}_{\mathrm{ex}}^r$ is denoted $M_r^D$ and $a_r^D$ equals $2^{5/2-r/2}K_r$.)

We will later, in Section 16, see that both the linear recurrence (4)–(5) and the quadratic recurrence (6)–(7) follow from the same asymptotic series (98) involving Airy functions. We will in the sequel see several variations and analogues for other Brownian areas of both these recursions.

Flajolet, Poblete and Viola [18] and Flajolet and Louchard [17] give the formula

$$\mathbb{E}\,\mathcal{A}^k = \frac{2\sqrt{\pi}}{\Gamma((3k-1)/2)} \Omega_k, \qquad k \geq 0, \tag{8}$$

(and use it as a definition of the distribution of $\mathcal{A}$; they further use $\mu^{(k)}$ and $\mu_k$, respectively, for $\mathbb{E}\,\mathcal{A}^k$). Here $\Omega_k$ is defined by Flajolet and Louchard [17] by $\Omega_0 := -1$ and the recursion

$$2\Omega_k = (3k-4)k\Omega_{k-1} + \sum_{j=1}^{k-1} \binom{k}{j} \Omega_j \Omega_{k-j}, \qquad k \geq 1. \tag{9}$$

In particular, $\Omega_1 = 1/2$.

It is easily seen that (8)–(9) are equivalent to (6)–(7) and the relation

$$K_k = \frac{\Omega_k}{2^{k+1} k!}, \qquad k \geq 0. \tag{10}$$

Flajolet, Poblete and Viola [18] define further

$$\omega_k := \frac{\Omega_k}{k!} = 2^{k+1} K_k, \qquad k \geq 0, \tag{11}$$

$$\omega_k^* := 2^{2k-1}\omega_k = \frac{2^{2k-1}}{k!}\Omega_k = 2^{3k} K_k, \qquad k \geq 0. \tag{12}$$

By (12), the recursion (9) translates to

$$\omega_k^* = 2(3k-4)\omega_{k-1}^* + \sum_{j=1}^{k-1} \omega_j^* \omega_{k-j}^*, \qquad k \geq 1; \tag{13}$$

with $\omega_0^* := -1/2$ and $\omega_1^* = 1$. Note that $\omega_k^*$ thus is an integer for $k \geq 1$. The numbers $\omega_k^*$ are the same as $\omega_{k0}^*$ in Janson [21]. (This is a special case of $\omega_{kl}^*$ in [21], but we will only need the case $l = 0$.) The sequence $(\omega_k^*)$ is called the Wright–Louchard–Takács sequence in [18].

By a special case of Janson [21, Theorem 3.3],

$$\mathbb{E}\,\xi^k = \frac{2^{2-5k/2}\sqrt{\pi}\,k!}{\Gamma((3k-1)/2)} \omega_k^*, \qquad k \geq 0, \tag{14}$$

which is equivalent to (6) by (3) and (12).

Further, (4) is equivalent to (8) and

$$\gamma_k = 18^k \Omega_k, \qquad k \geq 0. \tag{15}$$

It is easily seen from (9) that $2^k \Omega_k$ is an integer for all $k \geq 0$, and thus by (15), $\gamma_k$ is an integer for $k \geq 0$.

## 3. Other Brownian representations

Let $B_{\mathrm{br}}$ denote a Brownian bridge. By Vervaat [54], $B_{\mathrm{br}}(\cdot) - \min_{[0,1]} B_{\mathrm{br}}$ has the same distribution as a random translation of $B_{\mathrm{ex}}$, regarding $[0,1]$ as a circle. (That is, the translation by $u$ is defined as $B_{\mathrm{ex}}(\lfloor \cdot - u \rfloor)$.) In particular, these random functions have the same area:

$$\mathcal{B}_{\mathrm{ex}} := \int_0^1 B_{\mathrm{ex}}(t)\,\mathrm{d}t \overset{\mathrm{d}}{=} \int_0^1 B_{\mathrm{br}}(t)\,\mathrm{d}t - \min_{0 \leq t \leq 1} B_{\mathrm{br}}(t). \tag{16}$$

This equivalence was noted by Takács [44; 45], who considered the functional of a Brownian bridge in (16) in connection with a problem in railway traffic.

Darling [12] considered the maximum (denoted $G$ by him) of the process

$$Y(t) := B_{\mathrm{br}}(t) - \int_0^1 B_{\mathrm{br}}(s)\,\mathrm{d}s, \qquad 0 \leq t \leq 1. \tag{17}$$

Since $B_{\mathrm{br}}$ is symmetric, $B_{\mathrm{br}} \overset{\mathrm{d}}{=} -B_{\mathrm{br}}$, we obtain from (16) also

$$\mathcal{B}_{\mathrm{ex}} \overset{\mathrm{d}}{=} \max_{0 \leq t \leq 1} B_{\mathrm{br}}(t) - \int_0^1 B_{\mathrm{br}}(t)\,\mathrm{d}t = \max_{0 \leq t \leq 1} Y(t). \tag{18}$$

Thus Darling's $G := \max_t Y(t)$ [12] equals $\mathcal{B}_{\mathrm{ex}}$ (in distribution). $Y(t)$ is clearly a continuous Gaussian process with mean 0, and a straightforward calculation yields

$$\operatorname{Cov}\big(Y(s), Y(t)\big) = \tfrac{1}{2}(|t-s| - \tfrac{1}{2})^2 - \tfrac{1}{24}, \qquad s, t \in [0,1]. \tag{19}$$

It follows that if we regard $[0,1]$ as a circle, or if we extend $Y$ periodically to $\mathbb{R}$, $Y$ is a *stationary* Gaussian process, as observed by Watson [56]. Furthermore, by construction, the integral $\int_0^1 Y(t)\,\mathrm{d}t$ vanishes identically.

## 4. Graph enumeration

Let $C(n,q)$ be the number of connected graphs with $n$ given (labelled) vertices and $q$ edges. Recall Cayley's formula $C(n, n-1) = n^{n-2}$ for every $n \geq 1$. Wright [58] proved that for any fixed $k \geq -1$, we have the analogous asymptotic formula

$$C(n, n+k) \sim \rho_k n^{n+(3k-1)/2} \qquad \text{as } n \to \infty, \tag{20}$$

for some constants $\rho_k$ given by

$$\rho_k = \frac{2^{(1-3k)/2}\pi^{1/2}}{\Gamma(3k/2+1)}\sigma_k, \qquad k \ge -1, \tag{21}$$

with other constants $\sigma_k$ given by $\sigma_{-1} = -1/2$, $\sigma_0 = 1/4$, $\sigma_1 = 5/16$, and the quadratic recursion relation

$$\sigma_{k+1} = \frac{3(k+1)}{2}\sigma_k + \sum_{j=1}^{k-1}\sigma_j\sigma_{k-j}, \qquad k \ge 1. \tag{22}$$

Note the equivalent recursion formula

$$\sigma_{k+1} = \frac{3k+2}{2}\sigma_k + \sum_{j=0}^{k}\sigma_j\sigma_{k-j}, \qquad k \ge -1. \tag{23}$$

Wright gives in the later paper [59] the same result in the form

$$\rho_k = \frac{2^{(1-5k)/2}3^k\pi^{1/2}(k-1)!}{\Gamma(3k/2)}d_k, \qquad k \ge 1, \tag{24}$$

(although he now uses the notation $f_k = \rho_k$; we have further corrected a typo in [59, Theorem 2]), where $d_1 = 5/36$ and

$$d_{k+1} = d_k + \sum_{j=1}^{k-1}\frac{d_j d_{k-j}}{(k+1)\binom{k}{j}}, \qquad k \ge 1. \tag{25}$$

See also Bender, Canfield and McKay [5, Corollaries 1 and 2], which gives the result using the same $d_k$ and further numbers $w_k$ defined by $w_0 = \pi/\sqrt{6}$ and

$$w_k = \frac{(8/3)^{1/2}\pi(k-1)!}{\Gamma(3k/2)}\Bigl(\frac{27k}{8e}\Bigr)^{k/2}d_k, \qquad k \ge 1, \tag{26}$$

so that

$$\rho_k = \frac{3^{1/2}}{2\pi^{1/2}}\Bigl(\frac{e}{12k}\Bigr)^{k/2}w_k. \qquad k \ge 0. \tag{27}$$

(Wright [59] and Bender et al. [5] further consider extensions to the case $k \to \infty$, which does not interest us here.)

In the form

$$\rho_k = \frac{2^{-(5k+1)/2}3^{k+1}\pi^{1/2}k!}{\Gamma(3k/2+1)}d_k, \tag{28}$$

(24) holds for all $k \ge 0$, with $d_0 = 1/6$.

Wright's two versions (21), (22) [58] and (24), (25) [59] are equivalent and we have the relation

$$\sigma_k = \Bigl(\frac{3}{2}\Bigr)^{k+1}k!\,d_k, \qquad k \ge 0. \tag{29}$$

Next, define $c_k$, $k \geq 1$, as in Janson, Knuth, Łuczak and Pittel [24, §8]; $c_k$ is the coefficient for the leading term in an expansion of the generating function for connected graphs (or multigraphs) with $n$ vertices and $n+k$ edges. (Note that $c_k = c_{k0} = \hat{c}_{k0}$ in [24, §8]. $c_k$ is denoted $c_{k,-3k}$ in Wright [58] and $b_k$ in Wright [59].) We have by Wright [59, §5], or by comparing (25) and (32) below,

$$c_k = \left(\frac{3}{2}\right)^k (k-1)!\, d_k, \qquad k \geq 1. \tag{30}$$

From [24, (8.12)] (which is equivalent to Wright [58, (7)]) follows the recursion

$$3rc_r = \tfrac{1}{2}(3r-1)(3r-3)c_{r-1} + \tfrac{9}{2}\sum_{j=0}^{r-1} j(r-1-j)c_j c_{r-1-j}, \qquad r \geq 1, \tag{31}$$

where $jc_j$ is interpreted as $1/6$ when $j=0$. In other words, $c_1 = 5/24$ and

$$rc_r = \tfrac{1}{2}r(3r-3)c_{r-1} + \tfrac{3}{2}\sum_{j=1}^{r-2} j(r-1-j)c_j c_{r-1-j}, \qquad r \geq 2. \tag{32}$$

By (29), (30) is equivalent to

$$\sigma_k = \frac{3}{2}kc_k, \qquad k \geq 1 \tag{33}$$

(and also for $k=0$ with the interpretation $0c_0 = 1/6$ again).

By [24, §§3 and 8], (20) holds for $k \geq 1$ with

$$\rho_k = \frac{2^{(1-3k)/2}\sqrt{\pi}}{\Gamma(3k/2)}c_k, \qquad k \geq 1, \tag{34}$$

which clearly is equivalent to (21) and (24) by (33) and (30).

Finally, we note that (20) can be written

$$C(n, n+k-1) \sim \rho_{k-1} n^{n+3k/2-2} \qquad k \geq 0. \tag{35}$$

## 5. The connection

The connection between Brownian excursion area and graph enumeration was found by Spencer [43], who gave a new proof of (35), and thus (20), that further shows

$$\rho_{k-1} = \frac{\mathbb{E}\,\mathcal{B}_{\mathrm{ex}}^k}{k!}, \qquad k \geq 0. \tag{36}$$

See also Aldous [2, §6].

By (6) and (21) we thus have

$$\sigma_{k-1} = 2^k K_k, \qquad k \geq 0, \tag{37}$$

and it is easily seen that (23) is equivalent to (7). Similarly, (9) is equivalent to (23) and

$$\Omega_k = 2k!\,\sigma_{k-1}, \qquad k \geq 0, \tag{38}$$

and (8) is equivalent to (21) by (2), (36) and (38).

## 6. Further relations

Further relations are immediately obtained by combining the ones above; we give some examples.

By (38) and (12), or by comparing (13) and (23),

$$\omega_k^* = 2^{2k}\sigma_{k-1}, \qquad k \ge 0, \tag{39}$$

and thus also, by (12),

$$\omega_k = 2\sigma_{k-1}, \qquad k \ge 0. \tag{40}$$

By (38) and (33), or comparing (31) and (9), we find

$$\Omega_k = 3(k-1)k!\, c_{k-1}, \qquad k \ge 2. \tag{41}$$

By (12) and (41), or by (39), (40) and (33), or by (31) and (13),

$$\omega_k^* = 2^{2k-1}\, 3(k-1)c_{k-1}, \qquad k \ge 2, \tag{42}$$

$$\omega_k = 3(k-1)c_{k-1}, \qquad k \ge 2. \tag{43}$$

By (29) and (37),

$$d_{k-1} = \Bigl(\frac{4}{3}\Bigr)^k \frac{K_k}{(k-1)!}, \qquad k \ge 1. \tag{44}$$

By (40) and (29), or by (11) and (44),

$$\omega_k = 2\Bigl(\frac{3}{2}\Bigr)^k (k-1)!\, d_{k-1}, \qquad k \ge 1. \tag{45}$$

By (21) and (39),

$$\rho_{k-1} = \frac{2^{2-3k/2}\sqrt{\pi}}{\Gamma((3k-1)/2)}\sigma_{k-1} = \frac{2^{2-7k/2}\sqrt{\pi}}{\Gamma((3k-1)/2)}\omega_k^*, \qquad k \ge 0. \tag{46}$$

Note further that (34) can be written

$$\rho_{k-1} = \frac{2^{2-3k/2}\sqrt{\pi}}{\Gamma(3(k-1)/2)}c_{k-1}, \qquad k \ge 2. \tag{47}$$

By (14) and (42), or by (3), (36) and (47),

$$\mathbb{E}\,\xi^k = \frac{2^{1-k/2}3\sqrt{\pi}\,k!}{\Gamma((3k-1)/2)}(k-1)c_{k-1} = \frac{2^{2-k/2}\sqrt{\pi}\,k!}{\Gamma(3(k-1)/2)}c_{k-1}, \qquad k \ge 2, \tag{48}$$

as claimed in Janson [22, Remark 2.5].

## 7. Asymptotics

Wright [59] proved that the limit $\lim_{k\to\infty} d_k$ exists, and gave the approximation 0.159155. The limit was later identified by Bagaev and Dmitriev [4] as $1/(2\pi)$, i.e.

$$d_k \to \frac{1}{2\pi} \qquad \text{as } k \to \infty. \tag{49}$$

See [24, p. 262] for further history and references.

It follows by Stirling's formula that for $w_k$ in (26),

$$w_k \to 1 \qquad \text{as } k \to \infty. \tag{50}$$

Hence, by (27)

$$\rho_k \sim \frac{3^{1/2}}{2\pi^{1/2}}\Big(\frac{e}{12k}\Big)^{k/2} \qquad \text{as } k \to \infty, \tag{51}$$

and, equivalently,

$$\rho_{k-1} \sim 3\pi^{-1/2}k^{1/2}\Big(\frac{e}{12k}\Big)^{k/2} \qquad \text{as } k \to \infty. \tag{52}$$

By (36) follows further, as stated in Takács [44; 45; 46; 47; 48], see also Janson and Louchard [25],

$$\mathbb{E}\,\mathcal{B}_{\text{ex}}^k \sim 3\sqrt{2}\,k\Big(\frac{k}{12e}\Big)^{k/2} \qquad \text{as } k \to \infty, \tag{53}$$

and equivalently, as stated by Flajolet and Louchard [17] and Chassaing and Janson [23], respectively,

$$\mathbb{E}\,\mathcal{A}^k \sim 2^{1/2}3k\Big(\frac{2k}{3e}\Big)^{k/2} \qquad \text{as } k \to \infty, \tag{54}$$

$$\mathbb{E}\,\xi^k \sim 2^{1/2}3k\Big(\frac{k}{3e}\Big)^{k/2} \qquad \text{as } k \to \infty. \tag{55}$$

By (30) and (49), or by [24, Theorem 8.2],

$$c_k \sim \frac{1}{2\pi}\Big(\frac{3}{2}\Big)^k (k-1)! \qquad \text{as } k \to \infty. \tag{56}$$

Further, by (29) and (49),

$$\sigma_k \sim \frac{1}{2\pi}\Big(\frac{3}{2}\Big)^{k+1} k! \qquad \text{as } k \to \infty, \tag{57}$$

and by (44) and (49), as stated in Takács [44; 45; 46; 47; 48],

$$K_k \sim \frac{1}{2\pi}\Big(\frac{3}{4}\Big)^k (k-1)! \qquad \text{as } k \to \infty. \tag{58}$$

| | | | | |
|---|---|---|---|---|
| $\mathbb{E}\,\mathcal{B}_{\mathrm{ex}}^0 = 1$ | $\mathbb{E}\,\mathcal{B}_{\mathrm{ex}} = \frac{\sqrt{2\pi}}{4}$ | $\mathbb{E}\,\mathcal{B}_{\mathrm{ex}}^2 = \frac{5}{12}$ | $\mathbb{E}\,\mathcal{B}_{\mathrm{ex}}^3 = \frac{15\sqrt{2\pi}}{128}$ | $\mathbb{E}\,\mathcal{B}_{\mathrm{ex}}^4 = \frac{221}{1008}$ |
| $\mathbb{E}\,\mathcal{A}^0 = 1$ | $\mathbb{E}\,\mathcal{A} = \sqrt{\pi}$ | $\mathbb{E}\,\mathcal{A}^2 = \frac{10}{3}$ | $\mathbb{E}\,\mathcal{A}^3 = \frac{15\sqrt{\pi}}{4}$ | $\mathbb{E}\,\mathcal{A}^4 = \frac{884}{63}$ |
| $\gamma_0 = -1$ | $\gamma_1 = 9$ | $\gamma_2 = 405$ | $\gamma_3 = 65610$ | $\gamma_4 = 21749715$ |
| $K_0 = -\frac{1}{2}$ | $K_1 = \frac{1}{8}$ | $K_2 = \frac{5}{64}$ | $K_3 = \frac{15}{128}$ | $K_4 = \frac{1105}{4096}$ |
| $\Omega_0 = -1$ | $\Omega_1 = \frac{1}{2}$ | $\Omega_2 = \frac{5}{4}$ | $\Omega_3 = \frac{45}{4}$ | $\Omega_4 = \frac{3315}{16}$ |
| $\omega_0 = -1$ | $\omega_1 = \frac{1}{2}$ | $\omega_2 = \frac{5}{8}$ | $\omega_3 = \frac{15}{8}$ | $\omega_4 = \frac{1105}{128}$ |
| $\omega_0^* = -\frac{1}{2}$ | $\omega_1^* = 1$ | $\omega_2^* = 5$ | $\omega_3^* = 60$ | $\omega_4^* = 1105$ |
| $\rho_{-1} = 1$ | $\rho_0 = \frac{\sqrt{2\pi}}{4}$ | $\rho_1 = \frac{5}{24}$ | $\rho_2 = \frac{5\sqrt{2\pi}}{256}$ | $\rho_3 = \frac{221}{24192}$ |
| $\sigma_{-1} = -\frac{1}{2}$ | $\sigma_0 = \frac{1}{4}$ | $\sigma_1 = \frac{5}{16}$ | $\sigma_2 = \frac{15}{16}$ | $\sigma_3 = \frac{1105}{256}$ |
| | $d_0 = \frac{1}{6}$ | $d_1 = \frac{5}{36}$ | $d_2 = \frac{5}{36}$ | $d_3 = \frac{1105}{7776}$ |
| | | $c_1 = \frac{5}{24}$ | $c_2 = \frac{5}{16}$ | $c_3 = \frac{1105}{1152}$ |

TABLE 1
*Some numerical values for the Brownian excursion and for graph enumeration.*

Similarly by (57) and (39), (40) and (11),

$$\omega_k^* \sim \frac{1}{2\pi}\, 6^k (k-1)! \qquad \text{as } k \to \infty, \tag{59}$$

$$\omega_k \sim \frac{1}{\pi}\Big(\frac{3}{2}\Big)^k (k-1)! \qquad \text{as } k \to \infty, \tag{60}$$

$$\Omega_k \sim \frac{1}{\pi k}\Big(\frac{3}{2}\Big)^k (k!)^2 \sim 2\Big(\frac{3k^2}{2e^2}\Big)^k \qquad \text{as } k \to \infty. \tag{61}$$

## 8. Numerical values

Numerical values for small $k$ are given in Table 1. See further Louchard [34] ($\mathbb{E}\,\mathcal{B}_{\mathrm{ex}}^k$), Takács [44; 47; 48] ($\mathbb{E}\,\mathcal{B}_{\mathrm{ex}}^k$, $K_k$), Janson, Knuth, Łuczak and Pittel [24, p. 259 or 262] ($c_k$), Flajolet, Poblete and Viola [18, Table 1 and p. 503] ($\mathbb{E}\,\mathcal{A}^k$, $\Omega_k$, $\omega_k$, $\omega_k^*$), Flajolet and Louchard [17, Table 1] ($\mathbb{E}\,\mathcal{A}^k$, $\Omega_k$), and Janson [21, p. 343] ($\omega_k^*$).

## 9. Power series

Define the formal power series

$$C(z) := \sum_{r=1}^{\infty} c_r z^r. \tag{62}$$

By [24, §8 and (7.2)], we then have

$$e^{C(z)} = \sum_{r=0}^{\infty} e_r z^r, \tag{63}$$

where

$$\begin{aligned} e_r &:= \frac{(6r)!}{2^{5r}3^{2r}(3r)!\,(2r)!} = \Bigl(\frac{3}{2}\Bigr)^r \frac{\Gamma(r+5/6)\Gamma(r+1/6)}{2\pi r!} = \Bigl(\frac{3}{2}\Bigr)^r \frac{(5/6)^{\overline{r}}(1/6)^{\overline{r}}}{r!} \\ &= 18^{-r}\frac{\Gamma(3r+1/2)}{\Gamma(r+1/2)\,r!}. \end{aligned} \tag{64}$$

(The last formula follows by the triplication formula for the Gamma function, or by induction.) We have $e_0 = 1$, $e_1 = 5/24$, $e_2 = 385/1152$, $e_3 = 85085/82944$. The constants $e_r$ are the coefficient for the leading term in an expansion of the generating function for all graphs (or multigraphs) with $n$ vertices and $n + k$ edges, see [24, §7] (where $e_r$ is denoted $e_{r0}$).

## 10. Linear recursions

From (63) follows the linear recursion, see [24, §8],

$$c_r = e_r - \frac{1}{r}\sum_{j=1}^{r-1} j c_j e_{r-j}, \qquad r \ge 1, \tag{65}$$

where $e_k$ are given explicitly by (64).

By (41) together with (64) and simple calculations, (65) is equivalent to

$$18^r\Omega_r = \frac{12r}{6r-1}\frac{\Gamma(3r+1/2)}{\Gamma(r+1/2)} - \sum_{j=1}^{r-1}\binom{r}{j}\frac{\Gamma(3j+1/2)}{\Gamma(j+1/2)}18^{r-j}\Omega_{r-j}, \qquad r \ge 1, \tag{66}$$

given by Flajolet and Louchard [17], and, equivalently, see (15),

$$\gamma_r = \frac{12r}{6r-1}\frac{\Gamma(3r+1/2)}{\Gamma(r+1/2)} - \sum_{j=1}^{r-1}\binom{r}{j}\frac{\Gamma(3j+1/2)}{\Gamma(j+1/2)}\gamma_{r-j}, \qquad r \ge 1, \tag{67}$$

given by Louchard [34]; the latter was already given as (5) above. Changing the upper summation limit we can also write these as

$$18^r\Omega_r = \frac{6r+1}{6r-1}\frac{\Gamma(3r+1/2)}{\Gamma(r+1/2)} - \sum_{j=1}^{r}\binom{r}{j}\frac{\Gamma(3j+1/2)}{\Gamma(j+1/2)}18^{r-j}\Omega_{r-j}, \quad r \ge 0, \tag{68}$$

$$\gamma_r = \frac{6r+1}{6r-1}\,\frac{\Gamma(3r+1/2)}{\Gamma(r+1/2)} - \sum_{j=1}^{r}\binom{r}{j}\frac{\Gamma(3j+1/2)}{\Gamma(j+1/2)}\gamma_{r-j}, \qquad r \ge 0. \tag{69}$$

By (10), these are further equivalent to the linear recursion in Takács [44; 45]

$$K_r = \frac{6r+1}{2(6r-1)}\,\alpha_r - \sum_{j=1}^{r}\alpha_j K_{r-j}, \qquad r \ge 1, \tag{70}$$

where, cf. (64),

$$\alpha_j := 36^{-j}\frac{\Gamma(3j+1/2)}{\Gamma(j+1/2)\,j!} = 2^{-j}e_j, \qquad j \ge 0. \tag{71}$$

(Perman and Wellner [38] and Majumdar and Comtet [35] use $\gamma_j$ for $\alpha_j$.)

## 11. Semi-invariants

Takács [44] considers also the semi-invariants (cumulants) $\Lambda_n$ of $\mathcal{B}_{\mathrm{ex}}$, given recursively from the moments $M_k := \mathbb{E}\,\mathcal{B}_{\mathrm{ex}}^k$ by the general formula

$$\Lambda_n = M_n - \sum_{k=1}^{n-1}\binom{n-1}{k}M_k\Lambda_{n-k}, \qquad n \ge 1. \tag{72}$$

For example, $\Lambda_1 = \mathbb{E}\,\mathcal{B}_{\mathrm{ex}} = \sqrt{\pi/8}$, $\Lambda_2 = \operatorname{Var}(\mathcal{B}_{\mathrm{ex}}) = (10-3\pi)/24$ and $\Lambda_3 = \mathbb{E}\,\mathcal{B}_{\mathrm{ex}}^3 - 3\,\mathbb{E}\,\mathcal{B}_{\mathrm{ex}}^2\,\mathbb{E}\,\mathcal{B}_{\mathrm{ex}} + 2(\mathbb{E}\,\mathcal{B}_{\mathrm{ex}})^3 = (8\pi-25)\sqrt{2\pi}/128$. Since these do not have a nice form, they will not be considered further here.

## 12. The Airy function

There are many connections between the distribution of $\mathcal{B}_{\mathrm{ex}}$ and the transcendental *Airy function*; these connections have been noted in different contexts by several authors and, in particular, studied in detail by Flajolet and Louchard [17]. We describe many of them in the following sections.

The Airy function $\operatorname{Ai}(x)$ is defined by, for example, the conditionally convergent integral

$$\operatorname{Ai}(x) := \frac{1}{\pi}\int_0^\infty \cos\bigl(t^3/3 + xt\bigr)\,\mathrm{d}t, \qquad -\infty < x < \infty, \tag{73}$$

see e.g. Abramowitz and Stegun [1, §10.4] or Lebedev [32, §5.17], which also give the following representations using standard Bessel functions

$$\operatorname{Ai}(x) = \frac{1}{\pi}\Bigl(\frac{x}{3}\Bigr)^{1/2}K_{1/3}\Bigl(\frac{2x^{3/2}}{3}\Bigr), \qquad x > 0, \tag{74}$$

$$\operatorname{Ai}(x) = \frac{x^{1/2}}{3}\Bigl(I_{-1/3}\Bigl(\frac{2x^{3/2}}{3}\Bigr) - I_{1/3}\Bigl(\frac{2x^{3/2}}{3}\Bigr)\Bigr), \qquad x > 0, \tag{75}$$

$$\mathrm{Ai}(-x) = \frac{x^{1/2}}{3}\Bigl(J_{-1/3}\Bigl(\frac{2x^{3/2}}{3}\Bigr) + J_{1/3}\Bigl(\frac{2x^{3/2}}{3}\Bigr)\Bigr), \qquad x > 0. \tag{76}$$

$\mathrm{Ai}(x)$ is, up to normalization, the unique solution of the differential equation $\mathrm{Ai}''(x) = x\mathrm{Ai}(x)$ that is bounded for $x \ge 0$. Ai extends to an entire function. All zeros of $\mathrm{Ai}(x)$ lie on the negative real axis; they will appear several times below, and we denote them by $a_j = -|a_j|$, $j = 1, 2, \dots$, with $0 < |a_1| < |a_2|, \dots$. In other words,

$$\mathrm{Ai}(a_j) = \mathrm{Ai}(-|a_j|) = 0, \qquad j = 1, 2, \dots \tag{77}$$

We have the asymptotic formula $a_j \sim -(3\pi/2)^{2/3} j^{2/3}$, which can be refined to an asymptotic expansion [1, 10.4.94], [17]. (Abramowitz and Stegun [1, 10.4.94] and Louchard [34] use $a_j$ as we do, while Takács [44; 45; 46; 47; 48; 50] uses $a_j$ for $|a_j|$; Flajolet and Louchard [17] and Majumdar and Comtet [35] use $\alpha_j$ for $|a_j|$; Darling [12] denotes our $|a_j|$ by $\sigma_j$; his $\alpha_j$ are the zeros of $J_{1/3}(x) + J_{-1/3}(x)$, and are equal to $\frac{2}{3}|a_j|^{3/2}$ in our notation, see (76).)

We will later also need both the derivative Ai′ and the integral of Ai. It seems that there is no standard notation for the latter, and we will use

$$\mathrm{AI}(x) := \int_x^\infty \mathrm{Ai}(t)\,\mathrm{d}t = \frac{1}{3} - \int_0^x \mathrm{Ai}(t)\,\mathrm{d}t, \tag{78}$$

using $\int_0^\infty \mathrm{Ai}(x)\,\mathrm{d}x = 1/3$ [1, 10.4.82]. See further Appendix A.

## 13. Laplace transform

Let $\psi_{\mathrm{ex}}(t) := \mathbb{E}\, e^{-t\mathcal{B}_{\mathrm{ex}}}$ be the Laplace transform of $\mathcal{B}_{\mathrm{ex}}$; thus $\psi_{\mathrm{ex}}(-t) := \mathbb{E}\, e^{t\mathcal{B}_{\mathrm{ex}}}$ is the moment generating function of $\mathcal{B}_{\mathrm{ex}}$. It follows from (53) that $\psi_{\mathrm{ex}}(t)$ exists for all complex $t$, and is thus an entire function, with

$$\psi_{\mathrm{ex}}(t) = \sum_{k=0}^\infty \mathbb{E}\, \mathcal{B}_{\mathrm{ex}}^k \frac{(-t)^k}{k!}\,. \tag{79}$$

($\psi_{\mathrm{ex}}(t)$ is denoted $\phi(t)$ by Louchard [34], $G(2^{-3/2}t)$ by Flajolet and Louchard [17] and $\Psi_e(t)$ by Perman and Wellner [38].)

Darling [12] found (in the context of (18)) the formula, with $|a_j|$ as in (77),

$$\psi_{\mathrm{ex}}(t) = \sqrt{2\pi}\, t \sum_{j=1}^\infty \exp\bigl(-2^{-1/3}|a_j| t^{2/3}\bigr), \qquad t > 0 \tag{80}$$

also found by Louchard [34] using the formula, for $x \ge 0$,

$$\int_0^\infty \bigl(e^{-xt} - 1\bigr)\psi_{\mathrm{ex}}(t^{3/2}) \frac{\mathrm{d}t}{\sqrt{2\pi t^3}} = 2^{1/3}\left(\frac{\mathrm{Ai}'(2^{1/3}x)}{\mathrm{Ai}(2^{1/3}x)} - \frac{\mathrm{Ai}'(0)}{\mathrm{Ai}(0)}\right), \tag{81}$$

proved by Louchard [33]. We give a (related but somewhat different) proof of (81) in Appendix C.4. The relation (80) follows easily by Laplace inversion from (81) and the partial fraction expansion, see [17],

$$\frac{\mathrm{Ai}'(z)}{\mathrm{Ai}(z)} - \frac{\mathrm{Ai}'(0)}{\mathrm{Ai}(0)} = \sum_{j=1}^{\infty} \left( \frac{1}{z - a_j} + \frac{1}{a_j} \right). \tag{82}$$

A proof of (80) by methods from mathematical physics is given by Majumdar and Comtet [35].

Taking the derivative with respect to $x$ in (81), we find, for $x > 0$,

$$\int_0^\infty e^{-xt} \psi_{\mathrm{ex}}(t^{3/2}) \frac{\mathrm{d}t}{\sqrt{2\pi t}} = -2^{1/3} \frac{\mathrm{d}}{\mathrm{d}x} \left( \frac{\mathrm{Ai}'(2^{1/3}x)}{\mathrm{Ai}(2^{1/3}x)} \right) = \left( 2^{1/3} \frac{\mathrm{Ai}'(2^{1/3}x)}{\mathrm{Ai}(2^{1/3}x)} \right)^2 - 2x, \tag{83}$$

or, by the changes of variables $x \mapsto 2^{-1/3}x$, $t \mapsto 2^{1/3}x$, as in [34],

$$\int_0^\infty e^{-xt} \psi_{\mathrm{ex}}(\sqrt{2}\, t^{3/2}) \frac{\mathrm{d}t}{\sqrt{\pi t}} = -2 \frac{\mathrm{d}}{\mathrm{d}x} \left( \frac{\mathrm{Ai}'(x)}{\mathrm{Ai}(x)} \right) = 2 \left( \frac{\mathrm{Ai}'(x)}{\mathrm{Ai}(x)} \right)^2 - 2x. \tag{84}$$

## 14. Series expansions for distribution and density functions

Let $F_{\mathrm{ex}}(x) := \mathbb{P}(\mathcal{B}_{\mathrm{ex}} \le x)$ and $f_{\mathrm{ex}}(x) := F'(x)$ be the distribution and density functions of $\mathcal{B}_{\mathrm{ex}}$. (Takács [44; 45; 46; 47; 48] uses $W$ and $W'$ for these, while Flajolet and Louchard [17] use $W$ and $w$ for the distribution and density functions of $\mathcal{A} = 2^{3/2}\mathcal{B}_{\mathrm{ex}}$; thus their $W(x) = F_{\mathrm{ex}}(2^{-3/2}x)$ and $w(x) = 2^{-3/2} f_{\mathrm{ex}}(2^{-3/2}x)$.) Since $\mathcal{B}_{\mathrm{ex}} > 0$ a.s., we consider in this section $x > 0$ only.

Let, as above, $a_j$ denote the zeros of the Airy function. Darling [12] found by Laplace inversion from (80) (in our notation)

$$F_{\mathrm{ex}}(x) = 2\sqrt{\pi} \sum_{j=1}^{\infty} |a_j|^{-3/2} p\big(\sqrt{2}|a_j|^{-3/2} x\big), \tag{85}$$

where $p$ is the density of the positive stable distribution with exponent $2/3$, normalized to have Laplace transform $\exp(-t^{2/3})$; in the notation of Feller [15, Section XVII.6]

$$p(x) = p\big(x; \tfrac{2}{3}, -\tfrac{2}{3}\big) = \frac{1}{\pi x} \sum_{k=1}^{\infty} (-1)^{k-1} \sin(2k\pi/3) \frac{\Gamma(2k/3+1)}{k!} x^{-2k/3}. \tag{86}$$

By Feller [15, Lemma XVII.6.2] also

$$p(x) = x^{-5/3} p\big(x^{-2/3}; \tfrac{3}{2}, -\tfrac{1}{2}\big), \tag{87}$$

where $p(x; \frac{3}{2}, -\frac{1}{2}) = \frac{1}{\pi x} \sum_{k=1}^{\infty} (-1)^{k-1} \sin(2k\pi/3) \frac{\Gamma(2k/3+1)}{k!} x^k$ is the density of a spectrally negative stable distribution with exponent $3/2$. (Takács [44] defines the function $g(x) = 2^{-1/3} p(2^{-1/3}x; \frac{3}{2}, -\frac{1}{2})$; this is another such density function with a different normalization. Takacs [49] uses $g$ for $p$.) The stable density function $p$ can also be expressed by a Whittaker function $W$ or a confluent hypergeometric function $U$ [1, 13.1], [32, 9.10 and 9.13.11] (where $U$ is denoted $\Psi$),

$$p(x) = \Big(\frac{3}{\pi}\Big)^{1/2} x^{-1} \exp\Big(-\frac{2}{27x^2}\Big) W_{\frac{1}{2},\frac{1}{6}}\Big(\frac{4}{27x^2}\Big) \tag{88}$$

$$= \frac{2^{4/3}}{3^{3/2}\pi^{1/2}}\, x^{-7/3} \exp\Big(-\frac{4}{27x^2}\Big) U\Big(\frac{1}{6}, \frac{4}{3}, \frac{4}{27x^2}\Big). \tag{89}$$

Hence (85) can be written [44; 45; 46; 47; 48]

$$F_{\mathrm{ex}}(x) = 2^{7/6} 3^{-3/2} x^{-7/3} \sum_{j=1}^{\infty} a_j^2 \exp\Big(-\frac{2|a_j|^3}{27x^2}\Big) U\Big(\frac{1}{6}, \frac{4}{3}; \frac{2|a_j|^3}{27x^2}\Big). \tag{90}$$

$$= \frac{\sqrt{6}}{x} \sum_{j=1}^{\infty} v_j^{2/3} e^{-v_j} U\Big(\frac{1}{6}, \frac{4}{3}; v_j\Big), \tag{91}$$

with $v_j = 2|a_j|^3/27x^2$, which leads to [17; 44; 45; 46; 47; 48]

$$f_{\mathrm{ex}}(x) = \frac{2\sqrt{6}}{x^2} \sum_{j=1}^{\infty} v_j^{2/3} e^{-v_j} U\Big(-\frac{5}{6}, \frac{4}{3}; v_j\Big). \tag{92}$$

## 15. Asymptotics of distribution and density functions

Louchard [34] gave the two first terms in an asymptotic expansion of $f_{\mathrm{ex}}(x)$ and the first term for $F_{\mathrm{ex}}(x)$ as $x \to 0$; it was observed by Flajolet and Louchard [17] that full asymptotic expansions readily follow from (92) (for $f_{\mathrm{ex}}$; the result for $F_{\mathrm{ex}}$ follows similarly from (90), or by integration); note that only the term with $j = 1$ in (92) is significant as $x \to 0$ and use the asymptotic expansion of $U$ given in e.g. [32, (9.12.3)]. The first terms are (correcting typos in [34] and [17]), as $x \to 0$,

$$f_{\mathrm{ex}}(x) \sim e^{-\frac{2|a_1|^3}{27x^2}} \Big(\tfrac{8}{81}|a_1|^{9/2}x^{-5} - \tfrac{35}{27}|a_1|^{3/2}x^{-3} - \tfrac{35}{144}|a_1|^{-3/2}x^{-1} + \dots\Big), \tag{93}$$

$$F_{\mathrm{ex}}(x) \sim e^{-\frac{2|a_1|^3}{27x^2}} \Big(\tfrac{2}{3}|a_1|^{3/2}x^{-2} + \tfrac{1}{4}|a_1|^{-3/2} - \tfrac{105}{64}|a_1|^{-9/2}x^2 + \dots\Big). \tag{94}$$

As in all other asymptotic expansions in this paper, we do not claim here that there is a convergent infinite series on the right hand side; the notation (using $\sim$ instead of $=$) signifies only that if we truncate the sum after an arbitrary finite number of terms, the error is smaller order than the last term. (Hence, more precisely, the error is of the order of the first omitted term.) In fact, the series in (93) and (94) diverge for every $x > 0$ because the asymptotic series for $U$ does so.

For $x \to \infty$, it was observed by Csörgő, Shi and Yor [11] that (53) implies

$$-\ln \mathbb{P}(\mathcal{B}_{\mathrm{ex}} > x) = -\ln\bigl(1 - F_{\mathrm{ex}}(x)\bigr) \sim -6x^2, \qquad x \to \infty. \tag{95}$$

Another proof was given by Fill and Janson [16]. Much more precise results were obtained by Janson and Louchard [25] who found, as $x \to \infty$, asymptotic expansions beginning with

$$f_{\mathrm{ex}}(x) \sim \frac{72\sqrt{6}}{\sqrt{\pi}} x^2 e^{-6x^2} \left(1 - \frac{1}{9}\,x^{-2} - \frac{5}{1296}\,x^{-4} - \frac{25}{46656}\,x^{-6} + \dots\right), \tag{96}$$

$$1 - F_{\mathrm{ex}}(x) \sim \frac{6\sqrt{6}}{\sqrt{\pi}} x e^{-6x^2} \left(1 - \frac{1}{36}\,x^{-2} - \frac{1}{648}\,x^{-4} - \frac{7}{46656}\,x^{-6} + \dots\right). \tag{97}$$

## 16. Airy function expansions

Following Louchard [34], we expand the left hand side of (81) in an asymptotic series as $x \to \infty$, using (79), and find after replacing $2^{1/3}x$ by $x$ the asymptotic series

$$\frac{\mathrm{Ai}'(x)}{\mathrm{Ai}(x)} \sim -x^{1/2} + \sum_{k=1}^{\infty} (-1)^k \, \mathbb{E}\, \mathcal{B}_{\mathrm{ex}}^k \frac{\Gamma((3k-1)/2)}{2\sqrt{\pi}\,k!} 2^{k/2} x^{(1-3k)/2}. \tag{98}$$

Note that the infinite sum in (98) diverges by (53); cf. the discussion in Section 15. For the proof of (98), we thus cannot simply substitute (79) into (81); instead we substitute a truncated version (a finite Taylor series)

$$\psi_{\mathrm{ex}}(t) = \sum_{k=0}^{N} \mathbb{E}\, \mathcal{B}_{\mathrm{ex}}^k \frac{(-t)^k}{k!} + O(t^{N+1}), \tag{99}$$

where $N$ is finite but arbitrary. In similar situations for other Brownian areas in later sections, for example with (133) below for the Brownian bridge, this yields directly a sum of Gamma integrals and an asymptotic expansion of the desired type. In the present case, however, we have to work a little more to avoid divergent integrals. One possibility is to substitute (99) into the differentiated version (84), evaluate the integrals, subtract $x^{-1/2}$ from both sides, and integrate with respect to $x$. Another possibility is to write the left hand side of (81) as

$$\int_0^\infty \bigl(e^{-xt} - 1\bigr) \frac{\mathrm{d}t}{\sqrt{2\pi t^3}} - \int_0^\infty \bigl(\psi_{\mathrm{ex}}(t^{3/2}) - 1\bigr) \frac{\mathrm{d}t}{\sqrt{2\pi t^3}} + \int_0^\infty e^{-xt} \bigl(\psi_{\mathrm{ex}}(t^{3/2}) - 1\bigr) \frac{\mathrm{d}t}{\sqrt{2\pi t^3}} \tag{100}$$

and subsitute (99) into the third integral, noting that the first integral is $-\sqrt{2}\,x^{1/2}$ and the second is a constant (necessarily equal to $2^{1/3}\mathrm{Ai}'(0)/\mathrm{Ai}(0)$).

Using (6), we can rewrite (98) as

$$\frac{\mathrm{Ai}'(x)}{\mathrm{Ai}(x)} \sim \sum_{k=0}^{\infty}(-1)^k 2K_k x^{1/2-3k/2}, \qquad x\to\infty. \tag{101}$$

This is also, with a change of variables, given in Takács [45]. Flajolet and Louchard [17] give this expansion with the coefficients $(-1)^k\Omega_k/(2^k k!)$, which is equivalent by (10). They give also related expansions involving Bessel and confluent hypergeometric functions with coefficients $\Omega_k/k!$. An equivalent expansion with coeficients $3(k-1)c_{k-1}$, which equals $\Omega_k/k!$ by (41), was given by Voblyĭ [55], see [24, (8.14) and (8.15)].

The Airy function and its derivative have, as $x\to\infty$, the asymptotic expansions [1, 10.4.59 and 10.4.61] (we write $c_k'$ and $d_k'$ instead of $c_k$ and $d_k$ to avoid confusion with our $c_k$ and $d_k$ above)

$$\mathrm{Ai}(x) \sim \tfrac{1}{2\sqrt{\pi}} x^{-1/4} e^{-2x^{3/2}/3}\sum_{k=0}^{\infty}(-1)^k c_k'\Big(\frac{3}{2}\Big)^k x^{-3k/2}, \tag{102}$$

$$\mathrm{Ai}'(x) \sim -\tfrac{1}{2\sqrt{\pi}} x^{1/4} e^{-2x^{3/2}/3}\sum_{k=0}^{\infty}(-1)^k d_k'\Big(\frac{3}{2}\Big)^k x^{-3k/2}, \tag{103}$$

where

$$c_k' = \frac{\Gamma(3k+\frac12)}{54^k k!\,\Gamma(k+\frac12)}, \qquad d_k' = -\frac{6k+1}{6k-1}c_k', \qquad k\ge 0. \tag{104}$$

These asymptotic expansions can also be written using $e_k$ or $\alpha_k$, since, by (64), (71) and (104),

$$3^k c_k' = e_k, \tag{105}$$

$$(3/2)^k c_k' = \alpha_k, \tag{106}$$

$$(3/2)^k d_k' = \alpha_k' := -\frac{6k+1}{6k-1}\alpha_k. \tag{107}$$

We thus obtain from (101)–(103), after the change of variables $z=-x^{-3/2}\nearrow 0$, the equality for formal power series

$$\sum_{k=0}^{\infty} 2K_k z^k = -\frac{\sum_{k=0}^{\infty}(3/2)^k d_k' z^k}{\sum_{k=0}^{\infty}(3/2)^k c_k' z^k} = -\frac{\sum_{k=0}^{\infty}\alpha_k' z^k}{\sum_{k=0}^{\infty}\alpha_k z^k}. \tag{108}$$

These asymptotic expansions can also be written as hypergeometric series. We have, as is easily verified, the equalities for formal power series

$$\sum_{k=0}^{\infty} c_k' z^k = {}_2F_0\big(5/6, 1/6; z/2\big), \tag{109}$$

$$\sum_{k=0}^{\infty} d_k' z^k = {}_2F_0\big(7/6, -1/6; z/2\big). \tag{110}$$

Hence, (102) and (103) can be written, for $x \to \infty$,

$$\operatorname{Ai}(x) \sim \tfrac{1}{2\sqrt{\pi}} x^{-1/4} e^{-2x^{3/2}/3} {}_2F_0\Bigl(\frac{5}{6}, \frac{1}{6}; -\frac{3}{4x^{3/2}}\Bigr), \tag{111}$$

$$\operatorname{Ai}'(x) \sim -\tfrac{1}{2\sqrt{\pi}} x^{1/4} e^{-2x^{3/2}/3} {}_2F_0\Bigl(\frac{7}{6}, -\frac{1}{6}; -\frac{3}{4x^{3/2}}\Bigr), \tag{112}$$

and (108) can be written

$$\sum_{k=0}^{\infty} 2K_k z^k = -\frac{{}_2F_0\bigl(7/6, -1/6; 3z/4\bigr)}{{}_2F_0\bigl(5/6, 1/6; 3z/4\bigr)}. \tag{113}$$

Flajolet and Louchard [17] give the equivalent formula, see (10),

$$\sum_{k=0}^{\infty} \Omega_k \frac{w^k}{k!} = -\frac{{}_2F_0\bigl(7/6, -1/6; 3w/2\bigr)}{{}_2F_0\bigl(5/6, 1/6; 3w/2\bigr)}. \tag{114}$$

By the asymptotic expansion for Bessel functions [32, (5.11.10)]

$$I_\nu(x) \sim e^x (2\pi x)^{-1/2} {}_2F_0(1/2+\nu, 1/2-\nu; 1/2x), \qquad x \to \infty, \tag{115}$$

this can, as in Flajolet and Louchard [17], also be written (with arbitrary signs)

$$-\frac{I_{\pm 2/3}(x/3)}{I_{\pm 1/3}(x/3)} \sim \sum_{k=0}^{\infty} \Omega_k \frac{x^{-k}}{k!} = \sum_{k=0}^{\infty} \omega_k x^{-k}, \qquad x \to \infty. \tag{116}$$

The coefficients in this asymptotic series can be rewritten in various ways by the relations in Sections 2 and 6, for example by (41) as $3(k-1)c_{k-1}$ for $k \ge 2$, which yields [24, (8.15)].

If we in (108) multiply by the denominator and identify coefficients, we obtain the linear recursion (70) and the equivalent recursions (65)–(69); note that $e_k = 3^k c'_k$ and $\alpha_k = (3/2)^k c'_k$ by (105) and (106).

On the other hand, by differentiating (101) (which is allowed, e.g. because the asymptotic expansion holds in a sector in the complex plane),

$$\begin{aligned} \sum_{k=0}^{\infty} (-1)^k (1-3k) K_k x^{-1/2-3k/2} &\sim \frac{\operatorname{Ai}''(x)}{\operatorname{Ai}(x)} - \left(\frac{\operatorname{Ai}'(x)}{\operatorname{Ai}(x)}\right)^2 \\ &\sim x - \left(\sum_{k=0}^{\infty} (-1)^k 2K_k x^{1/2-3k/2}\right)^2. \end{aligned} \tag{117}$$

This gives the equation for formal power series

$$\sum_{k=0}^{\infty} (3k-1) K_k z^{k+1} = 1 - \left(\sum_{k=0}^{\infty} 2K_k z^k\right)^2, \tag{118}$$

which is equivalent to the quadratic recursion (7).

## 17. A continued fraction

As noted in [24, (8.16)], (116) yields (after substituting $z = 1/x$) the asymptotic expansion

$$2z + \cfrac{1}{8z + \cfrac{1}{14z + \cfrac{1}{20z + \cfrac{1}{26z + \dots}}}} \sim -\sum_{k=0}^{\infty} \Omega_k \frac{z^k}{k!} = -\sum_{k=0}^{\infty} \omega_k z^k, \qquad z \to 0, \tag{119}$$

where the continued fraction has denominators $(6n+2)z$.

## 18. Moments and Airy zeroes

Define, following Flajolet and Louchard [17],

$$\Lambda(s) := \sum_{j=1}^{\infty} |a_j|^{-s}, \qquad \operatorname{Re} s > 3/2. \tag{120}$$

Since $|a_j| \sim (3\pi j/2)^{2/3}$ [1, 10.4.94], the sum converges and $\Lambda$ is analytic for $\operatorname{Re} s > 3/2$. (Flajolet and Louchard [17] call $\Lambda$ the *root zeta function* of the Airy function Ai.)

By (80), for $\operatorname{Re} s > 0$ and with absolutely convergent sums and integrals,

$$\begin{aligned} \int_0^\infty t^{s-1} \psi_{\mathrm{ex}}(t)\,\mathrm{d}t &= \sqrt{2\pi} \int_0^\infty t^{s+1} \sum_{j=1}^{\infty} \exp\bigl(-2^{-1/3}|a_j| t^{2/3}\bigr) \frac{\mathrm{d}t}{t} \\ &= \frac{3\sqrt{2\pi}}{2} \sum_{j=1}^{\infty} \int_0^\infty u^{(3s+3)/2} \exp\bigl(-2^{-1/3}|a_j| u\bigr) \frac{\mathrm{d}u}{u} \\ &= \frac{3\sqrt{2\pi}}{2} 2^{(s+1)/2} \sum_{j=1}^{\infty} |a_j|^{-(3s+3)/2} \Gamma\Bigl(\frac{3s+3}{2}\Bigr) \\ &= 3\sqrt{\pi}\, 2^{s/2} \Gamma\Bigl(\frac{3s+3}{2}\Bigr) \Lambda\Bigl(\frac{3s+3}{2}\Bigr). \end{aligned} \tag{121}$$

In particular, this is finite for all $s > 0$, and thus Lemma B.1 in Appendix B applies and shows that $\mathcal{B}_{\mathrm{ex}}$ has negative moments of all orders. Since all (positive) moments too are finite, $\mathbb{E}\,\mathcal{B}_{\mathrm{ex}}^s$ is an entire function of $s$, and by (121) and (335), $\Lambda$ extends to a meromorphic function in $\mathbb{C}$ and, as shown by Flajolet and Louchard [17],

$$\mathbb{E}\,\mathcal{B}_{\mathrm{ex}}^s = \frac{3\sqrt{\pi}\, 2^{-s/2}}{\Gamma(-s)} \Gamma\Bigl(\frac{3-3s}{2}\Bigr) \Lambda\Bigl(\frac{3-3s}{2}\Bigr), \qquad s \in \mathbb{C}. \tag{122}$$

Since the left hand side has no poles, while the Gamma factors are non-zero, $\Lambda\bigl((3-3s)/2\bigr)$ can have poles only at the poles of $\Gamma(-s)$, i.e. at $s = 0, 1, \dots$; hence, $\Lambda(z)$ can have poles only at $z = \frac{3}{2}$, $0$, $-\frac{3}{2}$, $-3$, $-\frac{9}{2}$, .... Moreover, $\mathbb{E}\,\mathcal{B}_{\mathrm{ex}}^s$

is non-zero for real $s$; hence $\Gamma\bigl((3-3s)/2\bigr)\Lambda\bigl((3-3s)/2\bigr)$ and $\Gamma(-s)$ have exactly the same zeros (i.e., none) and the same poles. Consequently, $\Lambda(z)$ has simple poles at $z = \frac{3}{2}, -\frac{3}{2}, -\frac{9}{2}, \dots$, but is finite and non-zero at $0, -3, -6, \dots$ (where $\Gamma$ has a pole); furthermore, $\Lambda$ is zero at the other poles of $\Gamma$, i.e. at $-1$, $-2$, $-4$, $-5$, $-7$, .... This was shown by Flajolet and Louchard [17] (by a partly different approach, extending $\Lambda(s)$ by explicit formulas).

As shown in Flajolet and Louchard [17], as a simple consequence of (82), for $|z| < |a_1|$,

$$\frac{\mathrm{Ai}'(z)}{\mathrm{Ai}(z)} - \frac{\mathrm{Ai}'(0)}{\mathrm{Ai}(0)} = \sum_{k\ge1} \Lambda(k+1)(-z)^k, \tag{123}$$

and thus the values of $\Lambda(s)$ at positive integers $s = 2, 3, \dots$ can be computed from the Taylor series of Ai, given for example in [1, 10.4.2–5]. This and (122) gives explicit formulas for the negative moments $\mathcal{B}_{\mathrm{ex}}^{-s}$ when $s$ is an odd multiple of 1/3, including when $s$ is an odd integer; see Flajolet and Louchard [17]. Alternatively, these formulas follow from (340) in Appendix B. For example,

$$\mathbb{E}\,\mathcal{B}_{\mathrm{ex}}^{-1/3} = \frac{2^{13/6}3^{2/3}\pi^{5/2}}{\Gamma(1/3)^5} = 2^{-17/6}3^{19/6}\pi^{-5/2}\Gamma(2/3)^5 \approx 1.184, \tag{124}$$

$$\mathbb{E}\,\mathcal{B}_{\mathrm{ex}}^{-1} = 3\sqrt{2\pi}\left(1 - \frac{3^{5/2}\Gamma(2/3)^6}{4\pi^3}\right) \approx 1.693, \tag{125}$$

see further Section 29. We have here used the standard formula

$$\Gamma(1/3)\Gamma(2/3) = \frac{\pi}{\sin(\pi/3)} = \frac{2\pi}{\sqrt{3}}. \tag{126}$$

Note also that if $-1 < \operatorname{Re} s < -1/2$, then by (81), Fubini's theorem and (335) in Appendix B,

$$\begin{aligned}
2^{1/3}&\int_0^\infty x^{s-1}\left(\frac{\mathrm{Ai}'(2^{1/3}x)}{\mathrm{Ai}(2^{1/3}x)} - \frac{\mathrm{Ai}'(0)}{\mathrm{Ai}(0)}\right)\mathrm{d}x \\
&= \int_0^\infty\int_0^\infty x^{s-1}\bigl(e^{-xt}-1\bigr)\psi_{\mathrm{ex}}(t^{3/2})\frac{\mathrm{d}t}{\sqrt{2\pi t^3}}\,\mathrm{d}x \\
&= \Gamma(s)\int_0^\infty t^{-s}\psi_{\mathrm{ex}}(t^{3/2})\frac{\mathrm{d}t}{\sqrt{2\pi t^3}} \\
&= \frac{\Gamma(s)}{\sqrt{2\pi}}\cdot\frac{2}{3}\int_0^\infty u^{-2s/3-1/3}\psi_{\mathrm{ex}}(u)\frac{\mathrm{d}u}{u} \\
&= \frac{2\Gamma(s)}{3\sqrt{2\pi}}\Gamma\Bigl(-\frac{2s+1}{3}\Bigr)\,\mathbb{E}\,\mathcal{B}_{\mathrm{ex}}^{(2s+1)/3}.
\end{aligned} \tag{127}$$

By (122), this equals

$$2^{1/3-s/3}\Gamma(s)\Gamma(1-s)\Lambda(1-s) = 2^{(1-s)/3}\frac{\pi}{\sin(\pi s)}\Lambda(1-s), \tag{128}$$

as shown directly by Flajolet and Louchard [17].

## 19. Integrals of powers of $B_{\mathrm{ex}}$

Several related results are known for other functionals of a Brownian excursion or other variants of a Brownian motion. We describe some of them in this and the following sections, emphasizing the similarities with the results above, and in particular linear and quadratic recurrencies for moments and related formulas for generating functions.

The results above for moments of $\mathcal{B}_{\mathrm{ex}} = \int_0^1 B_{\mathrm{ex}}(t)\,\mathrm{d}t$ have been generalized to joint moments of the integrals $\int_0^1 B_{\mathrm{ex}}(t)^\ell\,\mathrm{d}t$, $\ell = 1, 2, \dots$, by Nguyên Thê [37] (mainly $\ell = 1$ and 2) and Richard [41] (all $\ell \geq 1$). These papers show that (6) and (7) (or equivalent formulas above) extend to these joint moments, with a quadratic recursion with multiple indices. We refer to these papers for details, and remark only that if we specialize to $\mathcal{B}_{\mathrm{ex}}$ ($\ell = 1$), then $A_{r,s}$ in Nguyên Thê [37] specializes to

$$A_{r,0} = 2^r K_r = \sigma_{r-1}, \qquad r \geq 0; \tag{129}$$

similarly, in Richard [41],

$$f_{k,0,\dots,0} = 2^{3-2k} K_k, \qquad k \geq 0, \tag{130}$$

$$m_{k,0,\dots,0} = 2^{3k/2}\,\mathbb{E}\,\mathcal{B}_{\mathrm{ex}}^k = \mathbb{E}\,\mathcal{A}^k, \qquad k \geq 0, \tag{131}$$

so $X_{1,n} \xrightarrow{\mathrm{d}} \mathcal{A}$, and $\eta_{k,0,\dots,0} = \mathbb{E}\,\mathcal{B}_{\mathrm{ex}}^k$. The recursions in both [37] and [41] specialize to (7) (or an equivalent quadratic recursion).

Janson [21] has given further similar results, including a similar quadratic recursion with double indices, for (joint) moments of $\mathcal{B}_{\mathrm{ex}}$ and another functional of a Brownian excursion; see also Chassaing and Janson [23].

## 20. Brownian bridge

Let $\mathcal{B}_{\mathrm{br}} := \int_0^1 |B_{\mathrm{br}}(t)|\,\mathrm{d}t$, the integral (or average) of the absolute value of a Brownian bridge. (Shepp [42] and Rice [40] use $\xi$ for $\mathcal{B}_{\mathrm{br}}$; Johnson and Killeen [28] use $L$; Takács [47] uses $\sigma$. Cifarelli [9] uses $C(x)$ for the distribution function of $2^{3/2}\mathcal{B}_{\mathrm{br}}$.) Comtet, Desbois and Texier [10] give also the alternative representation

$$\mathcal{B}_{\mathrm{br}} := \int_0^1 |B_{\mathrm{br}}(t)|\,\mathrm{d}t \overset{\mathrm{d}}{=} \int_0^1 |B_{\mathrm{br}}(t) - B_{\mathrm{br}}(1-t)|\,\mathrm{d}t, \tag{132}$$

which follows because $B_{\mathrm{br}}(t) - B_{\mathrm{br}}(1-t) \overset{\mathrm{d}}{=} B_{\mathrm{br}}(2t)$, $t \in [0, 1/2]$.

Cifarelli [9] and (independently) Shepp [42] found a formula for the Laplace transform $\psi_{\mathrm{br}}(t) := \mathbb{E}\,e^{-t\mathcal{B}_{\mathrm{br}}}$. (See Appendix C.2 for a proof.) An equivalent result, in physical formulation, was proved by path integral methods by Altshuler, Aronov and Khmelnitsky [3], see Comtet, Desbois and Texier [10]. Shepp's version is (with $\psi_{\mathrm{br}}$ denoted $\phi$ by Shepp), cf. (81):

$$\int_0^\infty e^{-xt}\psi_{\mathrm{br}}(\sqrt{2}\,t^{3/2})\frac{\mathrm{d}t}{\sqrt{t}} = -\sqrt{\pi}\frac{\mathrm{Ai}(x)}{\mathrm{Ai}'(x)}, \qquad x \geq 0, \tag{133}$$

while Cifarelli's version, which actually is stated using the Bessel function $K_{1/3}$, cf. (74), can be written

$$\int_0^\infty e^{-xt}\psi_{\mathrm{br}}(ut^{3/2})\frac{\mathrm{d}t}{\sqrt{2\pi t}} = -(2u)^{-1/3}\frac{\mathrm{Ai}(2^{1/3}u^{-2/3}x)}{\mathrm{Ai}'(2^{1/3}u^{-2/3}x)}, \qquad x \ge 0, \tag{134}$$

for arbitrary $u > 0$. Clearly, (133) is the special case $u = \sqrt{2}$ of (134), and the two formulas are equivalent by a simple change of variables. The moments of $\mathcal{B}_{\mathrm{br}}$ the are obtained by asymptotic expansions, cf. Section 16; Cifarelli [9] considers (134) as $u \to 0$ and Shepp [42] considers (133) as $x \to \infty$; these are obviously equivalent. Introduce $D_n$ defined by

$$\mathbb{E}\,\mathcal{B}_{\mathrm{br}}^n = \frac{\sqrt{\pi}\,2^{-n/2}n!}{\Gamma((3n+1)/2)}D_n, \qquad n \ge 0. \tag{135}$$

(Takács [47] denotes $\mathbb{E}\,\mathcal{B}_{\mathrm{br}}^n$ by $M_n^*$; Perman and Wellner [38] use $L_n$ for our $D_n$, $A_0$ for $\mathcal{B}_{\mathrm{br}}$ and $\mu_n$ for $\mathbb{E}\,\mathcal{B}_{\mathrm{br}}^n$; Nguyên Thê [36] denotes $\mathbb{E}\,\mathcal{B}_{\mathrm{me}}^n$ by $M_n^P$ and uses $Q_n$ for our $D_n$ and $a_n^P$ for $2^{-1/2-n/2}D_n$.) We then obtain from (133), cf. (98) and the discussion after it,

$$-\frac{\mathrm{Ai}(x)}{\mathrm{Ai}'(x)} \sim \sum_{k=0}^{\infty}(-1)^k D_k x^{-3k/2-1/2}, \qquad x \to \infty, \tag{136}$$

and thus by (102), (103), (106), (107), cf. (108), as formal power series,

$$\sum_{k=0}^{\infty} D_k z^k = \frac{\sum_{k=0}^{\infty}(3/2)^k c_k' z^k}{\sum_{k=0}^{\infty}(3/2)^k d_k' z^k} = \frac{\sum_{k=0}^{\infty}\alpha_k z^k}{\sum_{k=0}^{\infty}\alpha_k' z^k}. \tag{137}$$

Cifarelli [9] uses

$$m_k := 2^{3k/2}\,\mathbb{E}\,\mathcal{B}_{\mathrm{br}}^k = \mathbb{E}\big(2^{3/2}\mathcal{B}_{\mathrm{br}}\big)^k \tag{138}$$

and writes the result of the expansion as

$$m_k = \frac{\sqrt{\pi}\,k!}{\Gamma\big((3k+1)/2\big)}(-1)^k C_k, \qquad k \ge 0, \tag{139}$$

with, as formal power series (we use $a_k'$ and $b_k'$ for $a_k$ and $b_k$ in [9])

$$\sum_{k=0}^{\infty} C_k z^k = \frac{\sum_{k=0}^{\infty} a_k' z^k}{\sum_{k=0}^{\infty} b_k' z^k}, \quad \text{where } a_k' = (-3)^k c_k' \text{ and } b_k' = (-3)^k d_k'; \tag{140}$$

this is by (138) equivalent to (135) and (137) together with

$$C_k = (-2)^k D_k, \qquad k \ge 0. \tag{141}$$

Cifarelli [9] gives, instead of a recursion relation, the solution to (140) by the formula

$$(-1)^k C_k = \sum_{p=0}^{k}(-1)^{k-p}a_{k-p}'|B|_p = \sum_{p=0}^{k}|a_{k-p}'||B|_p, \tag{142}$$

where $|B|_p$ is given by the general determinant formula

$$|B|_k := \begin{vmatrix} b_1' & b_2' & \dots & b_{k-1}' & b_k' \\ 1 & b_1' & \dots & b_{k-2}' & b_{k-1}' \\ 0 & 1 & \dots & b_{k-3}' & b_{k-2}' \\ \vdots & \vdots & \ddots & \vdots & \vdots \\ 0 & 0 & \dots & 1 & b_1' \end{vmatrix}; \tag{143}$$

this is a general way (when $b_0' = 1$) to write the solution to

$$\sum_{k=0}^{\infty} |B|_k x^k = \left( \sum_{k=0}^{\infty} (-1)^k b_k' x^k \right)^{-1}. \tag{144}$$

Shepp [42] gives, omitting intermediate steps, the result of expanding (133) as (using $\xi$ for $\mathcal{B}_{\mathrm{br}}$)

$$\mathbb{E}\,\mathcal{B}_{\mathrm{br}}^n = \frac{\sqrt{\pi}}{2^{5n/2}3^{2n}\Gamma\bigl((3n+1)/2\bigr)}\bar{e}_n, \qquad n \ge 0, \tag{145}$$

where $\bar{e}_n$ satisfies the linear recursion

$$\bar{e}_n = \frac{\Gamma(3n+\frac12)}{\Gamma(n+\frac12)} + \sum_{k=1}^{n} \bar{e}_{n-k} \binom{n}{k} \frac{6k+1}{6k-1} \frac{\Gamma(3k+\frac12)}{\Gamma(k+\frac12)}, \qquad n \ge 0. \tag{146}$$

(Shepp writes $e_n$ but we use $\bar{e}_n$ to avoid confusion with $e_r$ in (64).)

We have, e.g. by (145), (135) and (141),

$$\bar{e}_n = 36^n n!\, D_n = (-18)^n n!\, C_n, \qquad n \ge 0. \tag{147}$$

Hence, (146) can, as in Perman and Wellner [38] (where further, as said above, $\alpha_n$ is denoted by $\gamma_n$), be rewritten as

$$D_n = \alpha_n - \sum_{i=1}^{n} \alpha_i' D_{n-i} = \alpha_n + \sum_{i=1}^{n} \frac{6i+1}{6i-1} \alpha_i D_{n-i}, \qquad n \ge 0, \tag{148}$$

which also follows from (137).

Takács [47] found, by different methods, a quadratic recurrence analogous to (6) and (7), namely (135) together with $D_0 = 1$ and

$$D_n = \frac{3n-2}{4} D_{n-1} - \frac{1}{2} \sum_{i=1}^{n-1} D_i D_{n-i}, \qquad n \ge 1. \tag{149}$$

This follows also by differentiating (136), which yields, as $x \to \infty$,

$$\sum_{k=0}^{\infty} (-1)^{k+1} \frac{3k+1}{2} D_k x^{-3k/2-3/2} \sim -1 + \frac{\mathrm{Ai}(x)\mathrm{Ai}''(x)}{\mathrm{Ai}'(x)^2} = -1 + x\Bigl(\frac{\mathrm{Ai}(x)}{\mathrm{Ai}'(x)}\Bigr)^2$$

$$\sim -1 + \left(\sum_{k=0}^{\infty}(-1)^k D_k x^{-3k/2}\right)^2 \tag{150}$$

and thus

$$1 + \sum_{k=0}^{\infty}\frac{3k+1}{2}D_k z^{k+1} = \left(\sum_{k=0}^{\infty}D_k z^k\right)^2. \tag{151}$$

Note further that (101) and (136) show that

$$\sum_{k=0}^{\infty}D_k z^k = \left(-2\sum_{k=0}^{\infty}K_k z^k\right)^{-1} = \frac{K_0}{\sum_{k=0}^{\infty}K_k z^k} \tag{152}$$

and thus

$$\sum_{i=0}^{k}D_i K_{k-i} = 0, \qquad k \geq 1, \tag{153}$$

which can be regarded as a linear recursion for $D_n$ given $K_n$, or conversely. Explicitly,

$$D_n = 2\sum_{i=0}^{n-1}D_i K_{n-i}, \qquad n \geq 1. \tag{154}$$

It follows easily from the asymptotics (58) that the term with $i = 0$ in (154) asymptotically dominates the sum of the others and thus

$$D_n \sim 2K_n \sim \frac{1}{\pi}\left(\frac{3}{4}\right)^n (n-1)! \sim \sqrt{\frac{2}{\pi n}}\left(\frac{3n}{4e}\right)^n \qquad \text{as } n \to \infty, \tag{155}$$

as given by Takács [47] (derived by him from (149) without further details). Consequently, cf. (53), by (135), or alternatively by Tolmatz [51] or Janson and Louchard [25],

$$\mathbb{E}\,\mathcal{B}_{\mathrm{br}}^n \sim \sqrt{2}\left(\frac{n}{12e}\right)^{n/2} \qquad \text{as } n \to \infty. \tag{156}$$

The entire function $\mathrm{Ai}'$ has all its zeros on the negative real axis (just as Ai, see Section 12), and we denote them by $a'_j = -|a'_j|$, $j = 1, 2, \dots$, with $0 < |a'_1| < |a'_2| < \dots$; thus, cf. (77),

$$\mathrm{Ai}'(a'_j) = \mathrm{Ai}'(-|a'_j|) = 0, \qquad j = 1, 2, \dots \tag{157}$$

(Abramowitz and Stegun [1, 10.4.95], Rice [40] and Johnson and Killeen [28] use $a'_j$ as we do; Takács [47; 49] uses $a'_j$ for our $|a'_j|$; Flajolet and Louchard [17] use $\alpha_j$ for $|a'_j|$; Kac [29] and Johnson and Killeen [28] use $\delta_j$ for $2^{-1/3}|a'_j|$.)

The meromorphic function $\mathrm{Ai}(z)/\mathrm{Ai}'(z)$ has the residue $\mathrm{Ai}(a'_j)/\mathrm{Ai}''(a'_j) = 1/a'_j$ at $a'_j$, and there is a convergent partial fraction expansion

$$\frac{\mathrm{Ai}(z)}{\mathrm{Ai}'(z)} = \sum_{j=1}^{\infty}\frac{1}{a'_j(z-a'_j)} = -\sum_{j=1}^{\infty}\frac{1}{|a'_j|(z+|a'_j|)}, \qquad z \in \mathbb{C}. \tag{158}$$

Inversion of the Laplace transform in (133) yields, as found by Rice [40],

$$\psi_{\mathrm{br}}(t) = 2^{-1/6}\pi^{1/2}\sum_{j=1}^{\infty}|a_j'|^{-1}t^{1/3}e^{-2^{-1/3}|a_j'|t^{2/3}}, \qquad t>0. \tag{159}$$

Johnson and Killeen [28] found by a second Laplace transform inversion the distribution function $F_{\mathrm{br}}(x) := \mathbb{P}(\mathcal{B}_{\mathrm{br}} \le x)$, here written as in Takács [47] with $u_j := |a_j'|^3/(27x^2)$,

$$\begin{aligned} F_{\mathrm{br}}(x) &= 18^{-1/6}\pi^{1/2}x^{-1}\sum_{j=1}^{\infty}e^{-u_j}u_j^{-1/3}\mathrm{Ai}\big((3u_j/2)^{2/3}\big) \\ &= 2^{-1/6}3^{2/3}\pi^{1/2}x^{-1/3}\sum_{j=1}^{\infty}\frac{1}{|a_j'|}e^{-|a_j'|^3/(27x^2)}\mathrm{Ai}\big(18^{-2/3}|a_j'|^2x^{-4/3}\big). \end{aligned} \tag{160}$$

Numerical values are given by Johnson and Killeen [28] and Takács [47].

The density function $f_{\mathrm{br}} = F_{\mathrm{br}}'$ can be obtained by termwise differentiation of (160), although no-one seems to have bothered to write it out explicitly; numerical values (obtained from (159)) are given by Rice [40].

Asymptotic expansions of $f_{\mathrm{br}}(x)$ and $F_{\mathrm{br}}(x)$ as $x\to 0$ follow from (160) and (102); only the term with $j=1$ in (160) is significant and the first terms of the asymptotic expansions are, as $x \to 0$,

$$f_{\mathrm{br}}(x) \sim e^{-\frac{2|a_1'|^3}{27x^2}}\left(\tfrac{2}{9}\,|a_1'|^{3/2}x^{-3} - \tfrac{5}{12}\,|a_1'|^{-3/2}x^{-1} + \tfrac{25}{64}\,|a_1'|^{-9/2}x + O\left(x^3\right)\right), \tag{161}$$

$$F_{\mathrm{br}}(x) \sim e^{-\frac{2|a_1'|^3}{27x^2}}\left(\tfrac{3}{2}\,|a_1'|^{-3/2} - \tfrac{45}{16}\,|a_1'|^{-9/2}x^2 + \tfrac{10395}{256}\,|a_1'|^{-15/2}x^4 + O\left(x^6\right)\right). \tag{162}$$

For $x\to\infty$ we have by Tolmatz [51] and Janson and Louchard [25] asymptotic expansions beginning with

$$f_{\mathrm{br}}(x) \sim \frac{2\sqrt{6}}{\sqrt{\pi}}\,e^{-6x^2}\left(1+\frac{1}{18}\,x^{-2}+\frac{1}{432}\,x^{-4}+O\left(x^{-6}\right)\right), \tag{163}$$

$$1-F_{\mathrm{br}}(x) \sim \frac{1}{\sqrt{6\pi}}\,x^{-1}e^{-6x^2}\left(1-\frac{1}{36}\,x^{-2}+\frac{1}{108}\,x^{-4}+O\left(x^{-6}\right)\right). \tag{164}$$

Nguyên Thê [37] extended some of these results to the joint Laplace transform and moments of $\mathcal{B}_{\mathrm{br}}$ and $\int_0^1 |B_{\mathrm{br}}|^2\,\mathrm{d}t$. Specializing to $\mathcal{B}_{\mathrm{br}}$ only, his moment formula is

$$\mathbb{E}\,\mathcal{B}_{\mathrm{br}}^r = \beta_{r,0} := \frac{\sqrt{\pi}\,2^{-3r/2}r!}{\Gamma((3k+1)/2)}B_{r,0}, \qquad r\ge 0, \tag{165}$$

with (correcting several typos in [37]) $B_{0,0}=1$ and

$$B_{r,0} = \frac{3r-2}{2}B_{r-1,0} - \tfrac{1}{2}\sum_{i=1}^{k-1}B_{i,0}B_{r-i,0}, \qquad r\ge 1. \tag{166}$$

| | | | | |
|---|---|---|---|---|
| $\mathbb{E}\,\mathcal{B}_{\mathrm{br}}^0 = 1$ | $\mathbb{E}\,\mathcal{B}_{\mathrm{br}} = \frac{1}{4}\sqrt{\frac{\pi}{2}}$ | $\mathbb{E}\,\mathcal{B}_{\mathrm{br}}^2 = \frac{7}{60}$ | $\mathbb{E}\,\mathcal{B}_{\mathrm{br}}^3 = \frac{21}{512}\sqrt{\frac{\pi}{2}}$ | $\mathbb{E}\,\mathcal{B}_{\mathrm{br}}^4 = \frac{19}{720}$ |
| $m_0 = 1$ | $m_1 = \frac{1}{2}\sqrt{\pi}$ | $m_2 = \frac{14}{15}$ | $m_3 = \frac{21}{32}\sqrt{\pi}$ | $m_4 = \frac{76}{45}$ |
| $\bar{e}_0 = 1$ | $\bar{e}_1 = 9$ | $\bar{e}_2 = 567$ | $\bar{e}_3 = 91854$ | $\bar{e}_4 = 28796229$ |
| $D_0 = 1$ | $D_1 = \frac{1}{4}$ | $D_2 = \frac{7}{32}$ | $D_3 = \frac{21}{64}$ | $D_4 = \frac{1463}{2048}$ |
| $C_0 = 1$ | $C_1 = -\frac{1}{2}$ | $C_2 = \frac{7}{8}$ | $C_3 = -\frac{21}{8}$ | $C_4 = \frac{1463}{128}$ |
| $\bar{e}_0^* = 1$ | $\bar{e}_1^* = 1$ | $\bar{e}_2^* = 7$ | $\bar{e}_3^* = 126$ | $\bar{e}_4^* = 4389$ |
| $D_0^* = 1$ | $D_1^* = 2$ | $D_2^* = 14$ | $D_3^* = 168$ | $D_4^* = 2926$ |

TABLE 2
*Some numerical values for the Brownian bridge.*

Clearly, these results are equivalent to (135) and (149), with

$$B_{r,0} = 2^r D_r = (-1)^r C_r, \qquad r \ge 0. \tag{167}$$

Nguyên Thê [37] gives also the relation

$$\sum_{i=0}^{r} A_{i,0} B_{r-i,0} = 0, \qquad r \ge 1, \tag{168}$$

which by (129) and (167) is equivalent to (153). Equivalently,

$$\sum_{k=0}^{\infty} A_{k,0} y^k = \frac{A_{0,0}}{\sum_{k=0}^{\infty} B_{k,0} y^k}. \tag{169}$$

Some numerical values are given in Table 2; see further Cifarelli [9] ($m_n$), Shepp [42] ($\mathbb{E}\,\mathcal{B}_{\mathrm{br}}^n$, $\bar{e}_n$), Takács [47] ($\mathbb{E}\,\mathcal{B}_{\mathrm{br}}^n$, $D_n$), Nguyên Thê [37] ($\mathbb{E}\,\mathcal{B}_{\mathrm{br}}^n$, $B_{n,m}$). We define $D_n^* := 2^{3n} D_n$ and $\bar{e}_n^* := 9^{-n}\bar{e}_n = 4^n n!\, D_n$ (considered by Shepp [42]); note that these are integers.

## 21. Brownian motion

Let $\mathcal{B}_{\mathrm{bm}} := \int_0^1 |B(t)|\,\mathrm{d}t$, the integral of the absolute value of a Brownian motion on the unit interval, and let $\psi_{\mathrm{bm}}(t) := \mathbb{E}\,e^{-t\mathcal{B}_{\mathrm{bm}}}$ be its Laplace transform. Kac [29] gave the formula

$$\psi_{\mathrm{bm}}(t) = \sum_{j=1}^{\infty} \kappa_j e^{-\delta_j t^{2/3}}, \qquad t > 0, \tag{170}$$

where $\delta_j = 2^{-1/3}|a_j'|$ and, with Kac's notation $P(y) = 2^{1/3}\mathrm{Ai}\bigl(-2^{1/3}y\bigr)$,

$$\kappa_j := \frac{1 + 3\int_0^{\delta_j} P(y)\,\mathrm{d}y}{3\delta_j P(\delta_j)} = \frac{1 + 3\int_{a_j'}^{0} \mathrm{Ai}(x)\,\mathrm{d}x}{3|a_j'|\mathrm{Ai}(a_j')}. \tag{171}$$

(Takács [49] uses $C_j$ for $\kappa_j$.) Using the function AI defined in (78), (171) can be written

$$\kappa_j = \frac{\mathrm{AI}(a_j')}{|a_j'|\mathrm{Ai}(a_j')}. \tag{172}$$

In fact, see Takács [49] and Perman and Wellner [38], or (362) in Appendix C.1 below,

$$\int_0^\infty e^{-xt}\psi_{\mathrm{bm}}\big(\sqrt{2}\,t^{3/2}\big)\,\mathrm{d}t = -\frac{\mathrm{AI}(x)}{\mathrm{Ai}'(x)}, \qquad x>0. \tag{173}$$

The meromorphic function $-\mathrm{AI}(z)/\mathrm{Ai}'(z)$ has residue $-\mathrm{AI}(a_j')/\mathrm{Ai}''(a_j') = \kappa_j$ at $a_j'$, and the partial fraction expansion

$$-\frac{\mathrm{AI}(z)}{\mathrm{Ai}'(z)} = \sum_{j=1}^\infty \frac{\kappa_j}{z-a_j'} = \sum_{j=1}^\infty \frac{\kappa_j}{z+|a_j'|}, \qquad z\in\mathbb{C}. \tag{174}$$

An inversion of the Laplace transform in (173) thus yields (170).

Takács [49] found by a second Laplace transform inversion, this time of (170), the density function $f_{\mathrm{bm}}(x)$ (by him denoted $h(x)$) of $\mathcal{B}_{\mathrm{bm}}$ as

$$f_{\mathrm{bm}}(x) = \sum_{j=1}^\infty \kappa_j\sqrt{2}|a_j'|^{-3/2}p\big(\sqrt{2}|a_j'|^{-3/2}x\big) \tag{175}$$

$$= \frac{\sqrt{3}}{\sqrt{\pi}\,x}\sum_{j=1}^\infty \kappa_j(v_j')^{2/3}e^{-v_j'}U\Big(\frac{1}{6},\frac{4}{3};v_j'\Big), \tag{176}$$

where $p$ is the stable density in (87), and $v_j' = 2|a_j'|^3/27x^2$ (denoted $v_j$ by Takács [49]), cf. (85), (91), (92).

Asymptotic expansions of $f_{\mathrm{bm}}(x)$ and the distribution function $F_{\mathrm{bm}}(x)$ as $x\to 0$ follow from (176) using the asymptotic expansion of $U$ given in e.g. [32, (9.12.3)], cf. Section 15; only the term with $j=1$ in (176) is significant and the first terms of the asymptotic expansions are, as $x\to 0$,

$$f_{\mathrm{bm}}(x) \sim \frac{\kappa_1}{\sqrt{2\pi}}e^{-\frac{2|a_1'|^3}{27x^2}}\Big(\tfrac{2}{3}|a_1'|^{3/2}x^{-2} + \tfrac{1}{4}|a_1'|^{-3/2} - \tfrac{105}{64}|a_1'|^{-9/2}x^2 + \dots\Big), \tag{177}$$

$$F_{\mathrm{bm}}(x) \sim \frac{\kappa_1}{\sqrt{2\pi}}e^{-\frac{2|a_1'|^3}{27x^2}}\Big(\tfrac{9}{2}|a_1'|^{-3/2}x - \tfrac{459}{16}|a_1'|^{-9/2}x^3 + \tfrac{145881}{256}|a_1'|^{-15/2}x^5 + \dots\Big). \tag{178}$$

For $x\to\infty$ we have by Tolmatz [52] and Janson and Louchard [25] asymptotic expansions beginning with

$$f_{\mathrm{bm}}(x) \sim \frac{\sqrt{6}}{\sqrt{\pi}}\,e^{-3x^2/2}\left(1+\frac{1}{18}\,x^{-2} - \frac{1}{162}\,x^{-4} + \dots\right), \tag{179}$$

$$1-F_{\mathrm{bm}}(x) \sim \frac{\sqrt{2}}{\sqrt{3\pi}}\,x^{-1}e^{-3x^2/2}\left(1-\frac{5}{18}\,x^{-2} + \frac{22}{81}\,x^{-4} + \dots\right). \tag{180}$$

Takács [49] also found recursion formulas for the moments. Define $L_n$ by

$$\mathbb{E}\,\mathcal{B}_{\mathrm{bm}}^n = \frac{2^{-n/2} n!}{\Gamma((3n+2)/2)} L_n, \qquad n \ge 0. \tag{181}$$

(Takács [49] denotes $\mathbb{E}\,\mathcal{B}_{\mathrm{bm}}^n$ by $\mu_n$; Perman and Wellner [38] use $K_n$ for our $L_n$, $A$ for $\mathcal{B}_{\mathrm{bm}}$ and $\nu_n$ for $\mathbb{E}\,\mathcal{B}_{\mathrm{bm}}^n$; Nguyên Thê [36] uses $a_n^B$ for $2^{-n/2}L_n$.) An asymptotic expansion of the left hand side of (173) yields, arguing as for (98) and (136),

$$-\frac{\mathrm{AI}(x)}{\mathrm{Ai}'(x)} \sim \sum_{k=0}^{\infty} (-1)^k L_k x^{-3k/2-1}, \qquad x \to \infty. \tag{182}$$

Recall the asymptotic expansions (102) and (103) of $\mathrm{Ai}(x)$ and $\mathrm{Ai}'(x)$; there is a similar expansion of $\mathrm{AI}(x)$ [49], see Appendix A,

$$\mathrm{AI}(x) \sim \tfrac{1}{2\sqrt{\pi}} x^{-3/4} e^{-2x^{3/2}/3} \sum_{k=0}^{\infty} (-1)^k \beta_k x^{-3k/2}, \qquad x \to \infty, \tag{183}$$

where $\beta_k$, $k \ge 0$, are given by $\beta_0 = 1$ and the recursion relation (obtained from comparing a formal differentiation of (183) with (102))

$$\beta_k = \alpha_k + \frac{3(2k-1)}{4} \beta_{k-1}, \qquad k \ge 1. \tag{184}$$

(Takács [49] further uses $h_k = (2/3)^k \beta_k$.) Hence, (182) yields, using (103) and (107), the equality (for formal power series)

$$\sum_{k=0}^{\infty} L_k z^k = \frac{\sum_{k=0}^{\infty} \beta_k z^k}{\sum_{k=0}^{\infty} (3/2)^k d'_k z^k} = \frac{\sum_{k=0}^{\infty} \beta_k z^k}{\sum_{k=0}^{\infty} \alpha'_k z^k}. \tag{185}$$

Multiplying by the denominator and identifying coefficients leads to the recursion by Takács [49],

$$L_n = \beta_n + \sum_{j=1}^{n} \frac{6j+1}{6j-1} \alpha_j L_{n-j}, \qquad n \ge 0. \tag{186}$$

Differentiation of (182) yields, using $\mathrm{Ai}''(x) = x\mathrm{Ai}(x)$,

$$\frac{\mathrm{Ai}(x)}{\mathrm{Ai}'(x)} + \frac{x\mathrm{Ai}(x)\mathrm{AI}(x)}{\mathrm{Ai}'(x)^2} \sim \sum_{k=0}^{\infty} (-1)^{k+1} \frac{3k+2}{2} L_k x^{-3k/2-2}, \qquad x \to \infty, \tag{187}$$

and thus by (136) and (182) the equality

$$-\sum_{k=0}^{\infty} D_k z^k + \sum_{k=0}^{\infty} D_k z^k \sum_{k=0}^{\infty} L_k z^k = \sum_{k=0}^{\infty} \frac{3k+2}{2} L_k z^{k+1} = \sum_{k=1}^{\infty} \frac{3k-1}{2} L_{k-1} z^k, \tag{188}$$

| | | | | |
|---|---|---|---|---|
| $\mathbb{E}\,\mathcal{B}_{\mathrm{bm}}^0 = 1$ | $\mathbb{E}\,\mathcal{B}_{\mathrm{bm}} = \frac{2}{3}\sqrt{\frac{2}{\pi}}$ | $\mathbb{E}\,\mathcal{B}_{\mathrm{bm}}^2 = \frac{3}{8}$ | $\mathbb{E}\,\mathcal{B}_{\mathrm{bm}}^3 = \frac{263}{630}\sqrt{\frac{2}{\pi}}$ | $\mathbb{E}\,\mathcal{B}_{\mathrm{bm}}^4 = \frac{903}{2560}$ |
| $\beta_0 = 1$ | $\beta_1 = \frac{41}{48}$ | $\beta_2 = \frac{9241}{4608}$ | $\beta_3 = \frac{5075225}{663552}$ | $\beta_4 = \frac{5153008945}{127401984}$ |
| $L_0 = 1$ | $L_1 = 1$ | $L_2 = \frac{9}{4}$ | $L_3 = \frac{263}{32}$ | $L_4 = \frac{2709}{64}$ |
| $L_0^* = 1$ | $L_1^* = 8$ | $L_2^* = 144$ | $L_3^* = 4208$ | $L_4^* = 173376$ |

TABLE 3
*Some numerical values for the Brownian motion.*

which is equivalent to $L_0 = 1$ and the recursion formula

$$\frac{3n-1}{2}L_{n-1} = \sum_{j=1}^{n} D_{n-j}L_j, \qquad n \ge 1 \tag{189}$$

or

$$L_n = \frac{3n-1}{2}L_{n-1} - \sum_{j=1}^{n-1} D_{n-j}L_j, \qquad n \ge 1. \tag{190}$$

Takács [49] proves also the asymptotic relations, see further Tolmatz [52] and Janson and Louchard [25],

$$\mathbb{E}\,\mathcal{B}_{\mathrm{bm}}^n \sim \sqrt{2}\Bigl(\frac{n}{3e}\Bigr)^{n/2}, \tag{191}$$

$$L_n \sim \sqrt{3}\Bigl(\frac{3n}{2e}\Bigr)^{n}. \tag{192}$$

Some numerical values are given in Table 3; see further Takács [49] and Perman and Wellner [38]. We define $L_n^* := 2^{3n}L_n$; these are integers by (190).

## 22. Brownian meander

Let $\mathcal{B}_{\mathrm{me}} := \int_0^1 |B_{\mathrm{me}}(t)|\,\mathrm{d}t$, the integral of a Brownian meander on the unit interval, and let $\psi_{\mathrm{me}}(t) := \mathbb{E}\,e^{-t\mathcal{B}_{\mathrm{me}}}$ be its Laplace transform.

Takács [50] and Perman and Wellner [38] give formulas equivalent to

$$\int_0^\infty e^{-xt}\psi_{\mathrm{me}}\bigl(\sqrt{2}\,t^{3/2}\bigr)\frac{\mathrm{d}t}{\sqrt{\pi t}} = \frac{\mathrm{AI}(x)}{\mathrm{Ai}(x)}, \qquad x > 0; \tag{193}$$

see also Appendix C.3.

Define $Q_n$ by

$$\mathbb{E}\,\mathcal{B}_{\mathrm{me}}^n = \frac{\sqrt{\pi}2^{-n/2}n!}{\Gamma((3n+1)/2)}Q_n, \qquad n \ge 0. \tag{194}$$

(Takács [50] denotes $\mathbb{E}\,\mathcal{B}_{\mathrm{me}}^n$ by $M_n$ and $\psi_{\mathrm{me}}$ by $\Psi$; Perman and Wellner [38] use $R_n$ for our $Q_n$, $A_{\mathrm{mean}}$ for $\mathcal{B}_{\mathrm{me}}$ and $\rho_n$ for $\mathbb{E}\,\mathcal{B}_{\mathrm{me}}^n$; Nguyên Thê [36] denotes $\mathbb{E}\,\mathcal{B}_{\mathrm{me}}^n$

by $M_n^F$ and uses $a_n^F$ for $2^{1/2-n/2}Q_n$; Majumdar and Comtet [35] denote $\mathbb{E}\,\mathcal{B}_{\mathrm{me}}^n$ by $a_n$ and $Q_n$ by $R_n$.) An asymptotic expansion of the left hand side of (193) yields, arguing as for (98), (136) and (182),

$$\frac{\mathrm{AI}(x)}{\mathrm{Ai}(x)} \sim \sum_{k=0}^{\infty}(-1)^k Q_k x^{-3k/2-1/2}, \qquad x \to \infty. \tag{195}$$

The asymptotic expansions (102) and (183) of $\mathrm{Ai}(x)$ and $\mathrm{AI}(x)$ yield, using (106), the equality (for formal power series)

$$\sum_{k=0}^{\infty} Q_k z^k = \frac{\sum_{k=0}^{\infty}\beta_k z^k}{\sum_{k=0}^{\infty}(3/2)^k c_k' z^k} = \frac{\sum_{k=0}^{\infty}\beta_k z^k}{\sum_{k=0}^{\infty}\alpha_k z^k}. \tag{196}$$

Differentiation of (195) yields

$$-1-\frac{\mathrm{Ai}'(x)\mathrm{AI}(x)}{\mathrm{Ai}(x)^2} \sim \sum_{k=0}^{\infty}(-1)^{k+1}\frac{3k+1}{2}Q_k x^{-3k/2-3/2}, \qquad x \to \infty, \tag{197}$$

and thus, using (101), the equality

$$\sum_{k=0}^{\infty}\frac{3k+1}{2}Q_k z^{k+1} = -1-2\sum_{k=0}^{\infty}K_k z^k \sum_{k=0}^{\infty}Q_k z^k, \tag{198}$$

which is equivalent to the recursion formula

$$Q_n = \frac{3n-2}{2}Q_{n-1} + 2\sum_{j=1}^{n}K_j Q_{n-j}, \qquad n \ge 1. \tag{199}$$

with $Q_0 = 1$, proved by Takács [50] by a different method, and the differential equation

$$3z^2\frac{\mathrm{d}}{\mathrm{d}z}\sum_{k=0}^{\infty}Q_k z^k + \Bigl(z + 4\sum_{k=0}^{\infty}K_k z^k\Bigr)\sum_{k=0}^{\infty}Q_k z^k = -2, \tag{200}$$

also given by Takács [50] (in a slightly different form, with $-z$ substituted for $z$). On the other hand, multiplying by $\sum_{k=0}^{\infty}\alpha_k z^k$ in (196) and identifying coefficients leads to another recursion by Takács [50] and Perman and Wellner [38]:

$$Q_n = \beta_n - \sum_{j=1}^{n}\alpha_j Q_{n-j}, \qquad n \ge 0. \tag{201}$$

Furthermore, (108), (185) and (196) yield the equality

$$-2\sum_{k=0}^{\infty}K_k z^k \sum_{k=0}^{\infty}L_k z^k = \sum_{k=0}^{\infty}Q_k z^k, \tag{202}$$

| | | | | |
|---|---|---|---|---|
| $\mathbb{E}\,\mathcal{B}_{\mathrm{me}}^0 = 1$ | $\mathbb{E}\,\mathcal{B}_{\mathrm{me}} = \frac{3}{4}\sqrt{\frac{\pi}{2}}$ | $\mathbb{E}\,\mathcal{B}_{\mathrm{me}}^2 = \frac{59}{60}$ | $\mathbb{E}\,\mathcal{B}_{\mathrm{me}}^3 = \frac{465}{512}\sqrt{\frac{\pi}{2}}$ | $\mathbb{E}\,\mathcal{B}_{\mathrm{me}}^4 = \frac{5345}{3696}$ |
| $Q_0 = 1$ | $Q_1 = \frac{3}{4}$ | $Q_2 = \frac{59}{32}$ | $Q_3 = \frac{465}{64}$ | $Q_4 = \frac{80175}{2048}$ |
| $Q_0^* = 1$ | $Q_1^* = 6$ | $Q_2^* = 118$ | $Q_3^* = 3720$ | $Q_4^* = 160350$ |

TABLE 4
*Some numerical values for the Brownian meander.*

and thus the relation

$$Q_n = -2\sum_{j=0}^{n} K_j L_{n-j}, \qquad n \geq 0. \tag{203}$$

By (152) we obtain from (202) also

$$\sum_{k=0}^{\infty} D_k z^k \sum_{k=0}^{\infty} Q_k z^k = \sum_{k=0}^{\infty} L_k z^k, \tag{204}$$

and thus the relation

$$L_n = \sum_{j=0}^{n} D_j Q_{n-j}, \qquad n \geq 0. \tag{205}$$

Takács [50] proves the asymptotic relations, see also Janson and Louchard [25],

$$\mathbb{E}\,\mathcal{B}_{\mathrm{me}}^n \sim \sqrt{3\pi n}\Big(\frac{n}{3e}\Big)^{n/2}, \tag{206}$$

$$Q_n \sim \sqrt{3}\Big(\frac{3n}{2e}\Big)^{n}. \tag{207}$$

Some numerical values are given in Table 4; see further Takács [50]. We define $Q_n^* := 2^{3n} Q_n$; these are integers by (203) (or by (199) and (12)).

The meromorphic function $\mathrm{AI}(z)/\mathrm{Ai}(z)$ has residue $r_j := \mathrm{AI}(a_j)/\mathrm{Ai}'(a_j)$ at $a_j$, and the partial fraction expansion

$$\frac{\mathrm{AI}(z)}{\mathrm{Ai}(z)} = \sum_{j=1}^{\infty} \frac{r_j}{z+|a_j|} = \sum_{j=1}^{\infty} r_j \Big(\frac{1}{z+|a_j|} - \frac{1}{|a_j|}\Big) + \frac{1}{3\mathrm{Ai}(0)}, \qquad z \in \mathbb{C}, \tag{208}$$

where the first sum is conditionally convergent only and the second sum is absolutely convergent, since $|a_j| \asymp j^{2/3}$ and $|r_j| \asymp j^{-1/6}$. (Note that $r_j$ alternates in sign.) In fact, the second version follows from an absolutely convergent partial fraction expansion of $\mathrm{AI}(z)/(z\mathrm{Ai}(z))$, using $\mathrm{AI}(0) = 1/3$, and the first then follows easily. A Laplace transform inversion (considering, for example, the Laplace transform of $\sqrt{t}\,\psi_{\mathrm{me}}(\sqrt{2}\,t^{3/2})$) yields the formula by Takács [50]

$$\psi_{\mathrm{me}}(t) = 2^{-1/6} t^{1/3} \sqrt{\pi} \sum_{j=1}^{\infty} r_j e^{-2^{-1/3}|a_j| t^{2/3}}, \qquad t > 0. \tag{209}$$

A second Laplace inversion, see Takács [50], yields the distribution function $F_{\mathrm{me}}(x)$ of $\mathcal{B}_{\mathrm{me}}$ as, using $u_j = |a_j|^3/(27x^2)$ and $R_j = |a_j|r_j = |a_j|\mathrm{AI}(a_j)/\mathrm{Ai}'(a_j)$,

$$\begin{aligned} F_{\mathrm{me}}(x) &= 18^{-1/6}\pi^{1/2}x^{-1}\sum_{j=1}^{\infty} R_j e^{-u_j} u_j^{-1/3}\mathrm{Ai}\big((3u_j/2)^{2/3}\big) \\ &= 2^{-1/6}3^{2/3}\pi^{1/2}x^{-1/3}\sum_{j=1}^{\infty} r_j e^{-|a_j|^3/(27x^2)}\mathrm{Ai}\big(18^{-2/3}a_j^2x^{-4/3}\big). \end{aligned} \tag{210}$$

Numerical values are given by Takács [50]. The density function $f_{\mathrm{me}}(x)$ can be obtained by termwise differentiation.

A physical proof of (209) by path integral methods is given by Majumdar and Comtet [35] (denoting our $r_j$ by $B(\alpha_j)$), who then conversely derive (193) from (209).

Asymptotic expansions of $f_{\mathrm{me}}(x)$ and $F_{\mathrm{me}}(x)$ as $x \to 0$ follow from (210) and (102); only the term with $j = 1$ in (210) is significant and the first terms of the asymptotic expansions are, as $x \to 0$,

$$f_{\mathrm{me}}(x) \sim r_1 e^{-\frac{2|a_1|^3}{27x^2}}\left(\tfrac{2}{9}|a_1|^{5/2}x^{-3} - \tfrac{5}{12}|a_1|^{-1/2}x^{-1} + \tfrac{25}{64}|a_1|^{-7/2}x + O\left(x^3\right)\right), \tag{211}$$

$$F_{\mathrm{me}}(x) \sim r_1 e^{-\frac{2|a_1|^3}{27x^2}}\left(\tfrac{3}{2}|a_1|^{-1/2} - \tfrac{45}{16}|a_1|^{-7/2}x^2 + \tfrac{10395}{256}|a_1|^{-13/2}x^4 + O\left(x^6\right)\right). \tag{212}$$

For $x \to \infty$, Janson and Louchard [25] found asymptotic expansions beginning with

$$f_{\mathrm{me}}(x) \sim 3\sqrt{3}\,xe^{-3x^2/2}\left(1 - \frac{1}{18}x^{-2} - \frac{1}{162}x^{-4} + \dots\right), \tag{213}$$

$$1 - F_{\mathrm{me}}(x) \sim \sqrt{3}\,e^{-3x^2/2}\left(1 - \frac{1}{18}x^{-2} + \frac{5}{162}x^{-4} + \dots\right). \tag{214}$$

## 23. Brownian double meander

Define

$$B_{\mathrm{dm}}(t) := B(t) - \min_{0\le u\le 1} B(u), \qquad t \in [0,1]. \tag{215}$$

This is a non-negative continuous stochastic process on $[0,1]$ that a.s. is 0 at a unique point $\tau \in [0,1]$ (the time of the minimum of $B(t)$ on $[0,1]$). It is well-known that $\tau$ has an arcsine ($= \mathrm{Beta}(\frac12,\frac12)$) distribution with density $\pi^{-1}(t(1-t))^{-1/2}$; moreover, given $\tau$, the processes $B_{\mathrm{dm}}(\tau - t)$, $t \in [0,\tau]$, and $B_{\mathrm{dm}}(\tau + t)$, $t \in [0,1-\tau]$, are two independent Brownian meanders on the respective intervals, see Williams [57], Denisov [13] and Bertoin, Pitman and Ruiz de Chavez [6]. (This may be seen as a limit of an elementary corresponding result for simple random walks.) We therefore call $B_{\mathrm{dm}}$ a *Brownian double meander*.

Let $\mathcal{B}_{\mathrm{dm}} := \int_0^1 B_{\mathrm{dm}}(t)\,\mathrm{d}t$ be the Brownian double meander area. By the definition, we have the formula

$$\mathcal{B}_{\mathrm{dm}} := \int_0^1 B_{\mathrm{dm}}(t)\,\mathrm{d}t = \int_0^1 B(t)\,\mathrm{d}t - \min_{0\le t\le 1} B(t); \tag{216}$$

note the analogy with (16) for $\mathcal{B}_{\mathrm{ex}}$ and $B_{\mathrm{br}}$. As in Section 3, there are further interesting equivalent forms of this. Since $B$ is symmetric, $B \overset{\mathrm{d}}{=} -B$, we also have

$$\mathcal{B}_{\mathrm{dm}} \overset{\mathrm{d}}{=} \max_{0\le t\le 1} B(t) - \int_0^1 B(t)\,\mathrm{d}t = \max_{0\le t\le 1} B^0(t), \tag{217}$$

where

$$B^0(t) := B(t) - \int_0^1 B(s)\,\mathrm{d}s, \qquad 0 \le t \le 1, \tag{218}$$

is a continuous Gaussian process with mean 0 and integral identically zero; its covariance function is, by straightforward calculation, given by

$$\operatorname{Cov}\bigl(B^0(s), B^0(t)\bigr) = \tfrac{1}{3} - \max(s,t) + \tfrac{1}{2}(s^2+t^2), \qquad s,t \in [0,1]. \tag{219}$$

Let $\psi_{\mathrm{dm}}(t) := \mathbb{E}\, e^{-t\mathcal{B}_{\mathrm{dm}}}$ be the Laplace transform of $\mathcal{B}_{\mathrm{dm}}$. Majumdar and Comtet [35] give the formula (using $C(\alpha_j)$ for our $r_j^2$), where $r_j = \mathrm{AI}(a_j)/\mathrm{Ai}'(a_j)$ as in Section 22,

$$\psi_{\mathrm{dm}}(t) = 2^{-1/3} t^{2/3} \sum_{j=1}^{\infty} r_j^2 e^{-2^{-1/3}|a_j| t^{2/3}}, \qquad t > 0; \tag{220}$$

they then derive from it

$$\int_0^\infty e^{-xt} \psi_{\mathrm{dm}}\bigl(\sqrt{2}\, t^{3/2}\bigr)\,\mathrm{d}t = \left(\frac{\mathrm{AI}(x)}{\mathrm{Ai}(x)}\right)^2, \qquad x > 0, \tag{221}$$

using the partial fraction expansion, cf. (208),

$$\left(\frac{\mathrm{AI}(z)}{\mathrm{Ai}(z)}\right)^2 = \sum_{j=1}^{\infty} \frac{r_j^2}{(z-a_j)^2} = \sum_{j=1}^{\infty} \frac{r_j^2}{(z+|a_j|)^2}, \qquad z \in \mathbb{C}. \tag{222}$$

(There are only quadratic terms in this partial fraction expansion, because $\mathrm{AI}'(z) = -\mathrm{Ai}(z)$ and $\mathrm{Ai}''(z) = z\mathrm{Ai}(z)$ vanish at the zeros of Ai; hence there is no constant term in the expansion of $\mathrm{AI}(z)/\mathrm{Ai}(z)$ at a pole $a_j$.) Conversely, we will prove (221) by other methods in Section 27, and then (220) follows by (222) and a Laplace transform inversion.

Majumdar and Comtet [35] found, by a Laplace transform inversion of (220), the density function $f_{\mathrm{dm}}(x)$ of $\mathcal{B}_{\mathrm{dm}}$ as

$$f_{\mathrm{dm}}(x) = \frac{2^{-1/3}}{\sqrt{3\pi}} x^{-7/3} \sum_{j=1}^{\infty} |a_j| r_j^2 e^{-v_j} \Bigl(U\Bigl(\frac{1}{6}, \frac{4}{3}; v_j\Bigr) + 2U\Bigl(-\frac{5}{6}, \frac{4}{3}; v_j\Bigr)\Bigr) \tag{223}$$

where $U$ is the confluent hypergeometric function and $v_j = 2|a_j|^3/27x^2$, cf. (91) and (92).

Asymptotic expansions of $f_{\mathrm{dm}}(x)$ and the distribution function $F_{\mathrm{dm}}(x)$ as $x \to 0$ follow as in [35] from (223), again using the asymptotic expansion of $U$ [32, (9.12.3)]; only the term with $j = 1$ in (223) is significant and the first terms of the asymptotic expansions are, as $x \to 0$,

$$f_{\mathrm{dm}}(x) \sim \sqrt{\tfrac{2}{\pi}} r_1^2 e^{-\frac{2|a_1|^3}{27x^2}} \Bigl( \tfrac{2}{27}|a_1|^{7/2}x^{-4} - \tfrac{17}{36}|a_1|^{1/2}x^{-2} + \tfrac{1}{192}|a_1|^{-5/2} + \dots \Bigr), \tag{224}$$

$$F_{\mathrm{dm}}(x) \sim \sqrt{\tfrac{2}{\pi}} r_1^2 e^{-\frac{2|a_1|^3}{27x^2}} \Bigl( \tfrac{1}{2}|a_1|^{1/2}x^{-1} + \tfrac{3}{16}|a_1|^{-5/2}x - \tfrac{315}{256}|a_1|^{-11/2}x^3 + \dots \Bigr), \tag{225}$$

For $x \to \infty$, Janson and Louchard [25] found asymptotic expansions beginning with

$$f_{\mathrm{dm}}(x) \sim \frac{2\sqrt{6}}{\sqrt{\pi}} e^{-3x^2/2} \left( 1 + \frac{1}{6}x^{-2} + \frac{1}{18}x^{-4} + \dots \right), \tag{226}$$

$$1 - F_{\mathrm{dm}}(x) \sim \frac{2\sqrt{2}}{\sqrt{3\pi}} x^{-1} e^{-3x^2/2} \left( 1 - \frac{1}{6}x^{-2} + \frac{2}{9}x^{-4} + \dots \right). \tag{227}$$

The weaker statement

$$-\ln \mathbb{P}(\mathcal{B}_{\mathrm{dm}} > x) = -\ln\bigl(1 - F_{\mathrm{dm}}(x)\bigr) \sim -3x^2/2, \qquad \text{as } x \to \infty, \tag{228}$$

had earlier been proved by Majumdar and Comtet [35] using moment asymptotics, see (235) below.

Define $W_n$ by, in analogy with (181),

$$\mathbb{E}\,\mathcal{B}_{\mathrm{dm}}^n = \frac{2^{-n/2} n!}{\Gamma((3n+2)/2)} W_n, \qquad n \ge 0. \tag{229}$$

An asymptotic expansion of the left hand side of (221) yields, arguing as for (98), (136), (182) and (195),

$$\left( \frac{\mathrm{AI}(x)}{\mathrm{Ai}'(x)} \right)^2 \sim \sum_{k=0}^{\infty} (-1)^k W_k x^{-3k/2-1}, \qquad x \to \infty. \tag{230}$$

By (183), (102), (106) and (196), this yields the equality (for formal power series)

$$\sum_{k=0}^{\infty} W_k z^k = \left( \frac{\sum_{k=0}^{\infty} \beta_k z^k}{\sum_{k=0}^{\infty} (3/2)^k c_k' z^k} \right)^2 = \left( \frac{\sum_{k=0}^{\infty} \beta_k z^k}{\sum_{k=0}^{\infty} \alpha_k z^k} \right)^2 = \left( \sum_{k=0}^{\infty} Q_k z^k \right)^2. \tag{231}$$

Consequently,

$$W_n = \sum_{j=0}^{n} Q_j Q_{n-j}, \qquad n \ge 0, \tag{232}$$

| | | | | |
|---|---|---|---|---|
| $\mathbb{E}\,\mathcal{B}_{\mathrm{dm}}^0 = 1$ | $\mathbb{E}\,\mathcal{B}_{\mathrm{dm}} = \sqrt{\frac{2}{\pi}}$ | $\mathbb{E}\,\mathcal{B}_{\mathrm{dm}}^2 = \frac{17}{24}$ | $\mathbb{E}\,\mathcal{B}_{\mathrm{dm}}^3 = \frac{123}{140}\sqrt{\frac{2}{\pi}}$ | $\mathbb{E}\,\mathcal{B}_{\mathrm{dm}}^4 = \frac{2963}{3840}$ |
| $W_0 = 1$ | $W_1 = \frac{3}{2}$ | $W_2 = \frac{17}{4}$ | $W_3 = \frac{1107}{64}$ | $W_4 = \frac{2963}{32}$ |
| $W_0^* = 1$ | $W_1^* = 12$ | $W_2^* = 272$ | $W_3^* = 8856$ | $W_4^* = 379264$ |

TABLE 5
*Some numerical values for the Brownian double meander.*

or, as given by Majumdar and Comtet [35] (where $\mathbb{E}\,\mathcal{B}_{\mathrm{dm}}^n$ is denoted by $\mu_n$),

$$\mathbb{E}\,\mathcal{B}_{\mathrm{dm}}^n = \frac{1}{\pi}\sum_{m=0}^{n}\binom{n}{m}B\Big(\frac{3m+1}{2},\frac{3(n-m)+1}{2}\Big)\,\mathbb{E}\,\mathcal{B}_{\mathrm{me}}^m\,\mathbb{E}\,\mathcal{B}_{\mathrm{me}}^{n-m}. \tag{233}$$

From (232) and (207) follows the asymptotic relation

$$W_n \sim 2Q_0Q_n = 2Q_n \sim 2\sqrt{3}\Big(\frac{3n}{2e}\Big)^n, \tag{234}$$

and thus by (229), see also Janson and Louchard [25],

$$\mathbb{E}\,\mathcal{B}_{\mathrm{dm}}^n \sim 2\sqrt{2}\Big(\frac{n}{3e}\Big)^{n/2}. \tag{235}$$

Some numerical values are given in Table 5; see also Majumdar and Comtet [35]. We define $W_n^* := 2^{3n}W_n$; these are integers by (232).

## 24. Positive part of a Brownian bridge

Let $x_+ := \max(x,0)$ and $x_- := (-x)_+$ and define

$$\mathcal{B}_{\mathrm{br}}^{\pm} := \int_0^1 B_{\mathrm{br}}(t)_\pm\,\mathrm{d}t; \tag{236}$$

thus $\mathcal{B}_{\mathrm{br}}^+$ is the average of the positive part and $\mathcal{B}_{\mathrm{br}}^-$ the average of the negative part of a Brownian bridge. (Perman and Wellner [38] use $A_0^+$ for $\mathcal{B}_{\mathrm{br}}^+$.) In particular,

$$\mathcal{B}_{\mathrm{br}} = \mathcal{B}_{\mathrm{br}}^+ + \mathcal{B}_{\mathrm{br}}^-, \tag{237}$$

and the difference is Gaussian:

$$\mathcal{B}_{\mathrm{br}}^+ - \mathcal{B}_{\mathrm{br}}^- = \int_0^1 B_{\mathrm{br}}(t)\,\mathrm{d}t \sim N(0,1/12). \tag{238}$$

By symmetry $\mathcal{B}_{\mathrm{br}}^- \overset{\mathrm{d}}{=} \mathcal{B}_{\mathrm{br}}^+$, so we concentrate on $\mathcal{B}_{\mathrm{br}}^+$. Let $\psi_{\mathrm{br}}^+(t) := \mathbb{E}\,e^{-t\mathcal{B}_{\mathrm{br}}^+}$ be its Laplace transform. Perman and Wellner [38] gave (using the notation $\Psi_0^+$ for $\psi_{\mathrm{br}}^+$) the formula

$$\int_0^\infty e^{-xt}\psi_{\mathrm{br}}^+\big(\sqrt{2}\,t^{3/2}\big)\frac{\mathrm{d}t}{\sqrt{\pi t}} = 2\frac{\mathrm{Ai}(x)}{x^{1/2}\mathrm{Ai}(x)-\mathrm{Ai}'(x)}, \qquad x>0, \tag{239}$$

see also Appendix C.2.

Define $D_n^+$ by

$$\mathbb{E}(\mathcal{B}_{\mathrm{br}}^+)^n = \frac{\sqrt{\pi}\, 2^{-n/2} n!}{\Gamma((3n+1)/2)} D_n^+, \qquad n \geq 0. \tag{240}$$

(Perman and Wellner [38] use $L_n^+$ for our $D_n^+$ and $\mu_n^+$ for $\mathbb{E}(\mathcal{B}_{\mathrm{br}}^+)^n$.) An asymptotic expansion of the left hand side of (239) yields, arguing as for e.g. (98) and (136),

$$\frac{2\mathrm{Ai}(x)}{x^{1/2}\mathrm{Ai}(x) - \mathrm{Ai}'(x)} \sim \sum_{k=0}^{\infty} (-1)^k D_k^+ x^{-3k/2-1/2}, \qquad x \to \infty, \tag{241}$$

and thus by (102), (103), (183), (106), (107),

$$\sum_{k=0}^{\infty} D_k^+ z^k = \frac{2\sum_{k=0}^{\infty} \alpha_k z^k}{\sum_{k=0}^{\infty} \alpha_k z^k + \sum_{k=0}^{\infty} \alpha'_k z^k} = \frac{\sum_{k=0}^{\infty} \alpha_k z^k}{\sum_{k=0}^{\infty} \frac{1}{1-6k} \alpha_k z^k}, \tag{242}$$

which leads to the recursion by Perman and Wellner [38]

$$D_n^+ = \alpha_n + \sum_{k=1}^{n} \frac{1}{6k-1} \alpha_k D_{n-k}^+, \qquad n \geq 0. \tag{243}$$

Using (108) and (137), we obtain from (242) also

$$\sum_{k=0}^{\infty} D_k^+ z^k = \frac{2}{1 - \sum_{k=0}^{\infty} 2K_k z^k} = \frac{2\sum_{k=0}^{\infty} D_k z^k}{1 + \sum_{k=0}^{\infty} D_k z^k}, \tag{244}$$

and the corresponding recursions

$$D_n^+ = \sum_{k=1}^{n} K_k D_{n-k}^+, \qquad n \geq 1, \tag{245}$$

$$D_n^+ = D_n - \tfrac{1}{2} \sum_{k=1}^{n} D_k D_{n-k}^+ \qquad n \geq 0. \tag{246}$$

Some numerical values are given in Table 6; see further Perman and Wellner [38] (but beware of typos for $k = 5$ and 7). We define $D_n^{+*} := 2^{3n} D_n^+$; these are integers by (245).

Tolmatz [53] gives the asymptotics, see also Janson and Louchard [25] and compare (156) and (155),

$$\mathbb{E}(\mathcal{B}_{\mathrm{br}}^+)^n \sim \tfrac{1}{2}\, \mathbb{E}\, \mathcal{B}_{\mathrm{br}}^n \sim \frac{1}{\sqrt{2}} \Big(\frac{n}{12e}\Big)^{n/2} \qquad \text{as } n \to \infty, \tag{247}$$

or, equivalently,

$$D_n^+ \sim \tfrac{1}{2} D_n \sim \frac{1}{2\pi} \Big(\frac{3}{4}\Big)^n (n-1)! \sim \frac{1}{\sqrt{2\pi n}} \Big(\frac{3n}{4e}\Big)^n \qquad \text{as } n \to \infty. \tag{248}$$

| | | | | |
|---|---|---|---|---|
| $\mathbb{E}(\mathcal{B}_{\mathrm{br}}^+)^0 = 1$ | $\mathbb{E}\,\mathcal{B}_{\mathrm{br}}^+ = \frac{1}{8}\sqrt{\frac{\pi}{2}}$ | $\mathbb{E}(\mathcal{B}_{\mathrm{br}}^+)^2 = \frac{1}{20}$ | $\mathbb{E}(\mathcal{B}_{\mathrm{br}}^+)^3 = \frac{71}{4096}\sqrt{\frac{\pi}{2}}$ | $\mathbb{E}(\mathcal{B}_{\mathrm{br}}^+)^4 = \frac{211}{18480}$ |
| $D_0^+ = 1$ | $D_1^+ = \frac{1}{8}$ | $D_2^+ = \frac{3}{32}$ | $D_3^+ = \frac{71}{512}$ | $D_4^+ = \frac{633}{2048}$ |
| $D_0^{+*} = 1$ | $D_1^{+*} = 1$ | $D_2^{+*} = 6$ | $D_3^{+*} = 71$ | $D_4^{+*} = 1266$ |

Table 6

*Some numerical values for the positive part of a Brownian bridge.*

Joint moments of $\mathcal{B}_{\mathrm{br}}^+$ and $\mathcal{B}_{\mathrm{br}}^-$ can be computed by the same method from the double joint Laplace transform computed by Perman and Wellner [38], see (369) in Appendix C.2 below. Taking $\lambda = 1$ in (369), the left hand side can be expanded in an asymptotic double series, for $\xi, \eta \searrow 0$,

$$\sim \sum_{k,l=0}^{\infty} (-1)^{k+l} \frac{2^{k/2+l/2}}{\sqrt{2\pi}\,k!\,l!} \Gamma\Big(\frac{3(k+l)+1}{2}\Big)\,\mathbb{E}\big((\mathcal{B}_{\mathrm{br}}^+)^k (\mathcal{B}_{\mathrm{br}}^-)^l\big)\xi^k\eta^l. \tag{249}$$

Using (101) and (108), we find from (369) the identity (for formal power series), after replacing $\xi$ by $-\xi$ and $\eta$ by $-\eta$,

$$\begin{aligned}\sum_{k,l=0}^{\infty} \frac{2^{k/2+l/2}}{\sqrt{\pi}\,k!\,l!} \Gamma\Big(\frac{3(k+l)+1}{2}\Big)\,\mathbb{E}\big((\mathcal{B}_{\mathrm{br}}^+)^k (\mathcal{B}_{\mathrm{br}}^-)^l\big)\xi^k\eta^l \\ = 2\left(\frac{\sum_{k=0}^\infty \alpha'_k\xi^k}{\sum_{k=0}^\infty \alpha_k\xi^k} + \frac{\sum_{k=0}^\infty \alpha'_k\eta^k}{\sum_{k=0}^\infty \alpha_k\eta^k}\right)^{-1} \\ = \frac{1}{-\sum_{k=0}^\infty K_k\xi^k - \sum_{k=0}^\infty K_k\eta^k}.\end{aligned} \tag{250}$$

Note that $\eta = \xi$ yields (137) (using (135) and (152)) and that $\eta = 0$ yields (244).

Define $D_{k,l}^{\pm}$ as the coefficient of $\xi^k\eta^l$ in the left hand side of (250). Thus

$$\mathbb{E}\big((\mathcal{B}_{\mathrm{br}}^+)^k (\mathcal{B}_{\mathrm{br}}^-)^l\big) = \frac{2^{-(k+l)/2}\sqrt{\pi}\,k!\,l!}{\Gamma\big((3(k+l)+1)/2\big)} D_{k,l}^{\pm}, \qquad k, l \ge 0 \tag{251}$$

and

$$\sum_{k,l=0}^{\infty} D_{k,l}^{\pm}\xi^k\eta^l = \frac{1}{-\sum_{k=0}^\infty K_k(\xi^k+\eta^k)} = \frac{1}{1-\sum_{k=1}^\infty K_k(\xi^k+\eta^k)}, \tag{252}$$

which yields the recursion formula, with $D_{0,0}^{\pm} = 1$, for $k + l > 0$,

$$D_{k,l}^{\pm} = \sum_{j=1}^{k} K_j D_{k-j,l}^{\pm} + \sum_{j=1}^{l} K_j D_{k,l-j}^{\pm}. \tag{253}$$

We have $D^{\pm}_{k,0} = D^{\pm}_{0,k} = D^{+}_{k}$ and

$$\sum_{k=0}^{n} D^{\pm}_{k,n-k} = D_n, \qquad n \geq 0. \tag{254}$$

We remark further that if $n = 2m$ is even, then (251) and (238) yield

$$\begin{aligned}\sum_{j=0}^{n} (-1)^{n-j} D^{\pm}_{j,n-j} &= \frac{2^{n/2}\Gamma\big((3n+1)/2\big)}{\sqrt{\pi}\, n!} \sum_{j=0}^{n} (-1)^{n-j} \binom{n}{j} \mathbb{E}\big(\mathcal{B}^{+}_{\mathrm{br}}\big)^{j} \big(\mathcal{B}^{-}_{\mathrm{br}}\big)^{n-j} \\ &= \frac{2^{n/2}\Gamma\big((3n+1)/2\big)}{\sqrt{\pi}\, n!} \mathbb{E}\big(\mathcal{B}^{+}_{\mathrm{br}} - \mathcal{B}^{-}_{\mathrm{br}}\big)^{n} \\ &= \frac{2^{n/2}\Gamma\big((3n+1)/2\big)}{\sqrt{\pi}\, n!} 12^{-n/2} \frac{n!}{2^{n/2}(n/2)!} \\ &= 2^{-8m} 3^{-m} \frac{(6m)!}{(3m)!\, m!}, \end{aligned} \tag{255}$$

while the sum vanishes by symmetry if $n$ is odd. Consequently, $\eta = -\xi$ in (250) or (252) yields

$$\sum_{m=0}^{\infty} 2^{-8m} 3^{-m} \frac{(6m)!}{(3m)!\, m!} \xi^{2m} = \frac{1}{-2\sum_{m=0}^{\infty} K_{2m}\xi^{2m}} = \frac{1}{1 - 2\sum_{m=1}^{\infty} K_{2m}\xi^{2m}}, \tag{256}$$

which can be written as a recursion relation for $K_{2m}$ with even indices only.

Some numerical values are given in Table 7. In particular, as found by Perman and Wellner [38],

$$\operatorname{Cov}(\mathcal{B}^{+}_{\mathrm{br}}, \mathcal{B}^{-}_{\mathrm{br}}) = \frac{1}{120} - \frac{\pi}{128}, \tag{257}$$

$$\operatorname{Corr}(\mathcal{B}^{+}_{\mathrm{br}}, \mathcal{B}^{-}_{\mathrm{br}}) = \frac{128 - 120\pi}{768 - 120\pi} \approx -0.636791. \tag{258}$$

We do not know of any formula for the density function $f^{+}_{\mathrm{br}}$ or distribution function $F^{+}_{\mathrm{br}}$, except the Laplace inversion formulas in Tolmatz [53] and Janson and Louchard [25] (which prove that $f^{+}_{\mathrm{br}}$ exists and is continuous). For $x \to \infty$ we have by Tolmatz [53] and Janson and Louchard [25] asymptotic expansions beginning with

$$f^{+}_{\mathrm{br}}(x) \sim \frac{\sqrt{6}}{\sqrt{\pi}} e^{-6x^2} \left(1 + \frac{1}{36} x^{-2} - \frac{7}{5184} x^{-4} + \dots\right), \tag{259}$$

$$1 - F^{+}_{\mathrm{br}}(x) \sim \frac{1}{2\sqrt{6\pi}} x^{-1} e^{-6x^2} \left(1 - \frac{1}{18} x^{-2} + \frac{65}{5184} x^{-4} + \dots\right). \tag{260}$$

We do not know any asymptotic results as $x \to 0$, but the results on existence of negative moments in Section 29 suggest that $f^{+}_{\mathrm{br}}(x) \asymp x^{-1/3}$. More precisely we conjecture, based on (317), that

$$f^{+}_{\mathrm{br}}(x) \overset{?}{\sim} \frac{2^{1/3} 3^{5/6} \Gamma(1/3)^5}{(2\pi)^3} x^{-1/3}, \qquad x \to 0. \tag{261}$$

| | | | |
|---|---|---|---|
| $\mathbb{E}\,\mathcal{B}_{\mathrm{br}}^{0,0} = 1$ | $\mathbb{E}\,\mathcal{B}_{\mathrm{br}}^{0,1} = \frac{\sqrt{2\pi}}{16}$ | $\mathbb{E}\,\mathcal{B}_{\mathrm{br}}^{0,2} = \frac{1}{20}$ | $\mathbb{E}\,\mathcal{B}_{\mathrm{br}}^{0,3} = \frac{71\sqrt{2\pi}}{8192}$ |
| $\mathbb{E}\,\mathcal{B}_{\mathrm{br}}^{1,0} = \frac{\sqrt{2\pi}}{16}$ | $\mathbb{E}\,\mathcal{B}_{\mathrm{br}}^{1,1} = \frac{1}{120}$ | $\mathbb{E}\,\mathcal{B}_{\mathrm{br}}^{1,2} = \frac{13\sqrt{2\pi}}{24576}$ | $\mathbb{E}\,\mathcal{B}_{\mathrm{br}}^{1,3} = \frac{1}{2880}$ |
| $\mathbb{E}\,\mathcal{B}_{\mathrm{br}}^{2,0} = \frac{1}{20}$ | $\mathbb{E}\,\mathcal{B}_{\mathrm{br}}^{2,1} = \frac{13\sqrt{2\pi}}{24576}$ | $\mathbb{E}\,\mathcal{B}_{\mathrm{br}}^{2,2} = \frac{43}{332640}$ | $\mathbb{E}\,\mathcal{B}_{\mathrm{br}}^{2,3} = \frac{17\sqrt{2\pi}}{1835008}$ |
| $\mathbb{E}\,\mathcal{B}_{\mathrm{br}}^{3,0} = \frac{71\sqrt{2\pi}}{8192}$ | $\mathbb{E}\,\mathcal{B}_{\mathrm{br}}^{3,1} = \frac{1}{2880}$ | $\mathbb{E}\,\mathcal{B}_{\mathrm{br}}^{3,2} = \frac{17\sqrt{2\pi}}{1835008}$ | $\mathbb{E}\,\mathcal{B}_{\mathrm{br}}^{3,3} = \frac{11}{3564288}$ |
| $D_{0,0}^{\pm} = 1$ | $D_{0,1}^{\pm} = \frac{1}{8}$ | $D_{0,2}^{\pm} = \frac{3}{32}$ | $D_{0,3}^{\pm} = \frac{71}{512}$ |
| $D_{1,0}^{\pm} = \frac{1}{8}$ | $D_{1,1}^{\pm} = \frac{1}{32}$ | $D_{1,2}^{\pm} = \frac{13}{512}$ | $D_{1,3}^{\pm} = \frac{77}{2048}$ |
| $D_{2,0}^{\pm} = \frac{3}{32}$ | $D_{2,1}^{\pm} = \frac{13}{512}$ | $D_{2,2}^{\pm} = \frac{43}{2048}$ | $D_{2,3}^{\pm} = \frac{255}{8192}$ |
| $D_{3,0}^{\pm} = \frac{71}{512}$ | $D_{3,1}^{\pm} = \frac{77}{2048}$ | $D_{3,2}^{\pm} = \frac{255}{8192}$ | $D_{3,3}^{\pm} = \frac{3025}{65536}$ |

TABLE 7
*Some numerical values for the positive and negative parts of a Brownian bridge.* $\mathbb{E}\,\mathcal{B}_{\mathrm{br}}^{k,l}$ *is an abbreviation for* $\mathbb{E}(\mathcal{B}_{\mathrm{br}}^{+})^k(\mathcal{B}_{\mathrm{br}}^{-})^l$.

## 25. Positive part of a Brownian motion

Define
$$\mathcal{B}_{\mathrm{bm}}^{\pm} := \int_0^1 B(t)_{\pm}\,\mathrm{d}t; \tag{262}$$
thus $\mathcal{B}_{\mathrm{bm}}^{+}$ is the average of the positive part and $\mathcal{B}_{\mathrm{bm}}^{-}$ the average of the negative part of a Brownian motion on $[0,1]$. (Perman and Wellner [38] use $A^+ = A^+(1)$ for $\mathcal{B}_{\mathrm{bm}}^{+}$.) In particular,
$$\mathcal{B}_{\mathrm{bm}} = \mathcal{B}_{\mathrm{bm}}^{+} + \mathcal{B}_{\mathrm{bm}}^{-}. \tag{263}$$
Note also that the difference is Gaussian:
$$\mathcal{B}_{\mathrm{bm}}^{+} - \mathcal{B}_{\mathrm{bm}}^{-} = \int_0^1 B(t)\,\mathrm{d}t \sim N(0,1/3). \tag{264}$$

Again, there is symmetry: $\mathcal{B}_{\mathrm{bm}}^{-} \overset{\mathrm{d}}{=} \mathcal{B}_{\mathrm{bm}}^{+}$, and we concentrate on $\mathcal{B}_{\mathrm{bm}}^{+}$. Let $\psi_{\mathrm{bm}}^{+}(t) := \mathbb{E}\,e^{-t\mathcal{B}_{\mathrm{bm}}^{+}}$ be its Laplace transform. Perman and Wellner [38] gave (using the notation $\Psi^+$ for $\psi_{\mathrm{bm}}^{+}$) the formula
$$\int_0^\infty e^{-xt}\psi_{\mathrm{bm}}^{+}\bigl(\sqrt{2}\,t^{3/2}\bigr)\,\mathrm{d}t = \frac{x^{-1/2}\mathrm{Ai}(x) + \mathrm{AI}(x)}{x^{1/2}\mathrm{Ai}(x) - \mathrm{Ai}'(x)}, \qquad x > 0, \tag{265}$$
see also Appendix C.1.

Define $L_n^+$ by
$$\mathbb{E}(\mathcal{B}_{\mathrm{bm}}^{+})^n = \frac{2^{-n/2}n!}{\Gamma((3n+2)/2)}L_n^+, \qquad n \geq 0. \tag{266}$$

| | | | | |
|---|---|---|---|---|
| $\mathbb{E}(\mathcal{B}_{\mathrm{bm}}^+)^0 = 1$ | $\mathbb{E}\,\mathcal{B}_{\mathrm{bm}}^+ = \dfrac{\sqrt{2}}{3\sqrt{\pi}}$ | $\mathbb{E}(\mathcal{B}_{\mathrm{bm}}^+)^2 = \dfrac{17}{96}$ | $\mathbb{E}(\mathcal{B}_{\mathrm{bm}}^+)^3 = \dfrac{251\sqrt{2}}{1260\sqrt{\pi}}$ | $\mathbb{E}(\mathcal{B}_{\mathrm{bm}}^+)^4 = \dfrac{6989}{40960}$ |
| $L_0^+ = 1$ | $L_1^+ = \dfrac{1}{2}$ | $L_2^+ = \dfrac{17}{16}$ | $L_3^+ = \dfrac{251}{64}$ | $L_4^+ = \dfrac{20967}{1024}$ |
| $L_0^{+*} = 1$ | $L_1^{+*} = 4$ | $L_2^{+*} = 68$ | $L_3^{+*} = 2008$ | $L_4^{+*} = 83868$ |

Table 8

*Some numerical values for the positive part of a Brownian motion.*

(Perman and Wellner [38] use $K_n^+$ for our $L_n^+$ and $\nu_n^+$ for $\mathbb{E}(\mathcal{B}_{\mathrm{bm}}^+)^n$.) An asymptotic expansion of the left hand side of (265) yields, arguing as for e.g. (98) and (182),

$$x^{-1}\frac{\mathrm{Ai}(x)+x^{1/2}\mathrm{AI}(x)}{\mathrm{Ai}(x)-x^{-1/2}\mathrm{Ai}'(x)} \sim \sum_{k=0}^{\infty}(-1)^k L_k^+ x^{-3k/2-1}, \qquad x\to\infty. \tag{267}$$

The asymptotic expansions (102), (103), (183) of $\mathrm{Ai}(x)$, $\mathrm{Ai}'(x)$ and $\mathrm{AI}(x)$ yield, using (106) and (107), the identity

$$\sum_{k=0}^{\infty} L_k^+ z^k = \frac{\sum_{k=0}^{\infty}\alpha_k z^k + \sum_{k=0}^{\infty}\beta_k z^k}{\sum_{k=0}^{\infty}\alpha_k z^k + \sum_{k=0}^{\infty}\alpha_k' z^k} = \frac{\sum_{k=0}^{\infty}(\alpha_k+\beta_k)z^k}{\sum_{k=0}^{\infty}\frac{2}{1-6k}\alpha_k z^k}, \tag{268}$$

which leads to the recursion by Perman and Wellner [38]

$$L_n^+ = \tfrac{1}{2}\alpha_n + \tfrac{1}{2}\beta_n + \sum_{k=1}^{n}\frac{1}{6k-1}\alpha_k L_{n-k}^+, \qquad n\ge 0. \tag{269}$$

Some numerical values are given in Table 8; see further Perman and Wellner [38]. We define $L_n^{+*} := 2^{3n}L_n^+$; these are integers as a consequence of (298) below.

Janson and Louchard [25] give the asymptotics, compare (191) and (192),

$$\mathbb{E}(\mathcal{B}_{\mathrm{bm}}^+)^n \sim \tfrac{1}{2}\,\mathbb{E}\,\mathcal{B}_{\mathrm{bm}}^n \sim \frac{1}{\sqrt{2}}\Big(\frac{n}{3e}\Big)^{n/2}, \tag{270}$$

$$L_n^+ \sim \tfrac{1}{2}L_n \sim \frac{\sqrt{3}}{2}\Big(\frac{3n}{2e}\Big)^n, \qquad \text{as } n\to\infty. \tag{271}$$

Joint moments of $\mathcal{B}_{\mathrm{bm}}^+$ and $\mathcal{B}_{\mathrm{bm}}^-$ can be computed by the same method from the double joint Laplace transform computed by Perman and Wellner [38], see (361) in Appendix C.1 below. Taking $\lambda = 1$ in (361), the left hand side can be expanded in an asymptotic double series, for $\xi,\eta \searrow 0$,

$$\sim \sum_{k,l=0}^{\infty}(-1)^{k+l}\frac{2^{k/2+l/2}}{k!\,l!}\Gamma\Big(\frac{3(k+l)}{2}+1\Big)\,\mathbb{E}\big((\mathcal{B}_{\mathrm{bm}}^+)^k(\mathcal{B}_{\mathrm{bm}}^-)^l\big)\xi^k\eta^l. \tag{272}$$

Using (102), (103), (183), (106), (107), (108), (196), we find from (361) the identity (for formal power series), after replacing $\xi$ by $-\xi$ and $\eta$ by $-\eta$,

$$\begin{aligned}\sum_{k,l=0}^{\infty} \frac{2^{k/2+l/2}}{k!\,l!}\Gamma\Big(\frac{3(k+l)}{2}+1\Big)\,\mathbb{E}\big((\mathcal{B}^+_{\mathrm{bm}})^k(\mathcal{B}^-_{\mathrm{bm}})^l\big)\xi^k\eta^l \\ = \left(\frac{\sum_{k=0}^\infty \beta_k\xi^k}{\sum_{k=0}^\infty \alpha_k\xi^k}+\frac{\sum_{k=0}^\infty \beta_k\eta^k}{\sum_{k=0}^\infty \alpha_k\eta^k}\right)\left(\frac{\sum_{k=0}^\infty \alpha'_k\xi^k}{\sum_{k=0}^\infty \alpha_k\xi^k}+\frac{\sum_{k=0}^\infty \alpha'_k\eta^k}{\sum_{k=0}^\infty \alpha_k\eta^k}\right)^{-1} \\ = \frac{\sum_{k=0}^\infty Q_k\xi^k+\sum_{k=0}^\infty Q_k\eta^k}{-2\sum_{k=0}^\infty K_k\xi^k-2\sum_{k=0}^\infty K_k\eta^k}.\end{aligned} \tag{273}$$

Note that $\eta=\xi$ yields (185) and that $\eta=0$ yields (268).

Define $L^{\pm}_{k,l}$ as the coefficient of $\xi^k\eta^l$ in the left hand side of (273). Thus

$$\mathbb{E}\big((\mathcal{B}^+_{\mathrm{bm}})^k(\mathcal{B}^-_{\mathrm{bm}})^l\big) = \frac{2^{-(k+l)/2}k!\,l!}{\Gamma\big(3(k+l)/2+1\big)}L^{\pm}_{k,l}, \qquad k,l\ge 0 \tag{274}$$

and

$$\sum_{k,l=0}^{\infty} L^{\pm}_{k,l}\xi^k\eta^l = \frac{\sum_{k=0}^\infty Q_k(\xi^k+\eta^k)}{-2\sum_{k=0}^\infty K_k(\xi^k+\eta^k)}, \tag{275}$$

which yields the recursion formula, with $L^{\pm}_{0,0}=1$,

$$L^{\pm}_{k,l} = \sum_{j=1}^{k} K_j L^{\pm}_{k-j,l} + \sum_{j=1}^{l} K_j L^{\pm}_{k,l-j} + \tfrac12\delta_{l0}Q_k + \tfrac12\delta_{k0}Q_l, \tag{276}$$

where $\delta$ is Kronecker's delta.

We have $L^{\pm}_{k,0}=L^{\pm}_{0,k}=L^+_k$ and

$$\sum_{k=0}^{n} L^{\pm}_{k,n-k} = L_n, \qquad n\ge 0. \tag{277}$$

We remark further that if $n=2m$ is even, then (274) and (264) yield

$$\begin{aligned}\sum_{j=0}^{n}(-1)^{n-j}L^{\pm}_{j,n-j} &= \frac{2^{n/2}\Gamma(3n/2+1)}{n!}\sum_{j=0}^{n}(-1)^{n-j}\binom{n}{j}\,\mathbb{E}\big(\mathcal{B}^+_{\mathrm{bm}}\big)^j\big(\mathcal{B}^-_{\mathrm{bm}}\big)^{n-j} \\ &= \frac{2^{n/2}\Gamma(3n/2+1)}{n!}\,\mathbb{E}\big(\mathcal{B}^+_{\mathrm{bm}}-\mathcal{B}^-_{\mathrm{bm}}\big)^n \\ &= \frac{2^{n/2}\Gamma(3n/2+1)}{n!}3^{-n/2}\frac{n!}{2^{n/2}(n/2)!} = \frac{(3m)!}{3^m m!},\end{aligned} \tag{278}$$

while the sum vanishes by symmetry if $n$ is odd. Consequently, $\eta=-\xi$ in (273) or (275) yields

$$\sum_{m=0}^{\infty}\frac{(3m)!}{3^m m!}\xi^{2m} = \frac{\sum_{m=0}^\infty Q_{2m}\xi^{2m}}{-2\sum_{m=0}^\infty K_{2m}\xi^{2m}}, \tag{279}$$

$$\mathbb{E}\,\mathcal{B}_{\mathrm{bm}}^{0,0}=1 \qquad \mathbb{E}\,\mathcal{B}_{\mathrm{bm}}^{0,1}=\frac{\sqrt{2}}{3\sqrt{\pi}} \qquad \mathbb{E}\,\mathcal{B}_{\mathrm{bm}}^{0,2}=\frac{17}{96} \qquad \mathbb{E}\,\mathcal{B}_{\mathrm{bm}}^{0,3}=\frac{251\sqrt{2}}{1260\sqrt{\pi}}$$

$$\mathbb{E}\,\mathcal{B}_{\mathrm{bm}}^{1,0}=\frac{\sqrt{2}}{3\sqrt{\pi}} \qquad \mathbb{E}\,\mathcal{B}_{\mathrm{bm}}^{1,1}=\frac{1}{96} \qquad \mathbb{E}\,\mathcal{B}_{\mathrm{bm}}^{1,2}=\frac{\sqrt{2}}{315\sqrt{\pi}} \qquad \mathbb{E}\,\mathcal{B}_{\mathrm{bm}}^{1,3}=\frac{149}{122880}$$

$$\mathbb{E}\,\mathcal{B}_{\mathrm{bm}}^{2,0}=\frac{17}{96} \qquad \mathbb{E}\,\mathcal{B}_{\mathrm{bm}}^{2,1}=\frac{\sqrt{2}}{315\sqrt{\pi}} \qquad \mathbb{E}\,\mathcal{B}_{\mathrm{bm}}^{2,2}=\frac{109}{368640} \qquad \mathbb{E}\,\mathcal{B}_{\mathrm{bm}}^{2,3}=\frac{31\sqrt{2}}{300300\sqrt{\pi}}$$

$$\mathbb{E}\,\mathcal{B}_{\mathrm{bm}}^{3,0}=\frac{251\sqrt{2}}{1260\sqrt{\pi}} \qquad \mathbb{E}\,\mathcal{B}_{\mathrm{bm}}^{3,1}=\frac{149}{122880} \qquad \mathbb{E}\,\mathcal{B}_{\mathrm{bm}}^{3,2}=\frac{31\sqrt{2}}{300300\sqrt{\pi}} \qquad \mathbb{E}\,\mathcal{B}_{\mathrm{bm}}^{3,3}=\frac{9391}{660602880}$$

$$L^{\pm}_{0,0}=1 \qquad L^{\pm}_{0,1}=\frac{1}{2} \qquad L^{\pm}_{0,2}=\frac{17}{16} \qquad L^{\pm}_{0,3}=\frac{251}{64}$$

$$L^{\pm}_{1,0}=\frac{1}{2} \qquad L^{\pm}_{1,1}=\frac{1}{8} \qquad L^{\pm}_{1,2}=\frac{3}{16} \qquad L^{\pm}_{1,3}=\frac{149}{256}$$

$$L^{\pm}_{2,0}=\frac{17}{16} \qquad L^{\pm}_{2,1}=\frac{3}{16} \qquad L^{\pm}_{2,2}=\frac{109}{512} \qquad L^{\pm}_{2,3}=\frac{279}{512}$$

$$L^{\pm}_{3,0}=\frac{251}{64} \qquad L^{\pm}_{3,1}=\frac{149}{256} \qquad L^{\pm}_{3,2}=\frac{279}{512} \qquad L^{\pm}_{3,3}=\frac{9391}{8192}$$

TABLE 9
*Some numerical values for the positive and negative parts of a Brownian motion.* $\mathbb{E}\,\mathcal{B}_{\mathrm{bm}}^{k,l}$ *is an abbreviation for* $\mathbb{E}(\mathcal{B}_{\mathrm{bm}}^{+})^k(\mathcal{B}_{\mathrm{bm}}^{-})^l$.

a relation between $Q_j$ and $K_k$ with even indices only.

Some numerical values are given in Table 9. In particular, as found by Perman and Wellner [38],

$$\operatorname{Cov}(\mathcal{B}_{\mathrm{bm}}^{+},\mathcal{B}_{\mathrm{bm}}^{-})=\frac{1}{96}-\frac{2}{9\pi}, \tag{280}$$

$$\operatorname{Corr}(\mathcal{B}_{\mathrm{bm}}^{+},\mathcal{B}_{\mathrm{bm}}^{-})=\frac{3\pi-64}{51\pi-64}\approx-0.567185. \tag{281}$$

We do not know of any formula for the density function $f_{\mathrm{bm}}^{+}$ or distribution function $F_{\mathrm{bm}}^{+}$, except the Laplace inversion formulas in Janson and Louchard [25] (which prove that $f_{\mathrm{bm}}^{+}$ exists and is continuous). For $x\to\infty$ we have by Janson and Louchard [25] asymptotic expansions beginning with

$$f_{\mathrm{bm}}^{+}(x)\sim\frac{\sqrt{3}}{\sqrt{2\pi}}\,e^{-3x^2/2}\left(1+\frac{1}{36}\,x^{-2}-\frac{5}{648}\,x^{-4}+\dots\right), \tag{282}$$

$$\mathbb{P}(\mathcal{B}_{\mathrm{bm}}^{+}>x)\sim\frac{1}{\sqrt{6\pi}}\,x^{-1}e^{-3x^2/2}\left(1-\frac{11}{36}\,x^{-2}+\frac{193}{648}\,x^{-4}+\dots\right). \tag{283}$$

We do not know any asymptotic results as $x\to0$, but the results on existence of negative moments in Section 29 suggest that $f_{\mathrm{bm}}^{+}(x)\asymp x^{-2/3}$. More precisely we conjecture, based on (319), that

$$f_{\mathrm{bm}}^{+}(x)\overset{?}{\sim}2^{-5/6}3^{1/6}\pi^{-3/2}\Gamma(1/3)x^{-2/3},\qquad x\to0. \tag{284}$$

## 26. Convolutions and generating functions

Many of the formulas above can be written as convolutions of sequences, as is done by Perman and Wellner [38]. This gives the simple formulas below, where a letter $X$ stands for the sequence $(X_n)_0^\infty$ and $X * Y$ is the sequence defined by $(X*Y)_n := \sum_{k=0}^n X_k Y_{n-k}$. Alternatively, the formulas can be interpreted as identities for formal power series, if we instead interpret $X$ as the generating function $\sum_{n=0}^\infty X_n z^n$ and $*$ as ordinary multiplication. We let $\mathbf{1}$ denote the sequence $(\delta_{n0}) = (1, 0, 0, \dots)$ with the generating function 1.

Note that the sequences denoted by various letters below, except $K$, all have $X_0 = 1$, and thus the generating functions have constant term 1; the exception is $K_0 = -\frac{1}{2}$, which explains why $-2K$ occurs frequently. (It would be more consistent to give a new name to $-2K_n$ and use it instead of $K_n$ in our formulas (*mutatis mutandis*), but we have kept $K_n$ for historical reasons, and because all numbers $K_n$ with $n \geq 1$ are positive. See also (311), which shows that $(6n-2)K_n$ is more natural for some other purposes.)

In the following list of formulas, we also give references to the corresponding equations above. Many of these formulas can also be found in Perman and Wellner [38].

$$-2\alpha * K = \alpha' \qquad (108), (70) \tag{285}$$

$$\alpha' * D = \alpha \qquad (137), (148) \tag{286}$$

$$-2K * D = \mathbf{1} \qquad (152), (153), (154) \tag{287}$$

$$\alpha' * L = \beta \qquad (185), (186) \tag{288}$$

$$\alpha * Q = \beta \qquad (196), (201) \tag{289}$$

$$-2K * L = Q \qquad (202), (203) \tag{290}$$

$$D * Q = L \qquad (204), (205) \tag{291}$$

$$Q * Q = W \qquad (231), (232) \tag{292}$$

$$\alpha * \alpha * W = \beta * \beta \qquad (231) \tag{293}$$

$$(\alpha + \alpha') * D^+ = 2\alpha \qquad (242), (243) \tag{294}$$

$$D^+ - 2K * D^+ = 2 \cdot \mathbf{1} \qquad (244), (245) \tag{295}$$

$$D^+ + D * D^+ = 2D \qquad (244), (246) \tag{296}$$

$$(\alpha + \alpha') * L^+ = \alpha + \beta \qquad (268), (269) \tag{297}$$

$$D^+ * Q + D^+ = 2L^+ \tag{298}$$

The final relation (298) follows from (294), (289) and (297) by simple algebra with generating functions.

## 27. Randomly stopped Brownian motion

Several of the formulas above have simple interpretations in terms of a Brownian motion stopped at a random time, which has been found and used by Perman and Wellner [38].

Let $\mathcal{B}_{\mathrm{bm}}(t) := \int_0^t |B(s)|\,\mathrm{d}s$; by Brownian scaling, $\mathcal{B}_{\mathrm{bm}}(t) \overset{\mathrm{d}}{=} t^{3/2}\mathcal{B}_{\mathrm{bm}}$. Let further $Z \sim \mathrm{Exp}(1)$ be an exponentially distributed random variable independent of the Brownian motion $B$. The area of the (reflected) Brownian motion stopped at $Z$ is thus

$$\mathcal{B}_{\mathrm{bm}}(Z) \overset{\mathrm{d}}{=} Z^{3/2}\mathcal{B}_{\mathrm{bm}} \tag{299}$$

with moments given by, see (181),

$$\mathbb{E}\big(\mathcal{B}_{\mathrm{bm}}(Z)\big)^n = \mathbb{E}\, Z^{3n/2}\,\mathbb{E}\,\mathcal{B}_{\mathrm{bm}}^n = 2^{-n/2} n! L_n. \tag{300}$$

Thus,

$$L_n = \frac{\mathbb{E}\big(\sqrt{2}\,\mathcal{B}_{\mathrm{bm}}(Z)\big)^n}{n!}, \qquad n \ge 0, \tag{301}$$

and the generating function $\sum_{k=0}^{\infty} L_k z^k$ in (185) is the exponential generating function of the moments of $\sqrt{2}\,\mathcal{B}_{\mathrm{bm}}(Z)$, i.e. the moment generating function of $\sqrt{2}\,\mathcal{B}_{\mathrm{bm}}$ interpreted as a formal power series. (This series diverges, so the moment generating function in the ordinary sense $\mathbb{E}\exp\big(z\sqrt{2}\,\mathcal{B}_{\mathrm{bm}}(Z)\big)$ is infinite for all $z > 0$.) Furthermore, the Laplace transform of $\sqrt{2}\,\mathcal{B}_{\mathrm{bm}}(Z)$ is given by, see (362) below, for $\xi > 0$,

$$\mathbb{E}\, e^{-\xi\sqrt{2}\,\mathcal{B}_{\mathrm{bm}}(Z)} = \int_0^\infty e^{-t}\psi_{\mathrm{bm}}\big(\sqrt{2}\,\xi t^{3/2}\big)\,\mathrm{d}t = -\frac{\mathrm{AI}(\xi^{-2/3})}{\xi^{2/3}\mathrm{Ai}'(\xi^{-2/3})}, \tag{302}$$

which is equivalent to (173) by a change of variables.

Similarly, with obvious definitions, (266) can be written

$$L_n^+ = \frac{\mathbb{E}\big(\sqrt{2}\,\mathcal{B}_{\mathrm{bm}}^+(Z)\big)^n}{n!}, \qquad n \ge 0, \tag{303}$$

and the generating function $\sum_{k=0}^{\infty} L_k^+ z^k$ in (268) is the moment generating function of $\sqrt{2}\,\mathcal{B}_{\mathrm{bm}}^+$ interpreted as a formal power series. The Laplace transform of $\sqrt{2}\,\mathcal{B}_{\mathrm{bm}}^+(Z)$ is given by, see (363) below, for $\xi > 0$,

$$\mathbb{E}\, e^{-\xi\sqrt{2}\,\mathcal{B}_{\mathrm{bm}}^+(Z)} = \int_0^\infty e^{-t}\psi_{\mathrm{bm}}^+\big(\sqrt{2}\,\xi t^{3/2}\big)\,\mathrm{d}t = \frac{\mathrm{Ai}(\xi^{-2/3}) + \xi^{-1/3}\mathrm{AI}(\xi^{-2/3})}{\mathrm{Ai}(\xi^{-2/3}) - \xi^{1/3}\mathrm{Ai}'(\xi^{-2/3})}, \tag{304}$$

which is equivalent to (265) by a change of variables.

Furthermore, if $\mathcal{B}_{\mathrm{br}}(t)$, $\mathcal{B}_{\mathrm{br}}^+(t)$ and $\mathcal{B}_{\mathrm{bm}}(t)$ are defined in the obvious way using a Brownian bridge and a Brownian meander on the interval $[0,t]$, and $Z' \sim \Gamma(1/2)$, then (135), (240) and (194) can be written (since $\mathbb{E}(Z')^{3n/2} = \Gamma(3n/2 + 1/2)/\Gamma(1/2)$)

$$D_n = \frac{\mathbb{E}\big(\sqrt{2}\,\mathcal{B}_{\mathrm{br}}(Z')\big)^n}{n!}, \qquad n \ge 0, \tag{305}$$

$$D_n^+ = \frac{\mathbb{E}\big(\sqrt{2}\,\mathcal{B}_{\mathrm{br}}^+(Z')\big)^n}{n!}, \qquad n \ge 0, \tag{306}$$

$$Q_n = \frac{\mathbb{E}\bigl(\sqrt{2}\,\mathcal{B}_{\mathrm{me}}(Z')\bigr)^n}{n!}, \qquad n \ge 0, \tag{307}$$

and thus $\sum_k D_k z^k$, $\sum_k D_k^+ z^k$, $\sum_k Q_k z^k$ in (137), (242), (196) are the moment generating functions of $\sqrt{2}\,\mathcal{B}_{\mathrm{br}}(Z')$, $\sqrt{2}\,\mathcal{B}_{\mathrm{br}}^+(Z')$ and $\sqrt{2}\,\mathcal{B}_{\mathrm{me}}(Z')$, respectively, interpreted as (divergent) formal power series. Since $Z'$ has the density $\pi^{-1/2}t^{-1/2}e^{-t}$, $t > 0$, the Laplace transforms of these variables are given by, using (370), (371), (383) below, for $\xi > 0$,

$$\mathbb{E}\, e^{-\xi\sqrt{2}\,\mathcal{B}_{\mathrm{br}}(Z')} = \int_0^\infty e^{-t}\psi_{\mathrm{br}}\bigl(\sqrt{2}\,\xi t^{3/2}\bigr)\frac{\mathrm{d}t}{\sqrt{\pi t}} = -\xi^{-1/3}\frac{\mathrm{Ai}(\xi^{-2/3})}{\mathrm{Ai}'(\xi^{-2/3})}, \tag{308}$$

$$\mathbb{E}\, e^{-\xi\sqrt{2}\,\mathcal{B}_{\mathrm{br}}^+(Z')} = \int_0^\infty e^{-t}\psi_{\mathrm{br}}^+\bigl(\sqrt{2}\,\xi t^{3/2}\bigr)\frac{\mathrm{d}t}{\sqrt{\pi t}} = 2\frac{\mathrm{Ai}(\xi^{-2/3})}{\mathrm{Ai}(\xi^{-2/3}) - \xi^{1/3}\mathrm{Ai}'(\xi^{-2/3})}, \tag{309}$$

$$\mathbb{E}\, e^{-\xi\sqrt{2}\,\mathcal{B}_{\mathrm{me}}(Z')} = \int_0^\infty e^{-t}\psi_{\mathrm{me}}\bigl(\sqrt{2}\,\xi t^{3/2}\bigr)\frac{\mathrm{d}t}{\sqrt{\pi t}} = \xi^{-1/3}\frac{\mathrm{AI}(\xi^{-2/3})}{\mathrm{Ai}(\xi^{-2/3})}, \tag{310}$$

which are equivalent to (133) and (134), (239) and (193).

The reason for the introduction of $Z'$ is that, as is well-known and explained in greater detail by Perman and Wellner [38], if $g_Z := \max\{t \le Z : B(t) = 0\}$, the last zero of $B$ before $Z$, then conditioned on $g_Z$, the restrictions of $B$ to the intervals $[0, g_Z]$ and $[g_Z, Z]$ form a Brownian bridge and a Brownian meander on these intervals, and the lengths $g_Z$ and $Z - g_Z$ of these intervals both have the distribution of $Z'$. Moreover, the restrictions to these two intervals are independent (also unconditionally, with $g_Z$ and $Z - g_Z$ independent), and thus $\mathcal{B}_{\mathrm{bm}}(Z) = \mathcal{B}_{\mathrm{br}}(g_Z) + \mathcal{B}_{\mathrm{me}}(Z - g_z)$, where the two variables on the right hand side are independent with the distributions of $\mathcal{B}_{\mathrm{br}}(Z')$ and $\mathcal{B}_{\mathrm{me}}(Z')$, which immediately yields (291). Considering only the positive part of $B$, we similarly obtain (298) [38].

For $\mathcal{B}_{\mathrm{ex}}$, the slightly different (6) can be written

$$K_n = \frac{\mathbb{E}\bigl(\sqrt{2}\,\mathcal{B}_{\mathrm{ex}}(Z')\bigr)^n}{(6n-2)\,n!}, \qquad n \ge 0. \tag{311}$$

Hence, the proper analogue of $L_n$, $L_n^+$, $D_n$, $D_n^+$ and $Q_n$ is $(6n-2)K_n$. Note that it follows from (101) that, as $x \to \infty$,

$$\sum_{k=0}^\infty (-1)^n(6n-2)K_n x^{-1/2-3n/2} \sim -2\frac{\mathrm{d}}{\mathrm{d}x}\left(\frac{\mathrm{Ai}'(x)}{\mathrm{Ai}(x)}\right) = 2\left(\frac{\mathrm{Ai}'(x)}{\mathrm{Ai}(x)}\right)^2 - 2x, \tag{312}$$

the Laplace transform appearing in (84); further, in analogy with (308)–(310), (84) is equivalent to, for $\xi > 0$,

$$\mathbb{E}\, e^{-\xi\sqrt{2}\,\mathcal{B}_{\mathrm{ex}}(Z')} = \int_0^\infty e^{-t}\psi_{\mathrm{ex}}\bigl(\sqrt{2}\,\xi t^{3/2}\bigr)\frac{\mathrm{d}t}{\sqrt{\pi t}} = 2\xi^{-1/3}\left(\frac{\mathrm{Ai}'(\xi^{-2/3})}{\mathrm{Ai}(\xi^{-2/3})}\right)^2 - 2\xi^{-1}. \tag{313}$$

For $\mathcal{B}_{\mathrm{dm}}$ we similarly let $B_{\mathrm{dm}}^{(t)}$ be a Brownian double meander on $[0,t]$. Then $B_{\mathrm{dm}}^{(Z)}$ has a.s. a unique zero at a time $\mu=\tau Z$ where $\tau\sim\mathrm{Beta}(\frac12,\frac12)$ is independent of $Z$; hence $\mu$ and $Z-\mu=(1-\tau)Z$ are two independent random variables with a $\Gamma(1/2)$ distribution. Moreover, see Section 23, given $\mu$ and $Z-\mu$, $B_{\mathrm{dm}}^{(Z)}$ consists of two independent Brownian meanders on the intervals $[0,\mu]$ and $[\mu,Z]$, with the first one reversed. Consequently, the Brownian double meander area $\mathcal{B}_{\mathrm{dm}}(Z):=\int_0^Z B_{\mathrm{dm}}^{(Z)}(t)\,\mathrm{d}t$ equals the sum of two independent Brownian meander areas, both having the same distribution as $\mathcal{B}_{\mathrm{me}}(Z')$. Hence, by (310),

$$\mathbb{E}\,e^{-\xi\sqrt2\,\mathcal{B}_{\mathrm{dm}}(Z)}=\Bigl(\mathbb{E}\,e^{-\xi\sqrt2\,\mathcal{B}_{\mathrm{me}}(Z')}\Bigr)^2=\xi^{-2/3}\left(\frac{\mathrm{AI}(\xi^{-2/3})}{\mathrm{Ai}(\xi^{-2/3})}\right)^2, \tag{314}$$

which is equivalent to (221) by a change of variables. Furthermore,

$$W_n=\frac{\mathbb{E}\bigl(\sqrt2\,\mathcal{B}_{\mathrm{dm}}(Z)\bigr)^n}{n!},\qquad n\ge0, \tag{315}$$

and the generating function $\sum_{k=0}^\infty W_k z^k$ in (231) is the moment generating function of $\sqrt2\,\mathcal{B}_{\mathrm{dm}}$ interpreted as a (divergent) formal power series.

## 28. A comparison

For easy reference, we collect in Table 10 the first two moments of the various Brownian areas treated above, together with the scale invariant ratio of the second moment and the square of the first.

Note that the variables $\mathcal{B}_{\mathrm{ex}}$, $\mathcal{B}_{\mathrm{dm}}$, $\mathcal{B}_{\mathrm{me}}$, $\mathcal{B}_{\mathrm{br}}$, $\mathcal{B}_{\mathrm{bm}}$ have quite small variances; if these variables are normalized to have means 1, their variances range from 0.061 to 0.325 (in this order, for which we see no intuitive reason), which means that these variables, and in particular $\mathcal{B}_{\mathrm{ex}}$, typically do not vary much from their means.

## 29. Negative moments

We gave in Section 18 formulas for negative moments of $\mathcal{B}_{\mathrm{ex}}$ due to Flajolet and Louchard [17]. These can be generalized using the results in Appendix B. First, (339) and (340) yield by (84), (133), (173), (193), (221), (239) and (265) expressions for the negative moments of order an odd multiple of $-\frac13$ for $\mathcal{B}_{\mathrm{ex}}$, $\mathcal{B}_{\mathrm{br}}$, $\mathcal{B}_{\mathrm{me}}$, $\mathcal{B}_{\mathrm{br}}^{+}$, and a multiple of $-\frac23$ for $\mathcal{B}_{\mathrm{bm}}$, $\mathcal{B}_{\mathrm{dm}}$, $\mathcal{B}_{\mathrm{bm}}^{+}$. For example, this yields the values in Table 11 (including a repetition of (124) and (125)); note that (126) enables us to rewrite the values in Table 11 in other forms.

More generally, by (341) and (84), (133), (173), (193), (221), where the right hand side is analytic at 0, the random variables $\mathcal{B}_{\mathrm{ex}}$, $\mathcal{B}_{\mathrm{br}}$, $\mathcal{B}_{\mathrm{bm}}$, $\mathcal{B}_{\mathrm{me}}$, $\mathcal{B}_{\mathrm{dm}}$ have negative moments of all orders; thus, $\mathbb{E}\,\mathcal{B}_{\mathrm{ex}}^s$, $\mathbb{E}\,\mathcal{B}_{\mathrm{br}}^s$, $\mathbb{E}\,\mathcal{B}_{\mathrm{bm}}^s$, $\mathbb{E}\,\mathcal{B}_{\mathrm{me}}^s$, $\mathbb{E}\,\mathcal{B}_{\mathrm{dm}}^s$ are entire functions of $s$.

| $X$ | $\mathbb{E}\,X$ | $\mathbb{E}\,X^2$ | $\mathbb{E}\,X^2/(\mathbb{E}\,X)^2$ |
|---|---|---|---|
| $\mathcal{B}_{\mathrm{ex}}$ | $\frac{\sqrt{2\pi}}{4}$ | $\frac{5}{12}$ | $\frac{10}{3\pi}\approx 1.06103$ |
| $\mathcal{B}_{\mathrm{br}}$ | $\frac{\sqrt{2\pi}}{8}$ | $\frac{7}{60}$ | $\frac{56}{15\pi}\approx 1.18836$ |
| $\mathcal{B}_{\mathrm{bm}}$ | $\frac{2\sqrt{2}}{3\sqrt{\pi}}$ | $\frac{3}{8}$ | $\frac{27\pi}{64}\approx 1.32536$ |
| $\mathcal{B}_{\mathrm{me}}$ | $\frac{3\sqrt{2\pi}}{8}$ | $\frac{59}{60}$ | $\frac{472}{135\pi}\approx 1.11291$ |
| $\mathcal{B}_{\mathrm{dm}}$ | $\frac{\sqrt{2}}{\sqrt{\pi}}$ | $\frac{17}{24}$ | $\frac{17\pi}{48}\approx 1.11265$ |
| $\mathcal{B}^+_{\mathrm{bm}}$ | $\frac{\sqrt{2}}{3\sqrt{\pi}}$ | $\frac{17}{96}$ | $\frac{51\pi}{64}\approx 2.50346$ |
| $\mathcal{B}^+_{\mathrm{br}}$ | $\frac{\sqrt{2\pi}}{16}$ | $\frac{1}{20}$ | $\frac{32}{5\pi}\approx 2.03718$ |

TABLE 10
*The first two moments for seven Brownian areas.*

For $\mathcal{B}^+_{\mathrm{br}}$ we see that the right hand side of (239) is finite at $x = 0$, but that its derivative is $\sim cx^{-1/2}$ as $x \to 0$ for some constant $c \neq 0$. Hence, (341) with $m = 1$, $\nu = 1/2$ and $\Psi(x)$ equal to the right hand side of (239), shows that for $0 < s < 3/2$, the moment $\mathbb{E}(\mathcal{B}^+_{\mathrm{br}})^{-1+2s/3}$ is finite if and only $s > 1/2$. Consequently, the negative moment $\mathbb{E}(\mathcal{B}^+_{\mathrm{br}})^{-s}$ is finite if and only if $s < 2/3$, and $\mathbb{E}(\mathcal{B}^+_{\mathrm{br}})^s$ is an analytic function in the half-plane $\operatorname{Re} s > -\frac{2}{3}$. More precisely, a differentiation of (239) shows that, for this $\Psi$,

$$-\Psi'(x) = \sqrt{\pi}\Big(\frac{\operatorname{Ai}(0)}{\operatorname{Ai}'(0)}\Big)^2 x^{-1/2} + O(1), \qquad x \to 0, \tag{316}$$

and thus the left hand side of (341) with $m = 1$ has the residue $\sqrt{\pi}\big(\operatorname{Ai}(0)/\operatorname{Ai}'(0)\big)^2$ at $s = 1/2$, and consequently (341) yields

$$\operatorname{Res}_{s=-\frac{2}{3}} \mathbb{E}(\mathcal{B}^+_{\mathrm{br}})^s = \frac{2^{1/3}}{\Gamma(2/3)}\Big(\frac{\operatorname{Ai}(0)}{\operatorname{Ai}'(0)}\Big)^2 = \frac{2^{1/3}3^{-2/3}\Gamma(1/3)^2}{\Gamma(2/3)^3} = \frac{2^{1/3}3^{5/6}\Gamma(1/3)^5}{(2\pi)^3}. \tag{317}$$

For $\mathcal{B}^+_{\mathrm{bm}}$ we see that the right hand side of (265) is $\sim cx^{-1/2}$ as $x \to 0$ for some constant $c > 0$. Hence, (338) with $\nu = 1$ and $\Psi(x)$ equal to the right hand side of (265) shows that for $0 < s < 1$, the moment $\mathbb{E}(\mathcal{B}^+_{\mathrm{bm}})^{-2(1-s)/3}$ is finite if and only $s > 1/2$. Consequently, the negative moment $\mathbb{E}(\mathcal{B}^+_{\mathrm{bm}})^{-s}$ is finite if and only if $s < 1/3$, and $\mathbb{E}(\mathcal{B}^+_{\mathrm{bm}})^s$ is an analytic function in the half-plane $\operatorname{Re} s > -\frac{1}{3}$. (This also follows from the considerations in Section 27, which imply that, for $s > 0$, $\mathbb{E}(\mathcal{B}^+_{\mathrm{bm}}(Z))^{-s}$ is finite if and only if $\mathbb{E}(\mathcal{B}^+_{\mathrm{br}}(Z'))^{-s} = \mathbb{E}(\mathcal{B}^+_{\mathrm{br}})^{-s}\,\mathbb{E}(Z')^{-3s/2}$ is, and $\mathbb{E}(Z')^{-3s/2} < \infty$ if and only if $s < 1/3$.) More precisely, (265) shows

| | | | |
|---|---|---|---|
| $\mathbb{E}\,\mathcal{B}_{\mathrm{ex}}^{-1/3}$ | $= \dfrac{2^{1/6}3\sqrt{\pi}}{\Gamma(1/3)}\cdot\left(\dfrac{\mathrm{Ai}'(0)}{\mathrm{Ai}(0)}\right)^2$ | $= 2^{-17/6}3^{19/6}\pi^{-5/2}\Gamma(2/3)^5$ | $\approx 1.184$ |
| $\mathbb{E}\,\mathcal{B}_{\mathrm{br}}^{-1/3}$ | $= \dfrac{2^{-5/6}3\sqrt{\pi}}{\Gamma(1/3)}\cdot\dfrac{\mathrm{Ai}(0)}{-\mathrm{Ai}'(0)}$ | $= 2^{-11/6}3^{7/6}\pi^{-1/2}\Gamma(1/3)$ | $\approx 1.528$ |
| $\mathbb{E}\,\mathcal{B}_{\mathrm{me}}^{-1/3}$ | $= \dfrac{2^{-5/6}3\sqrt{\pi}}{\Gamma(1/3)}\cdot\dfrac{\mathrm{AI}(0)}{\mathrm{Ai}(0)}$ | $= 2^{-11/6}3^{7/6}\pi^{-1/2}\Gamma(2/3)^2$ | $\approx 1.046$ |
| $\mathbb{E}(\mathcal{B}_{\mathrm{br}}^{+})^{-1/3}$ | $= \dfrac{2^{1/6}3\sqrt{\pi}}{\Gamma(1/3)}\cdot\dfrac{\mathrm{Ai}(0)}{-\mathrm{Ai}'(0)}$ | $= 2^{-5/6}3^{7/6}\pi^{-1/2}\Gamma(1/3)$ | $\approx 3.056$ |
| $\mathbb{E}\,\mathcal{B}_{\mathrm{bm}}^{-2/3}$ | $= \dfrac{2^{-2/3}3}{\Gamma(2/3)}\cdot\dfrac{\mathrm{AI}(0)}{-\mathrm{Ai}'(0)}$ | $= 2^{-2/3}3^{1/3}\dfrac{\Gamma(1/3)}{\Gamma(2/3)}$ | $\approx 1.797$ |
| $\mathbb{E}\,\mathcal{B}_{\mathrm{dm}}^{-2/3}$ | $= \dfrac{2^{-2/3}3}{\Gamma(2/3)}\cdot\left(\dfrac{\mathrm{AI}(0)}{\mathrm{Ai}(0)}\right)^2$ | $= 2^{-2/3}3^{1/3}\Gamma(2/3)$ | $\approx 1.230$ |
| $\mathbb{E}(\mathcal{B}_{\mathrm{bm}}^{+})^{-2/3}$ | $= \infty$ | | |
| $\mathbb{E}\,\mathcal{B}_{\mathrm{ex}}^{-1}$ | $= 3\sqrt{2\pi}\left(1-\dfrac{3^{5/2}\Gamma(2/3)^6}{4\pi^3}\right)$ | | $\approx 1.693$ |
| $\mathbb{E}\,\mathcal{B}_{\mathrm{br}}^{-1}$ | $= 2^{-1/2}3\sqrt{\pi}$ | | $\approx 3.760$ |
| $\mathbb{E}\,\mathcal{B}_{\mathrm{me}}^{-1}$ | $= 2^{-1/2}3\sqrt{\pi}\left(1-\dfrac{3^{1/2}\Gamma(2/3)^3}{2\pi}\right)$ | | $\approx 1.186$ |
| $\mathbb{E}(\mathcal{B}_{\mathrm{br}}^{+})^{-1}$ | $= \infty$ | | |
| $\mathbb{E}\,\mathcal{B}_{\mathrm{bm}}^{-4/3}$ | $= 2^{-4/3}3^{13/6}\pi^{-1}\Gamma(1/3)$ | | $\approx 3.658$ |
| $\mathbb{E}\,\mathcal{B}_{\mathrm{dm}}^{-4/3}$ | $= 2^{-1/3}3^{13/6}\pi^{-1}\Gamma(2/3)^2\left(1-\dfrac{3^{1/2}\Gamma(2/3)^3}{2\pi}\right)$ | | $\approx 1.580$ |
| $\mathbb{E}\,\mathcal{B}_{\mathrm{bm}}^{-2}$ | $= 3+\dfrac{3^{1/2}\Gamma(1/3)^3}{2\pi}$ | | $\approx 8.300$ |
| $\mathbb{E}\,\mathcal{B}_{\mathrm{dm}}^{-2}$ | $= 6-\dfrac{3^{5/2}\Gamma(2/3)^3}{\pi}+\dfrac{27\,\Gamma(2/3)^6}{2\pi^2}$ | | $\approx 2.112$ |

TABLE 11
*Some negative fractional moments.*

that, for this $\Psi$,

$$\Psi(x) = -\frac{\mathrm{Ai}(0)}{\mathrm{Ai}'(0)} x^{-1/2} + O(1), \qquad x \to 0, \tag{318}$$

and thus the left hand side of (338) has the residue $-\mathrm{Ai}(0)/\mathrm{Ai}'(0)$ at $s = 1/2$, and consequently (338) yields

$$\begin{aligned} \operatorname{Res}_{s=-\frac{1}{3}} \mathbb{E}(\mathcal{B}_{\mathrm{bm}}^{+})^{s} &= \frac{2^{1/6}}{\Gamma(1/2)\Gamma(1/3)} \frac{\mathrm{Ai}(0)}{-\mathrm{Ai}'(0)} = \frac{2^{1/6}3^{-1/3}}{\sqrt{\pi}\,\Gamma(2/3)} \\ &= 2^{-5/6}3^{1/6}\pi^{-3/2}\Gamma(1/3). \end{aligned} \tag{319}$$

Furthermore, the right hand sides of (239) and (265) have expansions in powers of $x^{1/2}$ at $x = 0$ (beginning with $x^{-1/2}$ for (265)). Hence, in (341) with $\Psi(x)$ equal to one of these functions, the left hand side extends to a meromorphic function in $0 < \operatorname{Re} s < m+\nu$ with possible poles (all simple) only at $\frac{1}{2}, \frac{3}{2}, \dots$, and thus the same holds for the negative moment $\mathbb{E}\, X^{-2(m+\nu-s)/3}$ on the right hand side of (341). Since $m$ is an arbitrary non-negative integer and $\nu = 1/2$ for $\mathcal{B}_{\mathrm{br}}^{+}$ while $\nu = 1$ for $\mathcal{B}_{\mathrm{bm}}^{+}$, it follows that $\mathbb{E}(\mathcal{B}_{\mathrm{br}}^{+})^{s}$ extends to a meromorphic function in the complex plane with poles only at $-\frac{2}{3}, -\frac{4}{3}, -2, \dots$, i.e., multiples of $-\frac{2}{3}$, and that $\mathbb{E}(\mathcal{B}_{\mathrm{bm}}^{+})^{s}$ extends to a meromorphic function in the complex plane with poles only at $-\frac{1}{3}, -1, -\frac{5}{3}, \dots$, i.e., odd multiples of $-\frac{1}{3}$. Furthermore, all poles are simple. (We cannot prove that all these points really are poles; it is conceivable that some actually are regular points because some terms in the expansions of $\Psi(x)$ might vanish due to unexpexted cancellations, but we believe that all actually are poles. At least, $-\frac{2}{3}$ and $-\frac{1}{3}$ are poles by the results above.)

Recall also the relation (122) in Section 18 between moments of $\mathcal{B}_{\mathrm{ex}}$ and the root zeta function $\Lambda$ of Ai. For the Brownian bridge, there is an analogous relation between moments of $\mathcal{B}_{\mathrm{br}}$ and the root zeta function of Ai$'$ defined by, cf. (120),

$$\tilde{\Lambda}(s) := \sum_{j=1}^{\infty} |a_j'|^{-s}, \qquad \operatorname{Re} s > 3/2. \tag{320}$$

Since $|a_j'| \sim |a_j| \sim (3\pi j/2)^{2/3}$ [1, 10.4.95] (in fact, $|a_{j-1}| < |a_j'| < |a_j|$ for $j \geq 2$), the sum converges and $\tilde{\Lambda}(s)$ is analytic for $\operatorname{Re} s > 3/2$. By (335) in Appendix B

and (159), for $\operatorname{Re} s > 0$ and with absolutely convergent sums and integrals,

$$
\begin{aligned}
\Gamma(s)\,\mathbb{E}\,\mathcal{B}_{\mathrm{br}}^{-s} &= \int_0^\infty t^{s-1}\psi_{\mathrm{br}}(t)\,\mathrm{d}t \\
&= 2^{-1/6}\sqrt{\pi}\int_0^\infty t^{s+1/3}\sum_{j=1}^\infty |a_j'|^{-1}\exp\bigl(-2^{-1/3}|a_j'|t^{2/3}\bigr)\frac{\mathrm{d}t}{t} \\
&= 2^{-7/6}3\sqrt{\pi}\sum_{j=1}^\infty |a_j'|^{-1}\int_0^\infty u^{(3s+1)/2}\exp\bigl(-2^{-1/3}|a_j'|u\bigr)\frac{\mathrm{d}u}{u} \\
&= 2^{s/2-1}3\sqrt{\pi}\sum_{j=1}^\infty |a_j'|^{-(3s+3)/2}\Gamma\Bigl(\frac{3s+1}{2}\Bigr) \\
&= 3\sqrt{\pi}\,2^{s/2-1}\Gamma\Bigl(\frac{3s+1}{2}\Bigr)\tilde{\Lambda}\Bigl(\frac{3s+3}{2}\Bigr).
\end{aligned}
\tag{321}
$$

Since $\mathbb{E}\,\mathcal{B}_{\mathrm{br}}^{s}$ is an entire function, this shows that $\tilde{\Lambda}$ extends to a meromorphic function in the complex plane with, similarly to (122),

$$
\mathbb{E}\,\mathcal{B}_{\mathrm{br}}^{s} = \frac{3\sqrt{\pi}\,2^{-s/2-1}}{\Gamma(-s)}\Gamma\Bigl(\frac{1-3s}{2}\Bigr)\tilde{\Lambda}\Bigl(\frac{3-3s}{2}\Bigr), \qquad s \in \mathbb{C}. \tag{322}
$$

Since the left hand side has no poles and no real zeros, it follows, arguing as in Section 18, that $\tilde{\Lambda}(z)$ has simple poles at $z = \frac{3}{2}, -\frac{3}{2}, -\frac{9}{2}, \dots$, and zeros at 1, $-1, -2, -4, -5, -7, \dots$, and nowhere else.

It follows from (158) that, for $|z| < |a_1'|$,

$$
\frac{\mathrm{Ai}(z)}{\mathrm{Ai}'(z)} = \sum_{m\ge 0}(-1)^{m+1}\tilde{\Lambda}(m+2)z^m, \tag{323}
$$

and thus the values of $\tilde{\Lambda}(s)$ at positive integers $s = 2, 3, \dots$ can be computed from the Taylor series of Ai. (In particular, $\tilde{\Lambda}(2) = -\mathrm{Ai}(0)/\mathrm{Ai}'(0) = 3^{-1/3}\Gamma(1/3)/\Gamma(2/3)$ and $\tilde{\Lambda}(3) = 1$. Furthermore, as just said, $\tilde{\Lambda}(1) = 0$.) This and (322) give again the same explicit formulas that we have obtained from (340) for the negative moments $\mathcal{B}_{\mathrm{br}}^{-s}$ when $s$ is an odd multiple of $1/3$.

Note also that by (133), (338) and (321), we have the Mellin transform, for $0 < \operatorname{Re} s < 1/2$,

$$
\begin{aligned}
\int_0^\infty x^{s-1}\frac{\mathrm{Ai}(x)}{\mathrm{Ai}'(x)}\,\mathrm{d}x &= -\frac{2^{5/6+s/3}}{3\sqrt{\pi}}\Gamma(s)\Gamma\Bigl(\frac{1-2s}{3}\Bigr)\,\mathbb{E}\,\mathcal{B}_{\mathrm{br}}^{(2s-1)/3} \\
&= -\Gamma(s)\Gamma(1-s)\tilde{\Lambda}(2-s) \\
&= -\frac{\pi}{\sin(\pi s)}\tilde{\Lambda}(2-s),
\end{aligned}
\tag{324}
$$

which also follows from (158); cf. (127) and (128).

For Brownian motion, Brownian meander and Brownian double meander, more complicated analogues of the root zeta function $\Lambda(s)$ appear in the same way: for Brownian motion, using (170), $\sum_{j=1}^{\infty} \kappa_j |a'_j|^{-s}$; for Brownian meander, using (209), $\sum_{j=1}^{\infty} r_j |a_j|^{-s}$; for Brownian double meander, using (220), $\sum_{j=1}^{\infty} r_j^2 |a_j|^{-s}$. We do not pursue this further.

## Appendix A: The integrated Airy function

We define, as in (78), using $\int_0^\infty \mathrm{Ai}(x)\,\mathrm{d}x = 1/3$ [1, 10.4.82],

$$\mathrm{AI}(z) := \int_z^{+\infty} \mathrm{Ai}(t)\,\mathrm{d}t = \frac{1}{3} - \int_0^z \mathrm{Ai}(t)\,\mathrm{d}t; \tag{325}$$

this is well-defined for all complex $z$ and yields an entire function provided the first integral is taken along, for example, a path that eventually follows the positive real axis to $+\infty$. Note that $\mathrm{AI}'(z) = -\mathrm{Ai}(z)$ and that $\mathrm{AI}(0) = 1/3$. Along the real axis we have the limits, by [1, 10.4.82–83],

$$\lim_{x \to +\infty} \mathrm{AI}(x) = 0, \tag{326}$$

$$\lim_{x \to -\infty} \mathrm{AI}(x) = \int_{-\infty}^{\infty} \mathrm{Ai}(x)\,\mathrm{d}x = 1. \tag{327}$$

In terms of the functions Gi and Hi defined in [1], we have, see [1, 10.4.47–48],

$$\mathrm{AI}(z) = \pi\big(\mathrm{Ai}(z)\mathrm{Gi}'(z) - \mathrm{Ai}'(z)\mathrm{Gi}(z)\big) \tag{328}$$

$$= 1 + \pi\big(\mathrm{Ai}'(z)\mathrm{Hi}(z) - \mathrm{Ai}(z)\mathrm{Hi}'(z)\big). \tag{329}$$

Repeated integrations by parts give

$$\begin{aligned} \mathrm{AI}(z) &= \int_z^{+\infty} \mathrm{Ai}(w)\,\mathrm{d}w = \int_z^{+\infty} \frac{1}{w}\mathrm{Ai}''(w)\,\mathrm{d}w \\ &= -\frac{\mathrm{Ai}'(z)}{z} + \int_z^{+\infty} \frac{1}{w^2}\mathrm{Ai}'(w)\,\mathrm{d}w \\ &= -\frac{\mathrm{Ai}'(z)}{z} - \frac{\mathrm{Ai}(z)}{z^2} + \int_z^{+\infty} \frac{2}{w^3}\mathrm{Ai}(w)\,\mathrm{d}w \\ &= -\frac{\mathrm{Ai}'(z)}{z} - \frac{\mathrm{Ai}(z)}{z^2} + \int_z^{+\infty} \frac{2}{w^4}\mathrm{Ai}''(w)\,\mathrm{d}w = \dots; \end{aligned} \tag{330}$$

if we assume $|\arg(z)| < \pi - \delta$ for some $\delta > 0$, we can here integrate for example along the horizontal line $\{z + t : t \geq 0\}$, and it follows from the asymptotic expansions (102) and (103), which are valid for $|\arg(z)| < \pi - \delta$ [1, 10.4.59 and 10.4.61] that there is a similar asymptotic expansion, as $|z| \to \infty$,

$$\mathrm{AI}(z) \sim \frac{1}{2\sqrt{\pi}} z^{-3/4} e^{-2z^{3/2}/3} \sum_{k=0}^{\infty} (-1)^k \beta_k z^{-3k/2}, \tag{331}$$

valid for $|\arg(z)| < \pi - \delta$ and thus extending (183) (which is given in Takács [49]; see also [1, 10.4.82]). The coefficients $\beta_k$ can be found by this procedure, but it is easier to observe that we can formally differentiate (331) and obtain an asymptotic expansion of $-\mathrm{Ai}(z)$, and a comparison with (102) yields the recursion relation (184). See Table 3 for a few numerical values.

Similarly, along the negative real axis, and more generally for $-z$ where $|\arg(z)| < 2\pi/3 - \delta$ for some $\delta > 0$, there are asymptotic expansions [1, 10.4.60 and 10.4.62]

$$\mathrm{Ai}(-z) \sim \pi^{-1/2}\Bigl(\sin\Bigl(\frac{2}{3}z^{3/2}+\frac{\pi}{4}\Bigr)\sum_{k=0}^{\infty}(-1)^k\alpha_{2k}z^{-3k-1/4} - \cos\Bigl(\frac{2}{3}z^{3/2}+\frac{\pi}{4}\Bigr)\sum_{k=0}^{\infty}(-1)^k\alpha_{2k+1}z^{-3k-7/4}\Bigr), \tag{332}$$

$$\mathrm{Ai}'(-z) \sim -\pi^{-1/2}\Bigl(\cos\Bigl(\frac{2}{3}z^{3/2}+\frac{\pi}{4}\Bigr)\sum_{k=0}^{\infty}(-1)^k\alpha'_{2k}z^{-3k+1/4} + \sin\Bigl(\frac{2}{3}z^{3/2}+\frac{\pi}{4}\Bigr)\sum_{k=0}^{\infty}(-1)^k\alpha'_{2k+1}z^{-3k-5/4}\Bigr), \tag{333}$$

which using (330) lead to, as $|z| \to \infty$ with $|\arg(z)| < 2\pi/3 - \delta$,

$$\mathrm{AI}(-z) \sim 1 - \pi^{-1/2}\Bigl(\cos\Bigl(\frac{2}{3}z^{3/2}+\frac{\pi}{4}\Bigr)\sum_{k=0}^{\infty}(-1)^k\beta_{2k}z^{-3k-3/4} + \sin\Bigl(\frac{2}{3}z^{3/2}+\frac{\pi}{4}\Bigr)\sum_{k=0}^{\infty}(-1)^k\beta_{2k+1}z^{-3k-9/4}\Bigr); \tag{334}$$

given (at least along the negative real axis) by Takács [49, (17)] (using his $h_k = (2/3)^k\beta_k$); see also [1, 10.4.83].

## Appendix B: The Mellin transform of a Laplace transform

For any positive random variable $X$ with Laplace transform $\psi(t) := \mathbb{E}\, e^{-tX}$, and any $s > 0$, Fubini's theorem yields

$$\int_0^\infty t^{s-1}\psi(t)\,\mathrm{d}t = \mathbb{E}\int_0^\infty t^{s-1}e^{-tX}\,\mathrm{d}t = \Gamma(s)\,\mathbb{E}\,X^{-s}. \tag{335}$$

If one side is finite for $s = s_0 > 0$, and thus both sides are, then $\mathbb{E}\,X^{-s} < \infty$ for $s = s_0$ and thus for $0 \le s \le s_0$; hence both sides of (335) are finite and (335) holds for all $s$ with $0 < s \le s_0$, and more generally for all complex $s$ with $0 < \mathrm{Re}(s) \le s_0$, with the terms in (335) analytic in the strip $0 < \mathrm{Re}(s) < s_0$.

**Lemma B.1.** *The following are equivalent:*

(i) *$X$ has negative moments of all orders, i.e., $\mathbb{E}\,X^{-s} < \infty$ for all $s > 0$;*

(ii) $\int_0^\infty t^{s-1}\psi(t)\,\mathrm{d}t < \infty$ *for all* $s > 0$*;*
(iii) $\psi(t) = O(t^{-N})$ *as* $t \to \infty$ *for all* $N$*.*

*If these hold and further* $X$ *has moments of all (positive) orders, then* $\mathbb{E}\,X^{-s}$ *is an entire function of* $s$*; further,* (335) *holds for* $\operatorname{Re} s > 0$ *and the right hand side defines a meromorphic extension of* $\int_0^\infty t^{s-1}\psi(t)\,\mathrm{d}t$ *to the complex plane.*

*Proof.* (i) $\iff$ (ii) by (335).

Since $\psi(t) \le \psi(0) = 1$, (iii) $\Longrightarrow$ (ii).

Conversely, since $\psi$ is decreasing, (ii) implies for every $u > 0$

$$\infty > \int_0^\infty t^{s-1}\psi(t)\,\mathrm{d}t \ge \psi(u)\int_0^u t^{s-1}\,\mathrm{d}t = s^{-1}u^s\psi(u), \tag{336}$$

and thus $\psi(u) = O(u^{-s})$.

The last statement is obvious. □

Several formulas above, viz. (84), (133), (173), (193), (221), (239), (265), give a formula for a double Laplace transform of the form

$$\int_0^\infty e^{-xt}\psi\big(\sqrt{2}\,t^{3/2}\big)t^{\nu-1}\,\mathrm{d}t = \Psi(x), \qquad x > 0, \tag{337}$$

where $\nu = 1$ or $1/2$; see also Section 27 for an explanation of this form. (We have $\nu = 1$ for $\mathcal{B}_{\mathrm{bm}}$ and $\mathcal{B}_{\mathrm{bm}}^+$ and $\nu = 1/2$ for the others.) We then can argue as above twice (using (335) the second time) and obtain for $0 < s < \nu$, and if the result then is finite also for complex $s$ with $0 < \operatorname{Re} s < \nu$,

$$\begin{aligned}\int_0^\infty x^{s-1}\Psi(x)\,\mathrm{d}x &= \Gamma(s)\int_0^\infty t^{-s}\psi\big(\sqrt{2}\,t^{3/2}\big)t^{\nu-1}\,\mathrm{d}t\\ &= \Gamma(s)\frac{2}{3}\int_0^\infty 2^{-(\nu-s)/3}u^{2(\nu-s)/3-1}\psi(u)\,\mathrm{d}u\\ &= \frac{2^{1-(\nu-s)/3}}{3}\Gamma(s)\Gamma\Big(\frac{2(\nu-s)}{3}\Big)\,\mathbb{E}\,X^{-2(\nu-s)/3}.\end{aligned} \tag{338}$$

Furthermore, $x = 0$ in (337) similarly yields, with $\Psi(0) := \lim_{x\to 0}\Psi(x) \le \infty$,

$$\begin{aligned}\Psi(0) &= \int_0^\infty \psi\big(\sqrt{2}\,t^{3/2}\big)t^{\nu-1}\,\mathrm{d}t\\ &= \frac{2}{3}\int_0^\infty 2^{-\nu/3}u^{2\nu/3-1}\psi(u)\,\mathrm{d}u\\ &= \frac{2^{1-\nu/3}}{3}\Gamma\Big(\frac{2\nu}{3}\Big)\,\mathbb{E}\,X^{-2\nu/3},\end{aligned} \tag{339}$$

which can be regarded as a limiting case of (338) for $s \to 0$. In particular, $\Psi(0) < \infty$ if and only if $\mathbb{E}\,X^{-2\nu/3} < \infty$ (which holds for all Brownian areas treated here except $\mathcal{B}_{\mathrm{bm}}^+$).

More generally, by first taking the $m$:th derivative of (337), we obtain for every integer $m \geq 0$, by letting $x \to 0$ and with $\Psi^{(m)}(0) := \lim_{x\to 0} \Psi^{(m)}(x)$,

$$
\begin{aligned}
(-1)^m \Psi^{(m)}(0) &= \int_0^\infty \psi\bigl(\sqrt{2}\, t^{3/2}\bigr) t^{m+\nu-1}\,\mathrm{d}t \\
&= \frac{2^{1-(m+\nu)/3}}{3} \Gamma\Bigl(\frac{2(m+\nu)}{3}\Bigr)\,\mathbb{E}\, X^{-2m/3-2\nu/3}
\end{aligned}
\tag{340}
$$

and, for $0 < s < m+\nu$ (and for $0 < \operatorname{Re} s < m+\nu$ if the result is finite),

$$
\begin{aligned}
(-1)^m &\int_0^\infty x^{s-1} \Psi^{(m)}(x)\,\mathrm{d}x \\
&= \frac{2^{1-(m+\nu-s)/3}}{3} \Gamma(s) \Gamma\Bigl(\frac{2(m+\nu-s)}{3}\Bigr)\,\mathbb{E}\, X^{-2(m+\nu-s)/3}.
\end{aligned}
\tag{341}
$$

## Appendix C: Feynman-Kac formulas

Many of the results above are derived from formulas obtained by simple applications of the Feynman–Kac formula, following Kac [30; 31]. This is by now standard, but for completeness we give here a unified treatment of all the Brownian functionals discussed above. We basically follow the proof by Shepp [42] for the Brownian bridge. See also Jeanblanc, Pitman and Yor [27] for a treatment of several related Brownian functionals.

One version of the Feynman–Kac formula is the following. Let $B_x(t)$ denote the Brownian motion started at $x$; thus $B_x(t) = x + B(t)$; we continue to use $B(t)$ for $B_0(t)$.

**Lemma C.1.** *Let $V(x) \geq 0$ be a non-negative continuous function on $(-\infty,\infty)$, let $\lambda > 0$, and let $\phi_+$ and $\phi_-$ be $C^2$ solutions of the differential equation*

$$
\tfrac{1}{2}\phi''(x) = \bigl(V(x)+\lambda\bigr)\phi(x) \tag{342}
$$

*such that $\phi_+$ is bounded on $[0,\infty)$ and $\phi_-$ is bounded on $(-\infty,0]$. Let $w := \phi_+(0)\phi_-'(0) - \phi_-(0)\phi_+'(0)$ and assume that $w \neq 0$. Moreover, let $f$ be any bounded function on $\mathbb{R}$. Then, for all real $x$,*

$$
\begin{aligned}
\int_0^\infty & e^{-\lambda t}\,\mathbb{E}\Bigl(e^{-\int_0^t V(B_x(s))\,\mathrm{d}s} f\bigl(B_x(t)\bigr)\Bigr)\,\mathrm{d}t \\
&= \frac{2}{w}\left(\phi_+(x)\int_{-\infty}^x \phi_-(y) f(y)\,\mathrm{d}y + \phi_-(x)\int_x^\infty \phi_+(y) f(y)\,\mathrm{d}y\right).
\end{aligned}
\tag{343}
$$

*Proof.* Although this is well-known, we give a proof for completeness.

Note that if we define the Wronskian

$$
W(x) := \phi_+(x)\phi_-'(x) - \phi_-(x)\phi_+'(x), \tag{344}
$$

then (342) implies that $W'(x) = 0$ and thus $W(x)$ is constant; hence

$$
W(x) = w. \tag{345}
$$

Assume first that $f$ is continuous and has compact support, and define

$$\phi(x) := \phi_+(x) \int_{-\infty}^{x} \phi_-(y) f(y) \, \mathrm{d}y + \phi_-(x) \int_{x}^{\infty} \phi_+(y) f(y) \, \mathrm{d}y. \tag{346}$$

Then $\phi \in C^1(\mathbb{R})$, and differentiation yields

$$\phi'(x) = \phi'_+(x) \int_{-\infty}^{x} \phi_-(y) f(y) \, \mathrm{d}y + \phi'_-(x) \int_{x}^{\infty} \phi_+(y) f(y) \, \mathrm{d}y. \tag{347}$$

Hence also $\phi' \in C^1$, and thus $\phi \in C^2$. A second differentiation yields, using (342),

$$\phi''(x) = 2\big(V(x) + \lambda\big)\phi(x) - W(x)f(x) = 2\big(V(x) + \lambda\big)\phi(x) - wf(x). \tag{348}$$

Moreover, since $f$ has compact support, (346) shows that if $x$ is large enough, then $\phi(x) = \big(\int \phi_- f\big)\phi_+(x)$, and thus $\phi(x)$ is bounded for large $x$; similarly $\phi(x)$ is bounded for $x \to -\infty$. Consequently, $\phi$ is bounded on $\mathbb{R}$.

Let $X(t)$ and $Y(t)$ be the stochastic processes $X(t) := \lambda t + \int_0^t V\big(B_x(s)\big) \, \mathrm{d}s$ and

$$Y(t) := e^{-X(t)} \phi\big(B_x(t)\big) + \frac{w}{2} \int_0^t e^{-X(u)} f\big(B_x(u)\big) \, \mathrm{d}u. \tag{349}$$

A straightforward application of Itô's formula, using (348), shows that

$$\mathrm{d}Y(t) = e^{-X(t)} \phi'\big(B_x(t)\big) \, \mathrm{d}B_x(t), \tag{350}$$

and thus $Y(t)$ is a local martingale. Moreover, since $X_t \geq \lambda t$ and $\phi$ and $f$ are bounded, $Y(t)$ is uniformly bounded, and thus a bounded martingale. Hence,

$$\phi(x) = Y(0) = \mathbb{E}\, Y(\infty), \tag{351}$$

where, recalling that $\phi$ is bounded and $X(t) \geq \lambda t$,

$$Y(\infty) := \lim_{t \to \infty} Y(t) = \frac{w}{2} \int_0^\infty e^{-X(u)} f\big(B_x(u)\big) \, \mathrm{d}u. \tag{352}$$

The result (343) now follows by (351), (352), and (346), under our assumption that $f$ is continuous with compact support.

By a monotone class theorem, for example [20, Theorem A.1], (343) remains true for all bounded measurable $f$ with support in a finite interval $[-A, A]$, for any $A < \infty$. Fixing an $x$ with $\phi_-(x) \neq 0$ and letting $f(y) := \operatorname{sign}(\phi_+(y))\mathbf{1}[x < y < A]$, we see by letting $A \to \infty$ that $\int_x^\infty |\phi_+(y)| \, \mathrm{d}y < \infty$. Hence $\phi_+$ is integrable over each interval $(x, \infty)$, and similarly $\phi_-$ is integrable over each interval $(-\infty, x)$. It follows by dominated convergence, considering $f(x)\mathbf{1}[|x| < n] \to f(x)$, that (343) holds for all bounded measurable $f$ on $\mathbb{R}$. ☐

**Remark C.2.** Majumdar and Comtet [35] have given "physical proofs" of (80) and (209) using path integral techniques. This method is closely related to the Feynman–Kac method as formulated in Lemma C.1, which can be seen as follows, where we argue formally and ignore giving proper technical conditions on $V$ and $f$ for the validity of the argument.

Let $H$ be the differential operator $H\phi := -\frac{1}{2}\phi'' + V\phi$. (This is an unbounded positive self-adjoint operator in $L^2(\mathbb{R})$.) Then (348) can be written

$$H\phi + \lambda\phi = \frac{w}{2} f. \tag{353}$$

Hence, the right hand side of (343), which is $\frac{2}{w}\phi(x)$ by (346), equals $(H + \lambda)^{-1} f(x)$. Further, $(H+\lambda)^{-1} f = \int_0^\infty e^{-\lambda t} e^{-tH} f \,\mathrm{d}t$, so by the uniqueness theorem for Laplace transforms, (343) is equivalent to

$$\mathbb{E}\Big(e^{-\int_0^t V(B_x(s))\,\mathrm{d}s} f\big(B_x(t)\big)\Big) = e^{-tH} f(x), \tag{354}$$

which essentially is the path integral formula used by Majumdar and Comtet [35] in their proof.

We let $\xi, \eta > 0$ and choose

$$V(x) := \begin{cases} \sqrt{2}\,\xi x, & x \ge 0, \\ \sqrt{2}\,\eta |x|, & x < 0. \end{cases} \tag{355}$$

It is easily checked that $\mathrm{Ai}\big(\sqrt{2}\,\xi^{1/3} x + \xi^{-2/3}\lambda\big)$ solves (342) on $(0,\infty)$ and is bounded there; thus this can be extended to a solution $\phi_+$ on $\mathbb{R}$ that is bounded on $[0,\infty)$. Similarly, $\mathrm{Ai}\big(-\sqrt{2}\,\eta^{1/3} x + \eta^{-2/3}\lambda\big)$ on $(-\infty, 0)$ can be extended to a solution $\phi_-$ that is bounded on $(-\infty, 0]$. We thus have $\phi_\pm$ as in Lemma C.1 with

$$\phi_+(x) = \mathrm{Ai}\big(\sqrt{2}\,\xi^{1/3} x + \xi^{-2/3}\lambda\big), \qquad x \ge 0, \tag{356}$$

$$\phi_-(x) = \mathrm{Ai}\big(\sqrt{2}\,\eta^{1/3} |x| + \eta^{-2/3}\lambda\big), \qquad x \le 0. \tag{357}$$

As a consequence,

$$w = -\sqrt{2}\,\eta^{1/3} \mathrm{Ai}'(\eta^{-2/3}\lambda)\mathrm{Ai}(\xi^{-2/3}\lambda) - \sqrt{2}\,\xi^{1/3}\mathrm{Ai}'(\xi^{-2/3}\lambda)\mathrm{Ai}(\eta^{-2/3}\lambda). \tag{358}$$

We will now apply Lemma C.1 with this $V$ to the different cases.

### *C.1. Brownian motion*

Recall the definition (262) of $\mathcal{B}^{\pm}_{\mathrm{bm}}$ and define the joint Laplace transform, for $\xi, \eta \ge 0$,

$$\psi^{\pm}_{\mathrm{bm}}(\xi, \eta) := \mathbb{E}\, e^{-\xi \mathcal{B}^+_{\mathrm{bm}} - \eta \mathcal{B}^-_{\mathrm{bm}}}. \tag{359}$$

Thus $\psi_{\mathrm{bm}}(\xi) = \psi^{\pm}_{\mathrm{bm}}(\xi,\xi)$ and $\psi^{+}_{\mathrm{bm}}(\xi) = \psi^{\pm}_{\mathrm{bm}}(\xi, 0)$.

Assume that $\xi, \eta > 0$ and use Lemma C.1 with $V$ as in (355), $\phi_\pm$ as in (356) and (357), $f = 1$ and $x = 0$. By (355) and homogeneity,

$$\int_0^t V(B(s))\,\mathrm{d}s \overset{\mathrm{d}}{=} t^{3/2} \int_0^1 V(B(s))\,\mathrm{d}s = \sqrt{2}\, t^{3/2} \xi \mathcal{B}^+_{\mathrm{bm}} + \sqrt{2}\, t^{3/2} \eta \mathcal{B}^-_{\mathrm{bm}}, \tag{360}$$

and thus (343) yields, recalling (359) and using (356), (357), (358) and (78),

$$\begin{aligned}
&\int_0^\infty e^{-\lambda t}\psi^{\pm}_{\mathrm{bm}}\big(\sqrt{2}\,t^{3/2}\xi,\sqrt{2}\,t^{3/2}\eta\big)\,\mathrm{d}t\\
&=\frac{2}{w}\left(\phi_+(0)\int_{-\infty}^0\phi_-(y)\,\mathrm{d}y+\phi_-(0)\int_0^\infty\phi_+(y)\,\mathrm{d}y\right)\\
&=\frac{2}{w}\left(\mathrm{Ai}(\xi^{-2/3}\lambda)\frac{\mathrm{AI}(\eta^{-2/3}\lambda)}{\sqrt{2}\,\eta^{1/3}}+\mathrm{Ai}(\eta^{-2/3}\lambda)\frac{\mathrm{AI}(\xi^{-2/3}\lambda)}{\sqrt{2}\,\xi^{1/3}}\right)\\
&=\left(\frac{\mathrm{AI}(\xi^{-2/3}\lambda)}{\xi^{1/3}\mathrm{Ai}(\xi^{-2/3}\lambda)}+\frac{\mathrm{AI}(\eta^{-2/3}\lambda)}{\eta^{1/3}\mathrm{Ai}(\eta^{-2/3}\lambda)}\right)\left(-\frac{\xi^{1/3}\mathrm{Ai}'(\xi^{-2/3}\lambda)}{\mathrm{Ai}(\xi^{-2/3}\lambda)}-\frac{\eta^{1/3}\mathrm{Ai}'(\eta^{-2/3}\lambda)}{\mathrm{Ai}(\eta^{-2/3}\lambda)}\right)^{-1}.
\end{aligned}\tag{361}$$

In particular, the special case $\eta=\xi$ yields

$$\int_0^\infty e^{-\lambda t}\psi_{\mathrm{bm}}\big(\sqrt{2}\,t^{3/2}\xi\big)\,\mathrm{d}t=-\frac{\mathrm{AI}(\xi^{-2/3}\lambda)}{\xi^{2/3}\mathrm{Ai}'(\xi^{-2/3}\lambda)}.\tag{362}$$

Similarly, letting $\eta\to0$ we obtain, since $\mathrm{AI}(x)/\mathrm{Ai}(x)\sim x^{-1/2}$ and $\mathrm{Ai}'(x)/\mathrm{Ai}(x)\sim -x^{1/2}$ as $x\to\infty$ by (102), (103) and (183),

$$\begin{aligned}
&\int_0^\infty e^{-\lambda t}\psi^{+}_{\mathrm{bm}}\big(\sqrt{2}\,t^{3/2}\xi\big)\,\mathrm{d}t\\
&\qquad=\left(\frac{\mathrm{AI}(\xi^{-2/3}\lambda)}{\xi^{1/3}\mathrm{Ai}(\xi^{-2/3}\lambda)}+\lambda^{-1/2}\right)\left(-\frac{\xi^{1/3}\mathrm{Ai}'(\xi^{-2/3}\lambda)}{\mathrm{Ai}(\xi^{-2/3}\lambda)}+\lambda^{1/2}\right)^{-1}\\
&\qquad=\xi^{-1/3}\frac{\mathrm{AI}(\xi^{-2/3}\lambda)+\lambda^{-1/2}\xi^{1/3}\mathrm{Ai}(\xi^{-2/3}\lambda)}{-\xi^{1/3}\mathrm{Ai}'(\xi^{-2/3}\lambda)+\lambda^{1/2}\mathrm{Ai}(\xi^{-2/3}\lambda)}\\
&\qquad=\lambda^{-1}\frac{\mathrm{Ai}(\xi^{-2/3}\lambda)+\xi^{-1/3}\lambda^{1/2}\mathrm{AI}(\xi^{-2/3}\lambda)}{\mathrm{Ai}(\xi^{-2/3}\lambda)-\xi^{1/3}\lambda^{-1/2}\mathrm{Ai}'(\xi^{-2/3}\lambda)}.
\end{aligned}\tag{363}$$

Taking $\xi=1$ yields (265).

### *C.2. Brownian bridge*

Let $B^{(t)}_{x,y}$ denote the Brownian bridge on $[0,t]$ with boundary values $B^{(t)}_{x,y}(0)=x$ and $B^{(t)}_{x,y}(t)=y$; it can be defined by conditioning $B_x$ on $\{B_x(t)=y\}$, interpreted in the usual way (with a distribution that is a continuous function of $y$). The standard Brownian bridge $B_{\mathrm{br}}$ equals $B^{(1)}_{0,0}$.

The left hand side of (343) can be written, by conditioning on $B_x(t)$ inside the integral,

$$\int_{t=0}^\infty e^{-\lambda t}\int_{y=-\infty}^\infty \mathbb{E}\,e^{-\int_0^t V(B^{(t)}_{x,y}(s))\,\mathrm{d}s}f(y)e^{-(y-x)^2/(2t)}\frac{\mathrm{d}y}{\sqrt{2\pi t}}\,\mathrm{d}t\tag{364}$$

while the right hand side can be written, with $x \vee y := \max(x,y)$ and $x \wedge y := \min(x,y)$,

$$\frac{2}{w}\int_{-\infty}^{\infty} \phi_+(x \vee y)\phi_-(x \wedge y) f(y)\,\mathrm{d}y. \tag{365}$$

Hence, using also Fubini in (364),

$$\int_0^\infty e^{-\lambda t}\,\mathbb{E}\,e^{-\int_0^t V(B_{x,y}^{(t)}(s))\,\mathrm{d}s} e^{-(y-x)^2/(2t)}\frac{\mathrm{d}t}{\sqrt{2\pi t}} = \frac{2}{w}\phi_+(x \vee y)\phi_-(x \wedge y) \tag{366}$$

for a.e. $y$. However, both sides of (366) are continuous functions of $y$ (the left hand side by dominated convergence because $B_{x,y}^{(t)}(u) \overset{\mathrm{d}}{=} B_{x,0}^{(t)}(u) + yu/t$ in $C[0,t]$ and thus $B_{x,y_n}^{(t)} \overset{\mathrm{d}}{\longrightarrow} B_{x,y}^{(t)}$ in $C[0,t]$ as $y_n \to y$). Hence, (366) holds for every $y \in \mathbb{R}$.

Take $x = y = 0$ and $V$ as in (355). Then, by homogeneity, recalling (236),

$$\int_0^t V\big(B_{0,0}^{(t)}(s)\big)\,\mathrm{d}s \overset{\mathrm{d}}{=} t^{3/2}\int_0^1 V\big(B_{\mathrm{br}}(s)\big)\,\mathrm{d}s = \sqrt{2}\,t^{3/2}\xi\mathcal{B}_{\mathrm{br}}^+ + \sqrt{2}\,t^{3/2}\eta\mathcal{B}_{\mathrm{br}}^-. \tag{367}$$

Define the joint Laplace transform,

$$\psi_{\mathrm{br}}^{\pm}(\xi,\eta) := \mathbb{E}\,e^{-\xi\mathcal{B}_{\mathrm{br}}^+ - \eta\mathcal{B}_{\mathrm{br}}^-}; \tag{368}$$

thus $\psi_{\mathrm{br}}(\xi) = \psi_{\mathrm{br}}^{\pm}(\xi,\xi)$ and $\psi_{\mathrm{br}}^{+}(\xi) = \psi_{\mathrm{br}}^{\pm}(\xi,0)$. Then, (366) and (367) yield, using (356)–(358), for $\xi,\eta \geq 0$,

$$\begin{aligned}\int_0^\infty e^{-\lambda t}\psi_{\mathrm{br}}^{\pm}\big(\sqrt{2}\,t^{3/2}\xi, \sqrt{2}\,t^{3/2}\eta\big)\frac{\mathrm{d}t}{\sqrt{2\pi t}} &= \frac{2}{w}\phi_+(0)\phi_-(0)\\ &= \sqrt{2}\left(-\frac{\xi^{1/3}\mathrm{Ai}'(\xi^{-2/3}\lambda)}{\mathrm{Ai}(\xi^{-2/3}\lambda)} - \frac{\eta^{1/3}\mathrm{Ai}'(\eta^{-2/3}\lambda)}{\mathrm{Ai}(\eta^{-2/3}\lambda)}\right)^{-1}.\end{aligned} \tag{369}$$

In particular, taking $\xi = \eta > 0$,

$$\int_0^\infty e^{-\lambda t}\psi_{\mathrm{br}}\big(\sqrt{2}\,t^{3/2}\xi\big)\frac{\mathrm{d}t}{\sqrt{\pi t}} = -\frac{\mathrm{Ai}(\xi^{-2/3}\lambda)}{\xi^{1/3}\mathrm{Ai}'(\xi^{-2/3}\lambda)}, \tag{370}$$

which by simple changes of variables is equivalent to (133) ($\lambda = x$, $\xi = 1$) and (134) ($\lambda = x$, $\xi = u/\sqrt{2}$). Furthermore, letting $\eta \to 0$ we obtain

$$\int_0^\infty e^{-\lambda t}\psi_{\mathrm{br}}^{+}\big(\sqrt{2}\,t^{3/2}\xi\big)\frac{\mathrm{d}t}{\sqrt{\pi t}} = 2\lambda^{-1/2}\frac{\mathrm{Ai}(\xi^{-2/3}\lambda)}{\mathrm{Ai}(\xi^{-2/3}\lambda) - \xi^{1/3}\lambda^{-1/2}\mathrm{Ai}'(\xi^{-2/3}\lambda)}. \tag{371}$$

Taking $\xi = 1$ yields (239).

### *C.3. Brownian meander*

We apply again Lemma C.1 with $V$ given by (355); this time keeping $\xi > 0$ fixed but letting $\eta \to \infty$; we also assume $x > 0$. By monotone convergence, the left hand side of (343) tends to

$$\begin{aligned}&\int_0^\infty e^{-\lambda t}\,\mathbb{E}\Bigl(e^{-\sqrt{2}\,\xi\int_0^t B_x(s)\,\mathrm{d}s} f\bigl(B_x(t)\bigr);\, B_x(s) \ge 0 \text{ on } [0,t]\Bigr)\,\mathrm{d}t\\ &= \int_0^\infty e^{-\lambda t}\,\mathbb{E}\Bigl(e^{-\sqrt{2}\,\xi\int_0^t B_x(s)\,\mathrm{d}s} f\bigl(B_x(t)\bigr) \bigm| B_x(s) \ge 0 \text{ on } [0,t]\Bigr)\,\mathbb{P}_x(\tau_0 > t)\,\mathrm{d}t,\end{aligned} \tag{372}$$

where $\mathbb{P}_x(\tau_0 > t)$ is the probability that the hitting time $\tau_0 := \inf\{t > 0 : B_x(t) = 0\}$ is larger than $t$, for the Brownian motion $B_x$ started at $x$. For the right hand side of (343), let us assume for simplicity that $f(y) = 0$ for all $y < x$. We note that by (356) and (357), $\phi_+(y)$ does not depend on $\eta$ for $y \ge 0$, while $\phi_-(y)$ does. On the half-line $[0,\infty)$, $\phi_-$ is by (342) governed by the equation $\frac{1}{2}\phi''(x) = \bigl(\sqrt{2}\,\xi x + \lambda\bigr)\phi(x)$; let $\phi_0$ and $\phi_1$ be the solutions to this equation with

$$\phi_0(0) = 1, \qquad \phi_0'(0) = 0, \tag{373}$$

$$\phi_1(0) = 0, \qquad \phi_1'(0) = 1. \tag{374}$$

(Note that these do not depend on $\eta$.) Thus, for $x > 0$,

$$\phi_-(x) = \phi_-(0)\phi_0(x) + \phi_-'(0)\phi_1(x). \tag{375}$$

Furthermore, as $\eta \to \infty$, (357) implies

$$\phi_-(0) \to \mathrm{Ai}(0), \qquad \phi_-'(0) \sim -\sqrt{2}\,\eta^{1/3}\mathrm{Ai}'(0), \tag{376}$$

while, by (358),

$$w \sim -\sqrt{2}\,\eta^{1/3}\mathrm{Ai}'(0)\mathrm{Ai}(\xi^{-2/3}\lambda). \tag{377}$$

As $\eta \to \infty$, we thus have by (375), (376), (377),

$$\frac{\phi_-(x)}{w} \to \frac{\phi_1(x)}{\mathrm{Ai}(\xi^{-2/3}\lambda)} = \frac{\phi_1(x)}{\phi_+(0)}. \tag{378}$$

Returning to (343) we thus have, for $x > 0$ and assuming $f(y) = 0$ for $y < x$, see also [19, Corollary 2.1(ii)],

$$\begin{aligned}&\int_0^\infty e^{-\lambda t}\,\mathbb{E}\Bigl(e^{-\sqrt{2}\,\xi\int_0^t B_x(s)\,\mathrm{d}s} f\bigl(B_x(t)\bigr) \bigm| B_x(s) \ge 0 \text{ on } [0,t]\Bigr)\,\mathbb{P}_x(\tau_0 > t)\,\mathrm{d}t\\ &\qquad = 2\frac{\phi_1(x)}{\phi_+(0)}\int_x^\infty \phi_+(y)f(y)\,\mathrm{d}y.\end{aligned} \tag{379}$$

Now divide by $x$ and let $x \searrow 0$. The distribution of $B_x$ conditioned on $B_x(s) \ge 0$ on $[0,t]$ converges to the distribution of the Brownian meander $B_{\mathrm{me}}^{(t)}$ on $[0,t]$. Furthermore, by the well-known formula for the distribution of the hitting time, with $\Phi$ the normal distribution function,

$$\frac{1}{x}\,\mathbb{P}_x(\tau_0>t)=\frac{1}{x}\Bigl(2\Phi\Bigl(\frac{x}{\sqrt{t}}\Bigr)-1\Bigr)\to\frac{2}{\sqrt{2\pi t}},\tag{380}$$

which also gives the bound

$$\frac{1}{x}\,\mathbb{P}_x(\tau_0>t)=\frac{1}{x}\Bigl(2\Phi\Bigl(\frac{x}{\sqrt{t}}\Bigr)-1\Bigr)\le\frac{2}{\sqrt{2\pi t}},\tag{381}$$

while the expectation inside the integral on the left hand side of (379) is bounded by $\sup|f|$, so dominated convergence applies to the left hand side of (379). For the right hand side we have $\phi_1(x)/x\to\phi_1'(0)=1$. Consequently, (379) yields

$$\int_0^\infty e^{-\lambda t}\,\mathbb{E}\,e^{-\sqrt{2}\,\xi\int_0^t B_{\mathrm{me}}^{(t)}(s)\,\mathrm{d}s} f\bigl(B_{\mathrm{me}}^{(t)}(t)\bigr)\frac{2}{\sqrt{2\pi t}}\,\mathrm{d}t \\ =\frac{2}{\phi_+(0)}\int_0^\infty\phi_+(y)f(y)\,\mathrm{d}y\tag{382}$$

for every $f$ that is bounded and continuous on $(0,\infty)$ and vanishes on some interval $(0,\delta)$, and thus by a limit argument (for example a monotone class argument) for every bounded $f$ on $(0,\infty)$. Taking $f(x)=1/\sqrt{2}$ we find, using (356) and the Brownian scaling $B_{\mathrm{me}}^{(t)}(ts)\overset{\mathrm{d}}{=}t^{1/2}B_{\mathrm{me}}(s)$, $0\le s\le 1$,

$$\begin{aligned}\int_0^\infty e^{-\lambda t}\psi_{\mathrm{me}}\bigl(\sqrt{2}\,\xi t^{3/2}\bigr)\frac{\mathrm{d}t}{\sqrt{\pi t}}&=\frac{\sqrt{2}}{\phi_+(0)}\int_0^\infty\phi_+(y)\,\mathrm{d}y\\&=\xi^{-1/3}\frac{\mathrm{AI}(\xi^{-2/3}\lambda)}{\mathrm{Ai}(\xi^{-2/3}\lambda)}.\end{aligned}\tag{383}$$

Taking $\xi=1$ yields (193) and $\lambda=1$ yields (310).

Note also that, since $B_{\mathrm{me}}^{(t)}(t)$ has the density function $(y/t)e^{-y^2/(2t)}$, $y>0$, it follows from (382) that for every $y>0$,

$$\int_0^\infty e^{-\lambda t}\,\mathbb{E}\Bigl(e^{-\sqrt{2}\,\xi\int_0^t B_{\mathrm{me}}^{(t)}(s)\,\mathrm{d}s}\Bigm| B_{\mathrm{me}}^{(t)}(t)=y\Bigr)\frac{y}{t}e^{-y^2/(2t)}\frac{1}{\sqrt{2\pi t}}\,\mathrm{d}t \\ =\frac{\phi_+(y)}{\phi_+(0)}=\frac{\mathrm{Ai}\bigl(\sqrt{2}\,\xi^{1/3}y+\xi^{-2/3}\lambda\bigr)}{\mathrm{Ai}\bigl(\xi^{-2/3}\lambda\bigr)},\tag{384}$$

which is a relation for the areas under the Bessel bridges $\bigl(B_{\mathrm{me}}^{(t)}\mid B_{\mathrm{me}}^{(t)}(t)=y\bigr)$.

### *C.4. Brownian excursion*

Formula (384) holds for every $\lambda>0$, and by monotone convergence as $\lambda\to0$ for $\lambda=0$ too. Taking the difference, we obtain

$$\int_0^\infty \bigl(1-e^{-\lambda t}\bigr)\,\mathbb{E}\Bigl(e^{-\sqrt{2}\,\xi\int_0^t B_{\mathrm{me}}^{(t)}(s)\,\mathrm{d}s} \Bigm| B_{\mathrm{me}}^{(t)}(t)=y\Bigr)\frac{y}{t}e^{-y^2/(2t)}\frac{1}{\sqrt{2\pi t}}\,\mathrm{d}t$$
$$= \frac{\mathrm{Ai}\bigl(\sqrt{2}\,\xi^{1/3}y\bigr)}{\mathrm{Ai}(0)} - \frac{\mathrm{Ai}\bigl(\sqrt{2}\,\xi^{1/3}y+\xi^{-2/3}\lambda\bigr)}{\mathrm{Ai}\bigl(\xi^{-2/3}\lambda\bigr)}. \tag{385}$$

Let $\xi = 1/\sqrt{2}$, divide by $y$ and let $y \to 0$, noting that as $y \to 0$, $\bigl(B_{\mathrm{me}}^{(t)} \mid B_{\mathrm{me}}^{(t)}(t)=y\bigr) \xrightarrow{\mathrm{d}} B_{\mathrm{ex}}^{(t)}$, a Brownian excursion on $[0,t]$. We then find by dominated convergence on the left hand side, and simple differentiation with respect to $y$ on the right hand side, the formula (81), which as explained in Section 13 yields (83) and (84).

Alternatively, we may take $\xi = 1$ in (384), differentiate with respect to $\lambda$ and then divide by $y$ and let $y \to 0$. This yields

$$\int_0^\infty e^{-\lambda t}\,\mathbb{E}\,e^{-\sqrt{2}\int_0^t B_{\mathrm{ex}}^{(t)}(s)\,\mathrm{d}s}\frac{1}{\sqrt{2\pi t}}\,\mathrm{d}t = -\frac{\partial}{\partial y}\Big|_{y=0}\frac{\partial}{\partial\lambda}\frac{\mathrm{Ai}\bigl(\sqrt{2}\,y+\lambda\bigr)}{\mathrm{Ai}(\lambda)}$$
$$= -\frac{\partial}{\partial\lambda}\frac{\partial}{\partial y}\Big|_{y=0}\frac{\mathrm{Ai}\bigl(\sqrt{2}\,y+\lambda\bigr)}{\mathrm{Ai}(\lambda)} = -\frac{\partial}{\partial\lambda}\frac{\sqrt{2}\,\mathrm{Ai}'(\lambda)}{\mathrm{Ai}(\lambda)}, \tag{386}$$

which yields (84).

### *C.5. Brownian double meander*

We do not know how to apply the Feynman–Kac formula directly to $B_{\mathrm{dm}}$, but we have shown in Section 27 how the result (314) follows from the result (310), or equivalently (383), for $\mathcal{B}_{\mathrm{me}}$.

## Acknowledgments

I thank Philippe Flajolet and Guy Louchard for providing me with several references.

## References

[1] M. Abramowitz & I. A. Stegun, eds., *Handbook of Mathematical Functions.* Dover, New York, 1972.

[2] D. Aldous, Brownian excursions, critical random graphs and the multiplicative coalescent. *Ann. Probab.* **25** (1997), 812–854. MR1434128

[3] B. L. Altshuler, A. G. Aronov & D. E. Khmelnitsky, Effects of electron-electron collisions with small energy transfers on quantum localisation. *J. Phys. C: Solid State. Phys.* **15** (1982) 7367–7386.

[4] G. N. Bagaev & E. F. Dmitriev, Enumeration of connected labeled bipartite graphs. (Russian.) *Doklady Akad. Nauk BSSR* **28** (1984), 1061–1063. MR0775755

[5] E. A. Bender, E. R. Canfield & B. D. McKay, Asymptotic properties of labeled connected graphs. *Random Struct. Alg.* **3** (1992), no. 2, 183–202. MR1151361

[6] J. Bertoin, J. Pitman & J. Ruiz de Chavez, Constructions of a Brownian path with a given minimum. *Electron. Comm. Probab.* **4** (1999), 31–37. MR1703609

[7] A. N. Borodin & P. Salminen, *Handbook of Brownian motion—facts and formulae.* 2nd edition. Birkhäuser, Basel, 2002. MR1912205

[8] P. Chassaing & G. Louchard, Reflected Brownian bridge area conditioned on its local time at the origin. *J. Algorithms* **44** (2002), no. 1, 29–51. MR1932676

[9] D. M. Cifarelli, Contributi intorno ad un test per l'omogeneità tra due campioni. *Giorn. Econom. Ann. Econom. (N.S.)* **34** (1975), no. 3–4, 233–249. MR0433704

[10] A. Comtet, J. Desbois & C. Texier, Functionals of Brownian motion, localization and metric graphs. *J. Phys. A* **38** (2005), no. 37, R341–R383. MR2169321

[11] M. Csörgő, Z. Shi & M. Yor, Some asymptotic properties of the local time of the uniform empirical process. *Bernoulli* **5** (1999), no. 6, 1035–1058. MR1735784

[12] D. A. Darling, On the supremum of a certain Gaussian process. *Ann. Probab.* **11** (1983), no. 3, 803–806. MR0704564

[13] I. V. Denisov, Random walk and the Wiener process considered from a maximum point. (Russian.) *Teor. Veroyatnost. i Primenen.* **28** (1983), no. 4, 785–788. English translation: *Theory Probab. Appl.* **28** (1983), no. 4, 821–824. MR0726906

[14] R. T. Durrett, D. L. Iglehart & D. R. Miller, Weak convergence to Brownian meander and Brownian excursion. *Ann. Probab.* **5** (1977), no. 1, 117–129. MR0436353

[15] W. Feller, *An Introduction to Probability Theory and its Applications,* Vol. II. 2nd ed., Wiley, New York, 1971. MR0270403

[16] J. A. Fill & S. Janson, Precise logarithmic asymptotics for the right tails of some limit random variables for random trees. Preprint, 2007. `http://www.arXiv.org/math.PR/0701259`

[17] P. Flajolet & G. Louchard, Analytic variations on the Airy distribution. *Algorithmica* **31** (2001), 361–377. MR1855255

[18] P. Flajolet, P. Poblete & A. Viola, On the analysis of linear probing hashing. *Algorithmica* **22** (1998), no. 4, 490–515. MR1701625

[19] P. Groeneboom, Brownian motion with a parabolic drift and Airy functions. *Probab. Theory Related Fields* **81** (1989), no. 1, 79–109. MR0981568

[20] S. Janson, *Gaussian Hilbert Spaces.* Cambridge Univ. Press, Cambridge, 1997. MR1474726

[21] S. Janson, The Wiener index of simply generated random trees. *Random Struct. Alg.* **22** (2003), no. 4, 337–358. MR1980963

[22] S. Janson, Left and right pathlenghts in random binary trees. *Algorithmica,* **46** (2006), no. 3/4, 419–429.

[23] S. Janson & P. Chassaing, The center of mass of the ISE and the Wiener index of trees. *Electronic Comm. Probab.* **9** (2004), paper 20, 178–187. MR2108865

[24] S. Janson, D. E. Knuth, T. Łuczak & B. Pittel, The birth of the giant component. *Random Struct. Alg.* **3** (1993), 233–358.

[25] S. Janson & G. Louchard, Tail estimates for the Brownian excursion area and other Brownian areas. In preparation.

[26] S. Janson & N. Petersson. In preparation.

[27] M. Jeanblanc, J. Pitman & M. Yor, The Feynman–Kac formula and decomposition of Brownian paths. *Mat. Apl. Comput.* **16** (1997), no. 1, 27–52. MR1458521

[28] B. McK. Johnson & T. Killeen, An explicit formula for the C.D.F. of the $L_1$ norm of the Brownian bridge. *Ann. Probab.* **11** (1983), no. 3, 807–808. MR0704570

[29] M. Kac, On the average of a certain Wiener functional and a related limit theorem in calculus of probability. *Trans. Amer. Math. Soc.* **59** (1946), 401–414. MR0016570

[30] M. Kac, On distributions of certain Wiener functionals. *Trans. Amer. Math. Soc.* **65** (1949), 1–13. MR0027960

[31] M. Kac, On some connections between probability theory and differential and integral equations. *Proceedings of the Second Berkeley Symposium on Mathematical Statistics and Probability, 1950*, University of California Press, Berkeley and Los Angeles, 1951, pp. 189–215. MR0045333

[32] N. N. Lebedev, *Special Functions and their Applications.* Dover, New York, 1972. (Translated from Russian.) MR0350075

[33] G. Louchard, Kac's formula, Lévy's local time and Brownian excursion. *J. Appl. Probab.* **21** (1984), no. 3, 479–499. MR0752014

[34] G. Louchard, The Brownian excursion area: a numerical analysis. *Comput. Math. Appl.* **10** (1984), no. 6, 413–417. Erratum: *Comput. Math. Appl. Part A* **12** (1986), no. 3, 375. MR0783514

[35] S. N. Majumdar & A. Comtet, Airy distribution function: from the area under a Brownian excursion to the maximal height of fluctuating interfaces. *J. Stat. Phys.* **119** (2005), no. 3-4, 777–826. MR2151223

[36] M. Nguyên Thê, Area of Brownian motion with generatingfunctionology. *Discrete random walks (Paris, 2003)*, Discrete Math. Theor. Comput. Sci. Proc., AC, Assoc. Discrete Math. Theor. Comput. Sci., Nancy, 2003, pp. 229–242, MR2042390

[37] M. Nguyen The, Area and inertial moment of Dyck paths. *Combin. Probab. Comput.* **13** (2004), no. 4–5, 697–716. MR2095979

[38] M. Perman & J. A. Wellner, On the distribution of Brownian areas. *Ann. Appl. Probab.* **6** (1996), no. 4, 1091–1111. MR1422979

[39] D. Revuz & M. Yor, *Continuous martingales and Brownian motion.* 3rd edition. Springer-Verlag, Berlin, 1999. MR1725357

[40] S. O. Rice, The integral of the absolute value of the pinned Wiener process—calculation of its probability density by numerical integration. *Ann. Probab.* **10** (1982), no. 1, 240–243. MR0637390

[41] C. Richard, On $q$-functional equations and excursion moments. Preprint, 2005. `http://www.arXiv.org/math.CO/0503198`

[42] L. A. Shepp, On the integral of the absolute value of the pinned Wiener process. *Ann. Probab.* **10** (1982), no. 1, 234–239. Acknowledgment of priority: *Ann. Probab.* **19** (1991), no. 3, 1397. MR0637389

[43] J. Spencer, Enumerating graphs and Brownian motion. *Comm. Pure Appl. Math.* **50** (1997), no. 3, 291–294. MR1431811

[44] L. Takács, A Bernoulli excursion and its various applications. *Adv. in Appl. Probab.* **23** (1991), no. 3, 557–585. MR1122875

[45] L. Takács, On a probability problem connected with railway traffic. *J. Appl. Math. Stochastic Anal.* **4** (1991), no. 1, 1–27. MR1095187

[46] L. Takács, Conditional limit theorems for branching processes. *J. Appl. Math. Stochastic Anal.* **4** (1991), no. 4, 263–292. MR1136428

[47] L. Takács, Random walk processes and their applications to order statistics. *Ann. Appl. Probab.* **2** (1992), no. 2, 435–459. MR1161061

[48] L. Takács, On the total heights of random rooted binary trees. *J. Combin. Theory Ser. B* **61** (1994), no. 2, 155–166. MR1280604

[49] L. Takács, On the distribution of the integral of the absolute value of the Brownian motion. *Ann. Appl. Probab.* **3** (1993), no. 1, 186–197. MR1202522

[50] L. Takács, Limit distributions for the Bernoulli meander. *J. Appl. Probab.* **32** (1995), no. 2, 375–395. MR1334893

[51] L. Tolmatz, Asymptotics of the distribution of the integral of the absolute value of the Brownian bridge for large arguments. *Ann. Probab.* **28** (2000), no. 1, 132–139. MR1756000

[52] L. Tolmatz, The saddle point method for the integral of the absolute value of the Brownian motion. *Discrete random walks (Paris, 2003)*, Discrete Math. Theor. Comput. Sci. Proc., AC, Nancy, 2003, pp. 309–324. MR2042397

[53] L. Tolmatz, Asymptotics of the distribution of the integral of the positive part of the Brownian bridge for large arguments. *J. Math. Anal. Appl.* **304** (2005), no. 2, 668–682. MR2126559

[54] W. Vervaat, A relation between Brownian bridge and Brownian excursion. *Ann. Probab.* **7** (1979), no. 1, 143–149. MR0515820

[55] V. A. Voblyĭ, O koeffitsientakh Raĭta i Stepanova-Raĭta. *Matematicheskie Zametki* **42** (1987), 854–862. English translation: Wright and Stepanov-Wright coefficients. *Mathematical Notes* **42** (1987), 969–974. MR0934817

[56] G. S. Watson, Goodness-of-fit tests on a circle. *Biometrika* **48** (1961), 109–114. MR0131930

[57] D. Williams, Path decomposition and continuity of local time for one-dimensional diffusions. I. *Proc. London Math. Soc. (3)* **28** (1974), 738–768. MR0350881

[58] E. M. Wright, The number of connected sparsely edged graphs. *J. Graph Th.* **1** (1977), 317–330. MR0463026

[59] E. M. Wright, The number of connected sparsely edged graphs. III. Asymptotic results. *J. Graph Th.* **4** (1980), 393–407. MR0595607